\documentclass[reqno,10.5pt]{amsart}

\usepackage[latin1]{inputenc}
\usepackage[english]{babel}

\usepackage{amscd}

\usepackage{amsmath}
\usepackage{amsthm}
\usepackage{amssymb}

\usepackage[notref, notcite, final]{showkeys}
\usepackage{color}
\usepackage{graphicx}
\usepackage[hidelinks]{hyperref}
\usepackage{cleveref}
\usepackage{fullpage}
\usepackage{comment}
\usepackage[shortlabels]{enumitem}

\usepackage{dsfont}

\newcommand{\id}{\mathcal{I}}

\newcommand{\rr}{\mathbb{R}}
\newcommand{\cc}{\mathbb{C}}
\newcommand{\hh}{\mathbb{H}}

\newcommand{\Q}{\mathcal{Q}}
\renewcommand{\Re}{\mathrm{Re}}

\newcommand{\boundOP}{\mathcal{B}}

\newcommand{\dom}{\operatorname{dom}}

\usepackage[autostyle]{csquotes}

\newtheorem{theorem}{Theorem}[section]
\newtheorem{lemma}[theorem]{Lemma}
\newtheorem{proposition}[theorem]{Proposition}
\newtheorem{corollary}[theorem]{Corollary}

\theoremstyle{definition}
\newtheorem{definition}[theorem]{Definition}

\newtheorem{assumption}[theorem]{Assumption}

\theoremstyle{remark}
\newtheorem{remark}[theorem]{Remark}
\crefname{enumi}{}{}

\title{\bf Poly slice monogenic functions, Cauchy formulas\\  and  the $PS$-functional calculus}

\author[D. Alpay]{Daniel Alpay}
\address{(DA)
Faculty of Mathematics, Physics, and Computation\\
Schmid College of Science and Technology\\
Chapman University\\
One University Drive
Orange, California 92866\\
USA}
\email{alpay@chapman.edu}

\author[F. Colombo]{Fabrizio Colombo}
\address{(FC)
Politecnico di Milano\\Dipartimento di Matematica\\Via E. Bonardi, 9\\20133
Milano, Italy}
\email{fabrizio.colombo@polimi.it}
\author[K. Diki]{Kamal  Diki}
\address{(KD)
Politecnico di Milano\\Dipartimento di Matematica\\Via E. Bonardi, 9\\20133
Milano, Italy
} \email{kamal.diki@polimi.it}
\author[I. Sabadini]{Irene Sabadini}
\address{(IS)
 Politecnico di Milano\\Dipartimento di Matematica\\Via E. Bonardi, 9\\20133
Milano, Italy
} \email{irene.sabadini@polimi.it}

\begin{document}
\maketitle

\begin{abstract}
Since 2006 the theory of slice hyperholomorphic functions
and the related spectral theory on the $S$-spectrum have had a very fast development.
This new spectral theory based on the $S$-spectrum has applications, for example, in the formulation of quaternionic quantum mechanics, in Schur analysis and  in fractional diffusion problems.
  In this paper we introduce and study the theory of poly slice monogenic functions, also proving some Cauchy type integral formulas. Then we introduce the associated functional calculus, called $PS$-functional calculus, which is the polyanalytic version of the $S$-functional calculus and which is based on the notion of $S$-spectrum. We study some different formulations of the calculus and we prove some of its properties, among which the product rules.
\end{abstract}

\medskip
\noindent AMS Classification 30G30, 47A10, 47A60.

\noindent Keywords: Poly slice monogenic function, Cauchy formulas, $PS$-resolvent operators, modified $S$-resolvent operators, spectral theory on the $S$-spectrum.

\noindent {\em }
\date{today}
\tableofcontents

\section{Introduction}

The theory of polyanalytic functions is an interesting topic in complex analysis. It extends the concept of holomorphic functions to nullsolutions of higher order powers of the Cauchy-Riemann operator.
Precisely, $n$-analytic or polyanalytic functions are nullsolutions of the $n$-power of the Cauchy-Riemann operator.
They were introduced in 1908 by Kolossov see
\cite{kolossov} to study elasticity problems. This  stream of
research was then continued by his
student Muskhelishvili and led to the book \cite{Muskhelishvili}. A
rather complete introduction to polyanalytic functions is in
\cite{Balk1991,balk_ency}. In more recent times this class of functions was
studied by various authors and with no pretense of completeness we
mention the works of Abreu, Agranovsky,  Begehr, Feichtinger, Vasilevski
\cite{abreu, abreufeicht, agranovsky, begehr, vasilevski} and the references
therein.
 Some famous Hilbert spaces of holomorphic functions that were extended to the setting of polyanalytic functions are the Bergman and Fock spaces, see for example \cite{AF2014,Alpay2015,Balk1991} and the references therein.

\medskip
 Polyanalytic functions are important not only from the theoretical point of view,  but also in the theory of signals since they allow to encode $n$ independent analytic functions into a single polyanalytic one using a special decomposition. This idea is similar to the problem of multiplexing signals. This is related to the construction of the polyanalytic Segal-Bargmann transform mapping $L^2(\mathbb{R})$ onto the poly-Fock space, see \cite{AF2014}.

\medskip
In quantum mechanics polyanalytic functions are relevant for the study of the Landau levels associated to Schr\"odinger operators, see \cite{AF2014,AIM1997}.
They were used also in \cite{A2010} to study sampling and interpolation problems on Fock spaces using time frequency analysis techniques such as short-time Fourier transform (STFT) or Gabor transforms. This allows to extend Bargmann theory to the polyanalytic setting using Gabor analysis.

\medskip
The theory of slice hyperholomorphic functions  started its full development from
the beginning of this century \cite{CSSf,GS2}. It has had a quite fast developments due to several authors and the main results, regarding the quaternionic setting, are contained in the books \cite{SCHURBOOK,Entirebook, QUATAPPBOOK,GSSBOOK} and the references therein,
 while for the Clifford algebra setting we refer the interested reader to the book \cite{NONCOMMBOOK} and its bibliography.
Nowadays the function theory has expanded in several directions
but it is in operators theory where it has found its most profound applications and several monographs have been published in the last decade
 \cite{COF,SCHURBOOK,FRACTBOOK,CGKBOOK,NONCOMMBOOK}.
The slice monogenic functions were introduced in \cite{CS,Cauchy,CSSf,CSSd,CSSe} also in collaboration with D. C. Struppa, and in this paper we generalize this class of functions to the poly analytic setting.

\medskip
In order to state our results we need to
explain the context in which we work  and to highlight the importance of this branch of operator theory which is called quaternionic and Clifford operator theory.
First  of all we point out that the appropriate definition of the quaternionic spectrum for a quaternionic linear operator has been open problem at least since the
paper \cite{BF}  of  G. Birkhoff, J. von Neumann, in 1936  on the logic of quantum mechanics, where the authors
 showed that quantum mechanics can be formulated also on quaternions.
Moreover, consider  a generalization of the gradient operators such as
$$
T=ia(x_1,x_2,x_3)\partial_{x_1}+jb(x_1,x_2,x_3)\partial_{x_2}+kc(x_1,x_2,x_3)\partial_{x_3}
$$
where $a$, $b$ and $c$ are given real valued functions of the variables $(x_1,x_2,x_3)\in \mathbb{R}^3$,
and $i$,$j$,$k$ are the imaginary units of the quaternions.
It is very interesting to observe that the spectral theory
for vector operators like the gradient operator $\nabla$, or its generalizations such as the operator $T$ defined above,
 has been unclear since long time, even before 1936.

\medskip
Regarding the quaternionic spectral theorem we observe that
some attempts have been done after the paper of  G. Birkhoff, J. von Neumann, but
all the approaches suffered of the lack of an appropriate notion
 of quaternionic spectrum. The turning point came in
 2006 when it was introduced the $S$-spectrum and the $S$-functional calculus which are a cornerstone of quaternionic and Clifford operator theory.
The $S$-spectrum was identified by purely hypercomplex analysis
techniques and not on physical considerations, as it is widely explained in the introduction of \cite{CGKBOOK}.

\medskip
   The spectral theorem for quaternionic normal operators based on the $S$-spectrum was finally proved in 2015 by Alpay, Colombo and Kimsey and published in 2016, see \cite{6SpecThm1}. This theorem is the most important tool for the formulation of quaternionic quantum mechanics and more recently there have been
several efforts to study the perturbations of quaternionic normal operators
in \cite{Paula}, moreover, the theory of quaternionic spectral operators has been developed in \cite{JONAMEM}. The theory of characteristic operator function
has started its development in this setting not too long ago and the main
 advances can be found in the book \cite{COF}.

\medskip
There are several applications of the spectral theory on the $S$-spectrum to
   fractional diffusion and fractional evolution problems because it is possible to define the fractional powers of vector operators so that we can generate fractional Fourier laws, see \cite{FRACTBOOK}.
   With this strategy we are able to write the  fractional heat  equation modifying just the Fourier's law and preserving the conservation of energy law.

\medskip
Among the most developed areas in the slice hyperholomorphic setting
there is the theory of slice hyperholomorphic reproducing kernel Hilbert spaces
and more in general  quaternionic Schur analysis has been largely investigated in the last decade.
The material is spread over several papers, but the interested reader
can find several results in the book \cite{SCHURBOOK} and in the references therein.

  \medskip
  The importance of hypercomplex analysis in operators theory
  is not limited just to the slice hypercomplex setting.
  In fact, using the classical theory
    of monogenic functions, see \cite{BLU,DELSOSOU,Gurlebeck:2008},
       A. Mc Intosh and his collaborators \cite{jmcpw,mcp} developed
    the monogenic functional calculus, which is based on the notion of monogenic spectrum.
     This monogenic calculus contains as particular cases the Weyl functional calculus
    and the Taylor functional calculus for commuting operators, see the book \cite{jefferies}.
   In harmonic analysis, the monogenic functional calculus plays a crucial role
     as shown in the book \cite{taobook} and the references therein.
     Moreover, the monogenic function theory has applications in boundary value problems \cite{BOOKEBVP}, and in Clifford wavelets, singular integrals, and Hardy spaces as one can see in the book \cite{mitreabook}.

    \medskip To complete this short introduction on the hypercomplex spectral theories we recall that
    there is a link between the slice hyperholomorphic functions and the monogenic functions via the Fueter-Sce-Qian extension theorem and that there is a link between the two spectral theories.
     This link is
     the so-called $F$-functional calculus that generates a version of the monogenic functional
      calculus using the notion of $S$-spectrum, see \cite{SCEBOOK}.
      We recall that also
      the Radon transform is a bridge between the monogenic and generalized slice monogenic functions, see \cite{radon}.

\medskip
This introductory part explains how the results of this paper have to be seen in the framework of hyperholomorphic function theories and the associated spectral theories.

In fact, here we extend the slice monogenic function theory and its $S$-functional calculus to the
poly slice monogenic setting, namely to the set of (suitable) functions in the kernel of the $M$th-power of the Cauchy-Riemann operator. The quaternionic counterpart of the function theory
  started with the recent works \cite{Rsimmath,ADSP2020}, while the corresponding functional calculus
    is introduced in this paper for the first time. In fact, in this paper we begin a systematic study of the function theory, also proving the Cauchy formulas. These formulas can be written using different Cauchy kernels which extend the one in the complex case. However, the noncommutative context requires suitable techniques in order to prove the results. Furthermore, the components of the kernels in their poly slice monogenic decomposition have different behavior at infinity.

Then, we define the so-called $PS$-functional calculus. It
is the poly slice monogenic version of the $S$-functional calculus and it coincides with it when the order is $1$.
This calculus is based on the $S$-spectrum, see Definition \ref{sspectrum}, and it applies
to  $(n+1)$-tuples of noncommuting operators $(T_0,T_1,...,T_n)$ written as the
paravector operator $T=T_0+T_1e_1+...+T_ne_n$, where $e_1,...,e_n$ are the units of the Clifford algebra $\mathbb{R}_n$. The quaternionic case is obviously a particular case.
We prove several results and a crucial tool is given by suitable modified $S$-resolvent operators for which we could also prove the resolvent equations. For some results, like the product rules we assumed commutativity of the components of the operators.

\medskip
The contents of the paper are organized as follows.
In Section \ref{preliminary} we recall some preliminary results on the theory of
slice monogenic functions.
In Section \ref{slicepoly} we develop the theory of poly slice monogenic functions and we show some properties. In particular, we prove a slice monogenic integral representation of poly slice monogenic functions that will be used for a representation of the $PS$-functional calculus.
In Section \ref{sec4} we prove the Cauchy formulas and we define a product of poly slice monogenic functions.
 In Section \ref{PSSEC}
 we give the formulations of the $PS$-functional calculus via the $PS$-resolvent operators and the poly slice monogenic Cauchy formula.
 In Section \ref{modofresr}
 we define and study the formulations of the $PS$-functional calculus
 via some suitably modified $S$-resolvent operators.
 In Section \ref{equivprodrule} we show the
 equivalence of the two definitions of the $PS$-functional calculus
  and we prove the product rules.

\section{Preliminary material}\label{preliminary}

In this section we recall some preliminary material useful
 to extend the theory of polyanalytic functions to the slice monogenic setting.
The classical polyanalytic functions are those functions $f:\Omega\subseteq \mathbb{C} \to \mathbb{C}$ of class $\mathcal{C}^M(\Omega)$, for $M\in \mathbb{N}$, such that
 $$
 \overline{\partial_i}^M f(z)=0,\ \ \ \ {\rm for \ all} \ \ \ \  z=u+iv\in \Omega, \ \ \ \ \ i=\sqrt{-1}
 $$
 where
 \begin{equation}\label{compoly}
 \overline{\partial_i}^M =\frac{1}{2^M}(\partial_u+i\partial_v)^M.
\end{equation}

A Cauchy-type formula for polyanalytic functions appeared for the first
time in Th\'eodoresco's doctoral thesis \cite{teodoresco} and recalled
in the paper \cite{balk_ency}:
\begin{theorem} [\cite{balk_ency}, Theorem 1.3]
If a function $f$ is  polyanalytic of order $M$ in a closed domain $G$
bounded by
a rectifiable closed contour $\Gamma$, then the value of $f$ at any
point $z$ of the domain
$G$ is expressed, using values of the function itself and its formal
derivatives at
points $t$ of the boundary $\Gamma$, by the formula
\begin{equation}\label{TEODOR}
f(z)=\frac{1}{2\pi i}\sum_{\ell=0}^{M-1}\int_{\Gamma} \frac{1}{\ell! (t-z)}
(\bar z-\bar t)^\ell \frac{\partial^\ell f}{\partial \bar t^\ell}\, dt.
\end{equation}
\end{theorem}
The formula contains a finite sum in which appear the kernels
$$
\pi_\ell(z,t)= \frac{1}{\ell! (t-z)} (\bar z-\bar t)^\ell,\ \ \ \ell=0,...,M-1.
$$
Another polyanalytic Cauchy formula is given by
$$
f(z)=\int_{\partial \Omega} \sum_{\ell=0}^{M-1}(-2)^\ell L_\ell(w-z)\ dw\
\overline{\partial_i}^\ell f(w),
$$
where $\partial \Omega$ is the boundary of the smooth bounded domain $\Omega$ in $\mathbb{C}$,
 the infinitesimal arc length is denoted by $dw$, and we have set
$$
L_\ell(z):=\frac{1}{2\pi i \  z}\frac{({\rm Re}(z))^\ell}{\ell !}, \ \ \ \ell=0,...,M-1.
$$
We observe that when we define $\overline{\partial_i}^M $ in (\ref{compoly})
 without the prefactor  $1/2^M$, then in the Cauchy formula the coefficients $(-2)^\ell$ have to be replaced by $(-1)^\ell$.

\medskip
In the quaternionic setting or, more in general, in the Clifford algebra setting, one can extend the notion of holomorphic functions by considering functions in the kernel of a generalized Cauchy-Riemann operator, thus obtaining the so-called regular or monogenic functions, or of its $n$-power, thus obtaining poly-regular functions or poly-monogenic functions, see
\cite{brackx1976k, delanghe1978hypercomplex}.

\medskip
We now recall the main definitions on Clifford algebras and the main facts on slice monogenic functions
that are necessary to introduce and develop the theory of poly slice monogenic functions.
Let $\rr_n$ be the real Clifford algebra over $n$ imaginary units $e_1,\ldots ,e_n$
satisfying the relations $e_\ell e_m+e_me_\ell=0$,\  $\ell\not= m$, $e_\ell^2=-1.$
 An element in the Clifford algebra will be denoted by $\sum_A e_Ax_A$ where
$A=\{ \ell_1\ldots \ell_r\}\in \mathcal{P}\{1,2,\ldots, n\},\ \  \ell_1<\ldots <\ell_r$
 is a multi-index
and $e_A=e_{\ell_1} e_{\ell_2}\ldots e_{\ell_r}$, $e_\emptyset =1$.
An element $(x_0,x_1,\ldots,x_n)\in \mathbb{R}^{n+1}$  will be identified with the element
$
 x=x_0+\underline{x}=x_0+ \sum_{\ell=1}^nx_\ell e_\ell\in\mathbb{R}_n
$
called paravector and the real part $x_0$ of $x$ will also be denoted by $\Re(x)$.
The norm of $x\in\mathbb{R}^{n+1}$ is defined as $|x|^2=x_0^2+x_1^2+\ldots +x_n^2$.
 The conjugate of $x$ is defined by
$
\bar x=x_0-\underline x=x_0- \sum_{\ell=1}^nx_\ell e_\ell.
$
We denote by $\mathbb{S}$ the sphere
$$
\mathbb{S}=\{ \underline{x}=e_1x_1+\ldots +e_nx_n\ | \  x_1^2+\ldots +x_n^2=1\};
$$
for $\mathrm{j}\in\mathbb{S}$ we obviously have $\mathrm{j}^2=-1$.
Given an element $x=x_0+\underline{x}\in\rr^{n+1}$ we put
$
\mathrm{j}_x=\underline{x}/|\underline{x}|$ if $\underline{x}\not=0,
$
 and given an element $x\in\rr^{n+1}$. The set
$$
[x]:=\{y\in\rr^{n+1}\ :\ y=x_0+\mathrm{j} |\underline{x}|, \ \mathrm{j}\in \mathbb{S}\}
$$
is an $(n-1)$-dimensional sphere in $\mathbb{R}^{n+1}$.
The vector space $\mathbb{R}+\mathrm{j}\mathbb{R}$ passing through $1$ and
$\mathrm{j}\in \mathbb{S}$ will be denoted by $\mathbb{C}_\mathrm{j}$ and
an element belonging to $\mathbb{C}_\mathrm{j}$ will be indicated by $u+\mathrm{j}v$, for $u$, $v\in \mathbb{R}$.
With an abuse of notation we will write $x\in\mathbb{R}^{n+1}$.
Thus, if $U\subseteq\mathbb{R}^{n+1}$ is an open set,
a function $f:\ U\subseteq \mathbb{R}^{n+1}\to\mathbb{R}_n$ can be interpreted as
a function of the paravector $x$.

\medskip
In this paper we use the definition of poly slice monogenic functions that is the generalization of
slice monogenic functions in the spirit of the Fueter-Sce-Qian mapping theorem, see \cite{SCEBOOK}.
This definition is the most appropriate for operator theory
and the reason is widely explained in several papers and in the books \cite{FRACTBOOK,CGKBOOK}. The corresponding function theory in real alternative algebras is developed in \cite{GP2}. The same definition will be used for
poly slice monogenic functions in the next section.

\begin{definition}
Let $U\subseteq \mathbb{R}^{n+1}$.
 We say that $U$ is axially symmetric if $[x]\in U$ for every $x\in U$.
\end{definition}

\begin{definition}[Slice monogenic functions]\label{SHolDef}
 Let $U\subseteq\mathbb{R}^{n+1}$ be an axially symmetric open set and let $\mathcal{U} = \{ (u,v)\in\rr^2: u+ \mathbb{S} v\subset U\}$. A function $f:U\to \mathbb{R}_n$ is called a left
 slice function, if it is of the form
 \[
 f(x) = f_{0}(u,v) + \mathrm{j}f_{1}(u,v)\qquad \text{for } x = u + \mathrm{j} v\in U
 \]
with the two functions $f_{0},f_{1}: \mathcal{U}\to \mathbb{R}_n$ that satisfy the compatibility conditions
\begin{equation}\label{CCondslic}
f_{0}(u,-v) = f_{0}(u,v),\qquad f_{1}(u,-v) = -f_{1}(u,v).
\end{equation}
If in addition $f_{0}$ and $f_{1}$ are $\mathcal C^1$ and satisfy the Cauchy-Riemann equations
\begin{equation}\label{CRR}
\begin{split}
&\partial_u f_{0}(u,v) -\partial_vf_{1}(u,v))=0
\\
&
\partial_v f_{0}(u,v) +\partial_uf_{1}(u,v))=0
\end{split}
\end{equation}
 then $f$ is called left slice monogenic.
A function $f:U\to \mathbb{R}_n$ is called a right slice function if it is of the form
\[
f(x) = f_{0}(u,v) + f_{1}(u,v) \mathrm{j}\qquad \text{for } x = u+ \mathrm{j}v \in U
\]
with the two functions $f_{0},f_{1}: \mathcal{U}\to \mathbb{R}_n$ that satisfy (\ref{CCondslic}).
If $f_{0}$ and $f_{1}$ are $\mathcal C^1$ and satisfy the Cauchy-Riemann equations
(\ref{CRR})
 then $f$ is called right slice monogenic.

If $f$ is a left (or right) slice function such that $f_{0}$ and $f_{1}$ are real-valued, then $f$ is called intrinsic.

We denote the sets of left, right and intrinsic
 slice monogenic functions on $U$ by $\mathcal{SM}_L(U)$,
$\mathcal{SM}_R(U)$ and $\mathcal{N}(U)$, respectively.

\end{definition}

\begin{remark}
The set $\mathcal{N}(U)$ is contained in both $\mathcal{SM}_L(U)$ and $\mathcal{SM}_R(U)$.
\end{remark}
\begin{definition}
We define the notion of $\mathrm{j}$-derivative by means of the
operator:
$$
{\partial}_\mathrm{j}:=\frac{1}{2}\left(\partial_u-\mathrm{j}\partial_v\right).
$$
For consistency, we will
denote by
$$
\overline{\partial}_\mathrm{j}=\frac{1}{2}\left(\partial_u+\mathrm{j}\partial_v\right)
$$
the Cauchy- Riemann operator associated with the complex plane $\mathbb{C}_\mathrm{j}$, for $\mathrm{j}\in \mathbb{S}$.
\end{definition}
Using the notations we have just introduced, the condition of left slice monogenicity will be expressed, in short, by
$
\overline{\partial}_\mathrm{j} f=0.
$
Right slice monogenicity will be expressed, with an abuse of notation, by
$
f\overline{\partial}_\mathrm{j} =0.
$
\begin{definition}
Let $U$ be an open set in $\mathbb{R}^{n+1}$ and let $f:\
U\to\mathbb{R}_n$ be a slice monogenic function. Its slice derivative
$\partial_S$  is defined as
\begin{equation}\label{s-derivative}
\partial_S(f)=\left\{\begin{array}{ll}
\partial_\mathrm{j}(f)(x)\quad &x=u+\mathrm{j}v,\ v\not=0\\
\partial_uf (u)\quad &{u\in\mathbb{R}}.
 \end{array}\right.
\end{equation}
\end{definition}

\begin{lemma}[Splitting Lemma]\label{splittinglemmamon}
\index{Splitting Lemma}
Let $U\subseteq\mathbb{R}^{n+1}$  be an axially symmetric open set and let $f\in \mathcal{SM}_L(U)$. For every
$\mathrm{j}=\mathrm{j}_1\in\mathbb{S}$ let $\mathrm{j}_2,\ldots, \mathrm{j}_n$ be a completion to a
basis of  $\mathbb{R}_n$ satisfying the defining
relations $\mathrm{j}_r\mathrm{j}_s+\mathrm{j}_s\mathrm{j}_r=-2\delta_{rs}$. Then there exist $2^{n-1}$
holomorphic functions $\mathcal{F}_A\ : U\cap \mathbb{C}_\mathrm{j}\to \mathbb{C}_\mathrm{j}$ such that for
every $z=u+\mathrm{j}v$ we have
$$
f_\mathrm{j}(z)=\sum_{|A|=0}^{n-1} \mathcal{F}_A(z) \mathrm{j}_A,\quad \mathrm{j}_A=\mathrm{j}_{i_1}\ldots \mathrm{j}_{i_s},
$$
where $A=i_1\ldots i_s$ is a subset of $\{2,\ldots ,n\}$, with ${i_1}<\ldots <{i_s}$,
or, when $|A|=0$, $\mathrm{j}_\emptyset=1$.
When $f\in \mathcal{SM}_R(U)$, then the splitting lemma becomes
$$
f_\mathrm{j}(z)=\sum_{|A|=0}^{n-1} \mathrm{j}_A\mathcal{F}_A(z) ,\quad \mathrm{j}_A=\mathrm{j}_{i_1}\ldots \mathrm{j}_{i_s}.
$$

\end{lemma}

The following formula is an immediate consequence of the definition of slice functions, see \cite{GP2}.
\begin{theorem}[The Structure Formula or Representation Formula]\label{formulamon}
\index{Representation Formula}
Let $U\subseteq
\mathbb{R}^{n+1}$ be an axially symmetric open set.

(I) Let
$f\in \mathcal{SM}_L(U)$.
Then, for any  vector
$x =u+\mathrm{j}_x v\in U$, the following formula holds:
\begin{equation}\label{repL}
f(x)=\frac{1}{2}\Big[1 - \mathrm{j} \mathrm{j}_x \Big]f(u +\mathrm{j} v)+\frac{1}{2}\Big[1 +  \mathrm{j}\mathrm{j}_x \Big]f(u -\mathrm{j} v),  \ \ {\rm for \ all} \ \ u +\mathrm{j} v\in U, \ \ \mathrm{j}\in \mathbb{S}.
\end{equation}

(II) Let $f\in \mathcal{SM}_R(U)$.
Then, for any  vector
$x =u+\mathrm{j}_x v\in U$, the following formula holds:
\begin{equation}\label{repR}
f(x)=\frac{1}{2}f(u +\mathrm{j} v)\Big[1 -  \mathrm{j}\mathrm{j}_x \Big]+\frac{1}{2}f(u -\mathrm{j} v)\Big[1 +  \mathrm{j}\mathrm{j}_x \Big],  \ \ {\rm for \ all} \ \ u +\mathrm{j} v\in U, \ \ \mathrm{j}\in \mathbb{S}.
\end{equation}
\end{theorem}
\begin{remark}
Using the representation formula we can write the slice monogenic Cauchy kernels in terms of the complex Cauchy kernel. For example for $S^{-1}_L(s,x)$ we have
$$
S^{-1}_L(s,x)=\frac{1}{2}\Big[1 -  \mathrm{j}_x \mathrm{j}\Big]\frac{1}{s-z}+\frac{1}{2}\Big[1 +  \mathrm{j}_x \mathrm{j}\Big]\frac{1}{s-\overline{z}}
$$
where we set $x=u+\mathrm{j}_x v$,  $z=u+\mathrm{j}v$, $s=s_0+\mathrm{j}s_1$, and where $\mathrm{j}$ is the imaginary unit of the complex plane $\mathbb{C}_\mathrm{j}$.
\end{remark}
For slice monogenic functions we have two equivalent ways
to write the Cauchy kernels.

\begin{proposition}\label{secondAA}
If $x, s\in \mathbb{R}^{n+1}$ with $x\not\in [s]$, then
\begin{gather}\label{secondAAEQ}
-(x^2 -2x \Re  (s)+|s|^2)^{-1}(x-\overline s)=(s-\bar q)(s^2-2\Re (x)s+|x|^2)^{-1}
\end{gather}
and
\begin{gather}\label{secondAAEQ1}
 (s^2-2\Re (x)s+|x|^2)^{-1}(s-\bar x)=-(x-\bar s)(x^2-2\Re (s)x+|s|^2)^{-1} .
\end{gather}
\end{proposition}
So we can give the following definition to distinguish the two representations of the Cauchy kernels.
\begin{definition}
Let $x,s\in \mathbb{R}^{n+1}$ with $x\not\in [s]$.
\begin{itemize}
\item
We say that  $S_L^{-1}(s,x)$ is written in the form I if
$$
S_L^{-1}(s,x):=-(x^2 -2 \Re  (s) x+|s|^2)^{-1}(x-\overline s).
$$
\item
We say that $S_L^{-1}(s,x)$ is written in the form II if
$$
S_L^{-1}(s,x):=(s-\bar x)(s^2-2\Re (x) s+|x|^2)^{-1}.
$$
\item
We say that  $S_R^{-1}(s,x)$ is written in the form I if
$$
S_R^{-1}(s,x):=-(x-\bar s)(x^2-2\Re (s)x+|s|^2)^{-1} .
$$
\item
We say that $S_R^{-1}(s,x)$ is written in the form II if
$$
S_R^{-1}(s,x):=(s^2-2\Re (x)s+|x|^2)^{-1}(s-\bar x).
$$
\end{itemize}
\end{definition}
\begin{lemma} Let $x,s\in \mathbb{R}^{n+1}$ with $s\notin [x]$.
The left slice monogenic Cauchy kernel $S_L^{-1}(s,x)$ is left slice monogenic in $x$ and right slice monogenic in $s$.
The right slice monogenic Cauchy kernel $S_R^{-1}(s,x)$ is left slice monogenic in $s$ and right slice monogenic in $x$.
\end{lemma}

\begin{definition}[Slice Cauchy domain]\label{Slice Cauchy domain}
An axially symmetric open set $U\subset \mathbb{R}^{n+1}$ is called a slice Cauchy domain, if $U\cap\cc_\mathrm{j}$ is a Cauchy domain in $\cc_\mathrm{j}$ for any $\mathrm{j}\in\mathbb{S}$. More precisely, $U$ is a slice Cauchy domain if, for any $\mathrm{j}\in\mathbb{S}$, the boundary ${\partial( U\cap\cc_\mathrm{j})}$ of $U\cap\cc_\mathrm{j}$ is the union a finite number of non-intersecting piecewise continuously differentiable Jordan curves in $\cc_{\mathrm{j}}$.
\end{definition}
\begin{theorem}[Cauchy formulas]\label{Cauchy formulas}
\label{Cauchygenerale}
Let $U\subset\mathbb{R}^{n+1}$ be a slice Cauchy domain, let $\mathrm{j}\in\mathbb{S}$ and set  $ds_\mathrm{j}=ds (-\mathrm{j})$.
If $f$ is a (left) slice monogenic function on a set that contains $\overline{U}$ then
\begin{equation}\label{cauchynuovo}
 f(x)=\frac{1}{2 \pi}\int_{\partial (U\cap \mathbb{C}_\mathrm{j})} S_L^{-1}(s,x)\, ds_\mathrm{j}\,  f(s),\qquad\text{for any }\ \  x\in U.
\end{equation}
If $f$ is a right slice monogenic function on a set that contains $\overline{U}$,
then
\begin{equation}\label{Cauchyright}
 f(x)=\frac{1}{2 \pi}\int_{\partial (U\cap \mathbb{C}_\mathrm{j})}  f(s)\, ds_\mathrm{j}\, S_R^{-1}(s,x),\qquad\text{for any }\ \  x\in U.
 \end{equation}
These integrals  depend neither on $U$ nor on the imaginary unit $\mathrm{j}\in\mathbb{S}$.
\end{theorem}

\begin{theorem}[Cauchy formulas on unbounded slice Cauchy domains]
Let $U\subset \mathbb{R}^{n+1}$ be an unbounded slice Cauchy domain and let $\mathrm{j}\in\mathbb{S}$. If $f\in\mathcal{SM}_{L}(\overline{U})$ and $f(\infty) := \lim_{|x|\to\infty}f(x)$ exists and is finite, then
\[
f(x) = f(\infty) + \frac{1}{2\pi}\int_{\partial(U\cap\mathbb{C}_{\mathrm{j}})}S_L^{-1}(s,x)\,ds_{\mathrm{j}}\, f(s) \qquad\text{for any }\ \ x\in U.
\]
If $f\in\mathcal{SM}_{R}(\overline{U})$ and $f(\infty) := \lim_{|x|\to\infty}f(x)$ exists and is finite, then
\[
f(x) = f(\infty) + \frac{1}{2\pi}\int_{\partial(U\cap\mathbb{C}_j)}f(s)\,ds_j\,S_R^{-1}(s,x) \qquad\text{for any }\ \ x\in U.
\]
 \end{theorem}

\section{Poly slice monogenic functions}\label{slicepoly}

We can now give the definition of poly slice monogenic functions and develop the corresponding function theory, which will be used to define the $PS$-functional calculus.
This definition is the monogenic counterpart of Definition 3.17 in  \cite{Rsimmath}.

\begin{definition}[Poly slice  monogenic functions]\label{SHolDef}

Let $M\in \mathbb{N}$ and denote by $\mathcal{C}^M(U)$ the set of
continuously differentiable functions with all their derivatives up to order $M$ on an
axially symmetric open set $U\subseteq\mathbb{R}^{n+1}$.
We let $\mathcal{U} = \{ (u,v)\in\rr^2: u+ \mathbb{S} v\subset U\}$.
A function $F:U\to \mathbb{R}_n$
 is called a left
 slice function, if it is of the form
 \[
 F(x) = F_{0}(u,v) + \mathrm{j}F_{1}(u,v)\qquad \text{for } x = u + \mathrm{j} v\in U
 \]
with the two functions $F_{0},F_{1}: \mathcal{U}\to \mathbb{R}_n$ that satisfy the compatibility condition
\begin{equation}\label{CCond}
F_{0}(u,-v) = F_{0}(u,v),\qquad F_{1}(u,-v) = -F_{1}(u,v).
\end{equation}
If in addition $F_{0}$ and $F_{1}$ are in $\mathcal{C}^M(U)$ and satisfy the poly Cauchy-Riemann equations of order $M\in \mathbb{N}$
 \begin{align}\label{CR}
\frac{1}{2^M}(\partial_u+\mathrm{j}\partial_v)^M(F_{0}(u,v) + \mathrm{j}F_{1}(u,v))=0,\ \ \ {\rm for\ all} \ \ \mathrm{j}\in \mathbb{S}
\end{align}
 then $F$ is called left poly slice monogenic function  of order $M\in \mathbb{N}$.
A function $F:U\to \mathbb{R}_n$ is called a right slice function if it is of the form
\[
F(x) = F_{0}(u,v) + F_{1}(u,v) \mathrm{j}\qquad \text{for } x = u+ \mathrm{j}v \in U
\]
with two functions $F_{0},F_{1}: \mathcal{U}\to \mathbb{R}^{n+1}$ that satisfy \eqref{CCond}.
If in addition $F_{0}$ and $F_{1}$ are in $\mathcal{C}^M(U)$ and satisfy the poly Cauchy-Riemann equations of order $M\in \mathbb{N}$
\begin{align}\label{CRRIGHT}
(F_{0}(u,v) + F_{1}(u,v)\mathrm{j})\frac{1}{2^M}(\partial_u+\mathrm{j}\partial_v)^M=0,\ \ \ {\rm for\ all} \ \ \mathrm{j}\in \mathbb{S}
\end{align}
 then $F$ is called right poly slice monogenic  of order $M\in \mathbb{N}$.
We will denote by $\mathcal{PS}^M_L(U)$ and $\mathcal{PS}^M_R(U)$ the set of left and right poly slice monogenic functions on the open set $U$, respectively.

If $F$ is a left (or right) slice function such that $F_{0}$ and $F_{1}$ are real-valued, then $F$ is called intrinsic.
By $\mathcal{PN}^M(U)$ we denote the set of poly intrinsic slice monogenic functions.
\end{definition}
\begin{remark}
For the definition of slice monogenic functions we required that the pair $(f_0,f_1)$ satisfies the Cauchy-Riemann system. In the case of poly slice monogenic functions we follow the same idea by observing that (for functions $F$ of class $\mathcal{C}^M(U)$)
$$
\frac{1}{2^M}(\partial_u+\mathrm{j}\partial_v)^MF
=\sum_{k=0}^{M}{M\choose k}\left(\partial_u\right)^{M-k}\left(\partial_v\right)^k
\mathrm{j}^kF
=\sum_{k=0}^{M}{M\choose k}D_{M-k,k}\mathrm{j}^kF
$$
where $D_{M-k,k}=\left(\partial_u\right)^{M-k}\left(\partial_v\right)^k$.
Due to the arbitrarity of $\mathrm{j}$, the condition that $F$ is left poly slice monogenic translates into a
system of two differential equations of order $M$ for the pair $(F_0,F_1)$ of $\mathbb{R}_n$-valued functions
(which reduces to the Cauchy-Riemann system
when $M=1$):
\begin{equation}\label{polyCR}
\begin{split}
&\left(\sum^M_{k=0 (mod\, 4)}{M\choose k} D_{M-k,k} -\sum^M_{k=2 (mod\,
4)}{M\choose k} D_{M-k,k}\right)F_0(u,v)
\\
&+\left(-\sum^M_{k=1 (mod\, 4)}{M\choose k} D_{M-k,k}+\sum^M_{k=3 (mod\,
4)}{M\choose k} D_{M-k,k}\right) F_1(u,v)=0,
\end{split}
\end{equation}
\begin{equation*}
\begin{split}
&\left(\sum^M_{k=1 (mod\, 4)}{M\choose k}  D_{M-k,k} -\sum^M_{k=3 (mod\,
4)}{M\choose k} D_{M-k,k} \right)F_0(u,v)
\\
&+\left(\sum^M_{k=0 (mod\, 4)}{M\choose k} D_{M-k,k}-\sum^M_{k=2 (mod\,
4)}{M\choose k} D_{M-k,k}\right) F_1(u,v)=0.
\end{split}
\end{equation*}
With obvious meaning of the symbols, we have
\[
\begin{split}
   &\mathsf{D}_1 F_0(u,v)-\mathsf{D}_2 F_1(u,v)=0\\
   &\mathsf{D}_2 F_0(u,v)+\mathsf{D}_1 F_1(u,v)=0
   \end{split}
\]
However, since this system is rather complicated to write, we prefer to use the above Definition \ref{SHolDef}.
\end{remark}

We point out that the definition of
poly slice monogenic functions extends to functions with values in a Clifford Banach module in a very natural way.
As in the quaternionic case, one can give the very useful definition of (strong) slice monogenicity,
see \cite{SCHURBOOK} for vector-valued functions.

\begin{definition}[Poly slice monogenic functions vector-valued]\label{SHolStrong}
 Let $U\subseteq \mathbb{R}^{n+1}$ be an axially symmetric open set and let
 \[
 \mathcal{U} = \{ (u,v)\in \mathbb{R}^2: u+ \mathbb{S} v\subset U\}.
 \]
  A function $f:U\to X_L$ with values in a  left Clifford Banach module $X_L$ is called
  a left slice function, if is of the form
 \[
 F(x) = F_{0}(u,v) + \mathrm{j}F_{1}(u,v)\qquad \text{for } x = u + j v\in U
 \]
with two functions $F_{0},F_{1}: \mathcal{U} \to X_L$ that satisfy the compatibility condition \eqref{CCond}. If in addition $F_{0}$ and $F_{1}$ are in $\mathcal{C}^M(U)$ and satisfy the poly Cauchy-Riemann equations \eqref{CR}, then $F$ is called {\em (strongly) left poly slice monogenic}.

A function $f:U\to X_R$ with values in a  right Clifford-Banach module
 is called a right slice function if it is of the form
\[
F(x) = F_{0}(u,v) + F_{1}(u,v) j\qquad \text{for } x = u+ jv \in U
\]
with two functions $F_{0},F_{1}: \mathcal{U} \to X_R$ that satisfy the compatibility condition \eqref{CCond}. If in addition $F_{0}$ and $F_{1}$ are in $\mathcal{C}^M(U)$ and satisfy the poly Cauchy-Riemann equations \eqref{CRRIGHT}, then $f$ is called {\em (strongly) right poly slice monogenic}.
\end{definition}

We have the following proposition.
\begin{proposition}\label{cista}
Let $M\in \mathbb{N}$ and
let $F\in \mathcal{PS}^M_L(U)$ (resp. $\mathcal{PS}^M_R(U)$ or $\mathcal{PN}^M(U)$).
Then $F\in \mathcal{PS}^{M+M'}_L(U)$ (resp. $\mathcal{PS}^{M+M'}_R(U)$ or $\mathcal{PN}^{M+M'}(U)$) for all $M'\in \mathbb{N}$.
\end{proposition}
\begin{proof}
It is a direct consequences of the definition.
\end{proof}

\begin{remark}
The restriction of a function $F$ to the complex plane $\mathbb{C}_\mathrm{j}$ is denoted by
$F_\mathrm{j}$.
\end{remark}

\begin{theorem}[Poly splitting lemma and poly decomposition]\label{carac3}
Let $U \subseteq \mathbb{R}^{n+1}$ be an axially symmetric open  set.

{\rm (Ia) (Poly splitting lemma)}.
Let $F\in \mathcal{PS}^M_L(U)$. Then, for every
$\mathrm{j}=\mathrm{j}_1\in\mathbb{S}$ let $\mathrm{j}_2,\ldots, \mathrm{j}_n$ be a completion to a
basis of  $\mathbb{R}_n$ satisfying the defining
relations $\mathrm{j}_r\mathrm{j}_s+\mathrm{j}_s\mathrm{j}_r=-2\delta_{rs}$.
Then, there exist $2^{n-1}$
polyanalytic functions $\mathcal{F}_A\ : U\cap \mathbb{C}_\mathrm{j}\to \mathbb{C}_\mathrm{j}$ such that, for
every $z=u+\mathrm{j}v$, we have
$$
F_\mathrm{j}(z)=\sum_{|A|=0}^{n-1} \mathcal{F}_A(z) \mathrm{j}_A,\quad \mathrm{j}_A=\mathrm{j}_{i_1}\ldots \mathrm{j}_{i_s},
$$
where $A=i_1\ldots i_s$ is a subset of $\{2,\ldots ,n\}$, with ${i_1}<\ldots <{i_s}$,
or, when $|A|=0$, $\mathrm{j}_\emptyset=1$.

{\rm (Ib) (Poly decomposition)}.  The function $F\in\mathcal{PS}^M_L(U)$ if and only if there exist uniquely determined functions $f_0,...,f_{M-1}\in \mathcal{SM}_L(U)$  such that, for $f_{M-1}\not=0$, we have the following decomposition
\begin{equation}\label{strtpolyleft}
F(x)=\displaystyle\sum_{k=0}^{M-1}\overline{x}^kf_k(x), \ \ \ \forall x\in U.
\end{equation}

{\rm (IIa)}
Let $F\in \mathcal{PS}^M_R(U)$ .
 Then, for every
$\mathrm{j}_1,\ldots, \mathrm{j}_n\in\mathbb{S}$  as above
 there exist $2^{n-1}$
polyanalytic functions $\mathcal{F}_A\ : U\cap \mathbb{C}_\mathrm{j}\to \mathbb{C}_\mathrm{j}$ such that, for
every $z=u+\mathrm{j}v$, we have
$$
F_\mathrm{j}(z)=\sum_{|A|=0}^{n-1} \mathrm{j}_A\mathcal{F}_A(z) ,\quad \mathrm{j}_A=\mathrm{j}_{i_1}\ldots \mathrm{j}_{i_s},
$$
where $A=i_1\ldots i_s$ is a subset of $\{2,\ldots ,n\}$, with ${i_1}<\ldots <{i_s}$,
or, when $|A|=0$, $\mathrm{j}_\emptyset=1$.

{\rm (IIb)} The function $F\in \mathcal{PS}^M_R(U)$ if and only if there exist uniquely determined functions $f_0,...,f_{M-1}\in \mathcal{SM}_R(U)$ such that, for $f_{M-1}\not=0$, we have the following decomposition
\begin{equation}\label{strtpolyright}
F(x)=\displaystyle\sum_{k=0}^{M-1}f_k(x) \overline{x}^k, \ \ \ \ \forall x\in U.
\end{equation}

\end{theorem}
\begin{proof}
We consider the case (Ia) and (Ib) since the other ones follows with similar arguments.

\medskip
Step (Ia).
For every
$\mathrm{j}=\mathrm{j}_1\in\mathbb{S}$ let $\mathrm{j}_2,\ldots, \mathrm{j}_n$ be a completion to a
basis of  $\mathbb{R}_n$ satisfying the defining
relations $\mathrm{j}_r\mathrm{j}_s+\mathrm{j}_s\mathrm{j}_r=-2\delta_{rs}$.
Then there exist $2^{n-1}$
 functions $\mathcal{F}_A\ : U\cap \mathbb{C}_\mathrm{j}\to \mathbb{C}_\mathrm{j}$ such that for
every $z=u+\mathrm{j}v$
$$
F_\mathrm{j}(z)=\sum_{|A|=0}^{n-1} \mathcal{F}_A(z) \mathrm{j}_A,\quad \mathrm{j}_A=\mathrm{j}_{i_1}\ldots \mathrm{j}_{i_s},
$$
where $F_\mathrm{j}(z)$ is the restriction to $U\cap \mathbb{C}_\mathrm{j}$,
and $A=i_1\ldots i_s$ is a subset of $\{2,\ldots ,n\}$, with ${i_1}<\ldots <{i_s}$,
or, when $|A|=0$, $\mathrm{j}_\emptyset=1$.

But since $F$ is a left poly slice monogenic function, the functions
$\mathcal{F}_A(z)$ are complex polyanalytic functions of order $M$.

\medskip
Step (Ib).
We show that the function
$F(x)$ defined in (\ref{strtpolyleft})
is poly slice monogenic. First of all, we note that $F$ is a slice function since it is the sum of slice functions. Then, by the definition of poly slice  monogenicity and the product rule, we have that
$$
\frac{1}{2^M}(\partial_u+\mathrm{j}\partial_v)^M
F(u+\mathrm{j}v)=\sum_{k=0}^{M-1}\frac{1}{2^M}(\partial_u+\mathrm{j}\partial_v)^M
\Big((u-\mathrm{j}v)^kf_k(u+\mathrm{j}v)\Big)=0.
$$
Viceversa, let us assume that
$F(x)=F(u+\mathrm{j}v)=F_0(u,v)+\mathrm{j}F_1(u,v)$ is a left poly slice
monogenic function of order $M$, i.e. the pair $(F_0,F_1)$ is an even-odd pair satisfying
\eqref{polyCR}.
  By fixing a basis $e_1,\ldots, e_n$ of $\mathbb{R}_n$ we can write
$F_i=\sum_{|A|=0}^n F_{i,A} e_A$, $i=0,1$ where the functions $F_{i,A}$
are real-valued, and by the linear independence of the basis elements
$e_A$, the system \eqref{polyCR} can be rewritten in terms of the $2^n$
real components $F_{i,A}$ of $F_i$, $i=0,1$.
Thus if $F(u+\mathrm{j}v)=F_0(u,v)+\mathrm{j}F_1(u,v)$ is left poly
slice monogenic, each function $F_A=F_{0,A}+\mathrm{j} F_{1,A}$ is a slice function and polyanalytic. By
the classical result applied to the $\mathbb{C}_{\mathrm{j}}$-valued
function $F_A$, we have
$$
F_{A}(u+\mathrm{j}v)=\sum_{k=0}^{M-1} (u-\mathrm{j}v)^k
f_{k,A}(u+\mathrm{j}v)
$$
where the functions $f_{k,A}$ are $\mathbb{C}_{\mathrm j}$-valued,
satisfy the Cauchy-Riemann system and are even-odd in the variables $u,v$ by direct verification.
We thus obtain
\[
\begin{split}
F(u+\mathrm{j}v)&= \sum_{|A|=0}^n F_A(u+\mathrm{j}v) e_A =
\sum_{|A|=0}^n\sum_{k=0}^{M-1} (u-\mathrm{j}v)^k f_{k,A}(u+\mathrm{j}v)e_A\\
&= \sum_{k=0}^{M-1} (u-\mathrm{j}v)^k f_{k}(u+\mathrm{j}v)
\end{split}
\]
where we set  $f_k= \sum_{|A|=0}^n f_{k,A}e_A$. The functions $f_k$ are
evidently left slice monogenic and this concludes the proof.
 \end{proof}
\begin{definition}
The functions $f_0,...,f_{M-1}\in \mathcal{SM}_L(U)$ that appear in Theorem \ref{carac3}
in the poly decomposition (Ib)
$$
F(x)=\displaystyle\sum_{k=0}^{M-1}\overline{x}^kf_k(x),\ \ \ \ \ \forall x\in U,
$$
are called the (left monogenic) components of the  left poly slice monogenic function $F$.
Similarly, we will call the functions $f_0,...,f_{M-1}\in \mathcal{SM}_R(U)$ in (IIb) the components
 of the right poly slice monogenic function $F$.
\end{definition}
\begin{assumption}\label{assumption}
In the sequel, when dealing with a function $F$ in $\mathcal{PS}^M_L(U)$ or $\mathcal{PS}^M_R(U)$, we always assume that $F$ is written via its poly decomposition.
\end{assumption}
\begin{proposition}\label{poly-dec-int}
Let $U\subseteq \mathbb{R}^{n+1}$ be an axially symmetric  domain and let $M\in \mathbb{N}$. Then we have:

(I)
  $F\in \mathcal{PS}^M_L(U)$ is intrinsic if and only if all its left slice monogenic components are also  intrinsic.

(II) $F\in \mathcal{PS}^M_R(U)$ is intrinsic if and only if all its right slice monogenic components are also  intrinsic.
\end{proposition}
\begin{proof} We prove (I).
We use the poly decomposition in Theorem \ref{carac3} to write
 $$
 F(x)=\sum_{k=0}^{M-1}\overline{x}^kf_k(x),
 $$
 where $f_k$ are the left slice monogenic components for all $k=0,..,M-1$.
 First, we observe that if all the functions $f_k$ are intrinsic,
  then $F$ will preserve any complex plane $U\cap\mathbb{C}_\mathrm{j}$
  that is  to say that $F(U\cap \mathbb{C}_\mathrm{j})\subseteq \mathbb{C}_\mathrm{j}$, so it is intrinsic.

For the converse, let us assume that $F=F_0+\mathrm{j}F_1$ is intrinsic
so that $F_0$, $F_1$ are real valued. This means that in the decomposition
$F=\sum_{|A|=0}^n (F_{0,A}+\mathrm{j}F_{1,A})e_A$ there is only the term
corresponding to $e_\emptyset=1$ and expanding in the form
$F(u+\mathrm{j}v)=\sum_{k=0}^{M-1}(u-\mathrm{j}v)^k f_k(u+\mathrm{j}v)$
the functions $f_k$ are $\mathbb{C}_{\mathrm{j}}$-valued, see the proof of Theorem \ref{carac3} (Ib), and so they
take $\mathbb{C}_{\mathrm{j}}$ into itself and so they are intrinsic.
\end{proof}

\begin{remark}
{\rm

(I) Observe that since $f_k(x)$ is slice monogenic, for $\ell, k\in \mathbb{N}$, for $k\geq \ell$,
we have that
$$
\overline{\partial_\mathrm{j}}^\ell \big(\overline{x}^k f_k(x)\big)=
\overline{\partial_\mathrm{j}}^\ell\big(\overline{x}^k \big)f_k(x)=
k(k-1)(k-2)\ldots (k-\ell+1 )\overline{x}^{k-\ell} f_k(x)=\frac{k!}{(k-\ell)!}\overline{x}^{k-\ell} f_k(x),
$$
and $\overline{\partial_\mathrm{j}}^\ell \big(\overline{x}^k f_k(x)\big)=0$, for $k< \ell$.

(II) Similarly, for $\ell, k\in \mathbb{N}$, for $k\geq \ell$, we have
$$
 \big( f_k(x)\overline{x}^k\big)\overline{\partial_\mathrm{j}}^\ell=
f_k(x)\big(\overline{x}^k \big)\overline{\partial_\mathrm{j}}^\ell=
k(k-1)(k-2)\ldots (k-\ell+1 ) f_k(x)\overline{x}^{k-\ell}
=\frac{k!}{(k-\ell)!}f_k(x)\overline{x}^{k-\ell},
$$
and $ \big( f_k(x)\overline{x}^k\big)\overline{\partial_\mathrm{j}}^\ell=0$, for $k< \ell$.
}
\end{remark}

The next result was already proved for quaternions, we revise its proof in more details here for the sake of completeness.
\begin{proposition}\label{Fgprop}

Let $U\subseteq \mathbb{R}^{n+1}$ be an axially symmetric  open set and $M\in \mathbb{N}$.
Then, denoting by $Fg$ or $gF$ the pointwise product, we have the following statements.

(Ia)
Let $F\in \mathcal{PN}^M(U)$ and $g\in \mathcal{SM}_L(U)$.
Then,  $Fg$ belongs to $\mathcal{PM}^M_L(U)$.

(Ib)
Let $F\in \mathcal{PM}^M_L(U)$ and $g\in \mathcal{N}(U)$.
Then,  $gF$ belongs to $\mathcal{PM}^M_L(U)$.

(IIa)
Let $F\in \mathcal{PN}^M(U)$ and $g\in \mathcal{SM}_R(U)$.
Then,  $gF$ belongs to $\mathcal{PM}^M_R(U)$.

(IIb)
 Let
$F\in \mathcal{PM}^M_R(U)$ and $g\in \mathcal{N}(U)$.
Then,  $Fg$ belongs to $\mathcal{PM}^M_R(U)$.

(III)
Let $F\in \mathcal{PN}^M(U)$ and $g\in \mathcal{N}(U)$.
Then, $gF=Fg$  belongs to $\mathcal{PN}^M(U)$.

\end{proposition}
\begin{proof}
We prove just step (Ia). In much the same way we can prove the other points.
So we assume that $F\in \mathcal{PN}^M(U)$, $g\in \mathcal{SM}_L(U)$  and  $\mathrm{j} \in \mathbb{S}$ and we set $x=u+\mathrm{j}v$. We will prove that \begin{equation}\label{FgP}
\overline{\partial_\mathrm{j}}^M(Fg)(u+v\mathrm{j})=0.
\end{equation}
Indeed, first we note that since $F$ is intrinsic, we have $F(U\cap \mathbb{C}_\mathrm{j})\subset \mathbb{C}_\mathrm{j}.$ In particular, by the Leibniz rule, the equality
$$
\overline{\partial_\mathrm{j}}(Fg)(u+\mathrm{j}v)
=F\overline{\partial}_\mathrm{j}(g)(u+v\mathrm{j})+\overline{\partial_\mathrm{j}}(F)g(u+v\mathrm{j})
$$
holds.
We note that since $F$ is poly slice monogenic of order $M$ and $g$ is slice monogenic
 we have $\overline{\partial_\mathrm{j}}(g)=0 \text{  and } \overline{\partial_\mathrm{j}}(F)\neq 0.$ Thus, we obtain
$$
\overline{\partial_\mathrm{j}}(Fg)(u+\mathrm{j}v)=\overline{\partial_\mathrm{j}}(F)g(u+v\mathrm{j}).
$$
Then, since $F$ is intrinsic we can use the Leibniz rule $M$ times and get
$$
\overline{\partial_\mathrm{j}}^M(Fg)(u+\mathrm{j}v)=\overline{\partial_\mathrm{j}}^{M-1}(F)\overline{\partial_\mathrm{j}}(g)(u+v\mathrm{j})
+\overline{\partial_\mathrm{j}}^{M}(F)g(u+v\mathrm{j}).
$$
Therefore, it follows that the formula \eqref{FgP} holds since $F\in\mathcal{PN}^M(U)$ and $g\in\mathcal{SM}_L(U)$.
Hence, the pointwise product $Fg$ is poly slice monogenic of order $M$ on $U$.

\end{proof}

Using the explicit poly decomposition of a given poly slice monogenic function $F$ of order $M$, see Assumption \ref{assumption}, we can give it an integral representation using the Cauchy formula for the slice monogenic components $\{f_k(x)\}_{k=0,...,M-1}$.

 \begin{corollary}[Slice monogenic integral representation of poly slice monogenic functions]
\label{Cauchygenerale}
 Let $U\subseteq
\mathbb{R}^{n+1}$ be a slice Cauchy domain.
Let $\mathrm{j}\in\mathbb{S}$ and set  $ds_\mathrm{j}=ds (-\mathrm{j})$.

(I) If $F$ is a (left) poly slice monogenic function of order $M$
on a set that contains $\overline{U}$ then
\begin{equation}\label{cauchynuovoLLL}
F(x)=\displaystyle
\frac{1}{2 \pi}\int_{\partial (U\cap \mathbb{C}_\mathrm{j})}\sum_{k=0}^{M-1}\overline{x}^k S_L^{-1}(s,x)\, ds_\mathrm{j}\,  f_k(s),\qquad\text{for any }\ \  x\in U.
\end{equation}

(II) If $F$ is a right slice monogenic function on a set that contains $\overline{U}$,
then
\begin{equation}\label{CauchyrightRRR}
 F(x)=\displaystyle\frac{1}{2 \pi}\int_{\partial (U\cap \mathbb{C}_\mathrm{j})} \sum_{k=0}^{M-1} f_k(s)\, ds_\mathrm{j}\, S_R^{-1}(s,x) \overline{x}^k,\qquad\text{for any }\ \  x\in U.
 \end{equation}
The integrals in (\ref{cauchynuovoLLL}) and (\ref{CauchyrightRRR}) depend neither on $U$ nor on the imaginary unit $\mathrm{j}\in\mathbb{S}$.
\end{corollary}
\begin{proof}
We consider the case (\ref{cauchynuovoLLL}) since (\ref{CauchyrightRRR}) can be obtained with similar considerations.
So  (\ref{cauchynuovoLLL}) follows by replacing the Cauchy formula for the slice monogenic components  $\{f_k(x)\}_{k=0,...,M-1}$ given by
$$
f_k(x)=\frac{1}{2 \pi}\int_{\partial (U\cap \mathbb{C}_\mathrm{j})} S_L^{-1}(s,x)\, ds_\mathrm{j}\,  f_k(s),\qquad\text{for any }\ \  x\in U
$$
into the poly decomposition formula
$
F(x)=\displaystyle\sum_{k=0}^{M-1}\overline{x}^kf_k(x).
$
 \end{proof}
 Similarly we can prove the case of unbounded domains.
 \begin{theorem}[Slice monogenic integral representation of poly slice monogenic functions on unbounded slice Cauchy domains]
Let $U\subset \mathbb{R}^{n+1}$ be an unbounded slice Cauchy domain and let $\mathrm{j}\in\mathbb{S}$.
Let $F\in\mathcal{PS}_{L}(\overline{U})$
 with poly slice decomposition
 $$
F(x)=\sum_{k=0}^{M-1}\overline{x}^kf_k(x)
$$
and the components $f_k\in\mathcal{SM}_{L}(\overline{U})$ are such that the limits
$
\lim_{|x|\to\infty}f_k(x)=f_k(\infty)
$
exist and are finite for all $k=0,...,M-1$.
Then, we have
\[
F(x)=\sum_{k=0}^{M-1}\overline{x}^kf_k(\infty) + \frac{1}{2\pi}\int_{\partial(U\cap\mathbb{C}_{\mathrm{j}})}\sum_{k=0}^{M-1}\overline{x}^kS_L^{-1}(s,x)\,ds_{\mathrm{j}}\, f_k(s) \qquad\text{for any }\ \ x\in U.
\]

Let $F\in\mathcal{PS}_{R}(\overline{U})$
 with poly slice decomposition
 $$
F(x)=\sum_{k=0}^{M-1}f_k(x)\overline{x}^k
$$
and the components $f_k\in\mathcal{SM}_{R}(\overline{U})$ are such that the limits
$
\lim_{|x|\to\infty}f_k(x)=f_k(\infty)
$
exist and are finite for all $k=0,...,M-1$.
Then, we have
\[
F(x)=
\sum_{k=0}^{M-1}f_k(\infty)\overline{x}^k
 + \frac{1}{2\pi}\int_{\partial(U\cap\mathbb{C}_j)}\sum_{k=0}^{M-1}f(s)\,ds_j\,S_R^{-1}(s,x)\overline{x}^k \qquad\text{for any }\ \ x\in U.
\]
 \end{theorem}

\begin{remark}{\rm
Observe that we have two possibilities to write the Cauchy kernels, using the form I or the form II.
}
\end{remark}

The following result is a direct consequence of the form of the
 slice function, i.e.,
 \[
 F(x) = F_{0}(u,v) + \mathrm{j}F_{1}(u,v)\qquad \text{for } x = u + \mathrm{j} v\in U
 \]
with the two functions $F_{0},F_{1}: \mathcal{U}\to \mathbb{R}_n$ that satisfy the compatibility condition (see the definition of poly slice monogenic functions).
We recall it for the reader's convenience.

 \begin{theorem}[Poly Structure (or Poly Representation) Formula]\label{formulpolyamon}
Let $U\subseteq
\mathbb{R}^{n+1}$ be an axially symmetric domain.

\medskip
(I) Let $F\in \mathcal{PS}^M_L(U)$. Then,
 for any  vector
$x =u+\mathrm{j}_x v\in U$, the following formula holds:
\begin{equation}\label{polyleftrep}
F(x)=\frac{1}{2}\Big[1 -  \mathrm{j}_x \mathrm{j}\Big]F(u +\mathrm{j} v)+\frac{1}{2}\Big[1 +  \mathrm{j}_x \mathrm{j}\Big]F(u -\mathrm{j} v),  \ \ {\rm for \ all} \ \ u +\mathrm{j} v\in U, \ \ \mathrm{j}\in \mathbb{S}.
\end{equation}

\medskip
(II) Let $F\in \mathcal{PS}^M_R(U)$. Then,
 for any  vector
$x =u+\mathrm{j}_x v\in U$, the following formula holds:
\begin{equation}\label{polyleftrep}
F(x)=\frac{1}{2}F(u +\mathrm{j} v)\Big[1 -  \mathrm{j}_x \mathrm{j}\Big]+\frac{1}{2}F(u -\mathrm{j} v)\Big[1 +  \mathrm{j}_x \mathrm{j}\Big],  \ \ {\rm for \ all} \ \ u +\mathrm{j} v\in U, \ \ \mathrm{j}\in \mathbb{S}.
\end{equation}
\end{theorem}

In the following we use the notation $\overline{\mathbb{R}^{n+1}}=\mathbb{R}^{n+1}\cup\lbrace \infty \rbrace$. Then, we first recall the Runge's theorem for slice monogenic functions which was proved in \cite{RUNGE}. We refer to this paper also for the terminology.

\begin{theorem}\label{smRunge}
Let $K$ be an axially symmetric compact set in $\mathbb{R}^{n+1}$, and let $A$ be a set having a point in each connected component of $\overline{\mathbb{R}^{n+1}}\setminus K$. For any axially symmetric open set $U\supset K$, for every $f\in\mathcal{SM}_L(U)$ and for every $\varepsilon>0$ there exists a rational function $r$ whose poles are spheres in $A$ such that
$$|f(x)-r(x)|<\varepsilon,$$
for all $x\in K$. Similar considerations hold for  every $f\in\mathcal{SM}_R(U)$.
\end{theorem}
Now, we are going to use Theorem \ref{smRunge} to prove the poly slice monogenic counterpart of the Runge's theorem. First, we give the definition of  poly slice monogenic rational function.

\begin{definition}
We say that a left poly slice monogenic function  $R(x)$ is rational
if the left slice  monogenic components $r_k(x)$ in the decomposition
\begin{equation}
R(x)=\displaystyle \sum_{k=0}^{M-1}\overline{x}^kr_k(x), \quad \forall x\in U
\end{equation}
 are rational.
 We say that a right poly slice monogenic function  $R(x)$ is rational
if the right slice  monogenic components $r_k(x)$ in the decomposition
\begin{equation}
R(x)=\displaystyle \sum_{k=0}^{M-1}r_k(x) \overline{x}^k, \quad \forall x\in U
\end{equation}
are rational.
\end{definition}
\begin{theorem}[Poly slice Runge's theorem]\label{smpRunge}
Let $M\in \mathbb{N}$.
Let $K$ be an axially symmetric compact set in $\mathbb{R}^{n+1}$, and let $A$ be a set having a point in each connected component of $\overline{\mathbb{R}^{n+1}}\setminus K$. For any axially symmetric  domain $U\supset K$, for every $F\in\mathcal{PS}_L^M(U)$ and for every $\varepsilon>0$ there exists a left poly rational function $R$ whose poles are spheres in $A$ such that
$$|F(x)-R(x)|<\varepsilon,$$
for all $x\in K$. The same approximation holds for $F\in\mathcal{PS}_R^M(U)$. In the case $F$ is intrinsic the rational functions $R$ are also intrinsic.
\end{theorem}
\begin{proof}
We show just the case when $F\in\mathcal{PS}_L^M(U)$, the other statements follows in much the same way.
We note that $K$ is a compact, in particular $K$ is bounded so that it is contained in some ball centered at the origin, i.e., $K\subset \overline{B(0,\rho)},$
for some  $\rho>0.$ Furthermore, since $F$ is poly slice monogenic of order $M$ on the domain $U$, we know the validity of the poly decomposition (Ib) in Theorem \ref{carac3}
with all $\{f_k\}_{0 \leq k\leq M-1}$ that are slice monogenic functions on $U$ and $f_{M-1}\neq 0$.
In particular, since $K$ is contained in $U$, we also have
\begin{equation}
F(x)=\displaystyle \sum_{k=0}^{M-1}\overline{x}^kf_k(x), \quad \forall x\in K.
\end{equation}

Then, we can apply the Runge's theorem \ref{smRunge} on each slice monogenic function $f_k$. Thus, for all $k=0,...,M-1$, we know that, for every $\varepsilon>0$, there exist a rational function $r_k$ whose poles are spheres in $A$ such that
\begin{equation}\label{InequalityR1}
\displaystyle |f_k(x)-r_k(x)|<\varepsilon\left(\frac{1}{\sum_{j=0}^{M-1}\rho^j}\right),
\end{equation}
for every $x\in K$.
Now, we consider the poly slice monogenic rational function given by
$$
R(x)=\displaystyle \sum_{k=0}^{M-1}\overline{x}^kr_k(x), \quad \forall x\in U.
$$
Let $\varepsilon>0$. Then, we have
\[ \begin{split}
|F(x)-R(x)|& = \left|\sum_{k=0}^{M-1}\overline{x}^k(f_k(x)-r_k(x)\right|
\leq \sum_{k=0}^{M-1}|x|^k|f_k(x)-r_k(x)|,
\end{split}
\]
 for every $x\in K$.
Therefore, we apply the inequality \eqref{InequalityR1} to each slice monogenic component and get
\begin{equation}\label{InequalityR2}
 |F(x)-R(x)|< \dfrac{\varepsilon}{\sum_{k=0}^{M-1}\rho^k}\left(\sum_{k=0}^{M-1}|x|^k\right),
\end{equation}
 for every $x\in K$.
However, we know that $K\subset \overline{B(0,\rho)}$. So we have $|x|\leq\rho$, for any $x\in K$. In particular, this shows that
$$ \sum_{k=0}^{M-1}|x|^k\leq \sum_{k=0}^{M-1}\rho^k,$$
for every $x\in K$.
Therefore, inserting this fact in the inequality \eqref{InequalityR2} we obtain that
$|F(x)-R(x)|<\varepsilon,$
for every $x\in K$.
\end{proof}

\section{Cauchy formulas and product of poly slice monogenic functions}\label{sec4}

In this section we develop the Cauchy formulas for poly slice monogenic functions. These will be used in the next section to define the $PS$-functional calculus for noncommuting operators.

\begin{lemma}[Poly Cauchy integral formula]\label{polycaulem}
Let $U\subseteq\mathbb{R}^{n+1}$ be a slice Cauchy domain, assume that
 $F\in \mathcal{PS}^M_L(U)$  and $G\in \mathcal{PS}^M_R(U)$ for some $M\in \mathbb{N}$. Then we have
\begin{equation}\label{cucthpoly}
\int_{\partial (U\cap \mathbb{C}_\mathrm{j})}
 \sum_{\ell=0}^{M-1}(-1)^\ell G(s)\overline{\partial}_\mathrm{j}^{M-\ell-1} \, ds_\mathrm{j}\,
\overline{\partial}_\mathrm{j}^{\ell}F(s)=0,
\end{equation}
where $ds_\mathrm{j}=ds (-\mathrm{j})$ and  $\overline{\partial}_\mathrm{j}:=\frac{1}{2}(\partial_u+\mathrm{j}\partial_v)$ for $\mathrm{j}\in \mathbb{S}$.
\end{lemma}
\begin{proof}
By writing $G=\sum_{|A|=0}^n e_AG_A$, $F=\sum_{|A|=0}^n F_A e_A$, see the proof of Theorem \ref{carac3}, (Ib), the result follows from the analogue theorem for polyanalytic functions of a complex variable.
\end{proof}

\subsection{Cauchy formulas with kernels $P_\ell S_L^{-1}$ and $P_\ell S_R^{-1}$ }

\medskip

The poly slice monogenic Cauchy kernel are described in the next result:
\begin{definition}[Cauchy kernels $P_\ell S_L^{-1}$ and $P_\ell S_R^{-1}$ ]
Let $\ell\in \mathbb{N}$. We define
the left poly slice  monogenic Cauchy kernels $P_\ell S_L^{-1}$ and $P_\ell S_R^{-1}$ of order $\ell+1$ by
\[
\begin{split}
P_\ell S^{-1}_{L}(s,x):&=\frac{(Re(s-x))^\ell}{\ell!}S_L^{-1}(s,x)
\\
&
=-\frac{(Re(s-x))^\ell}{\ell!}
(x^2 -2 \Re  (s) x+|s|^2)^{-1}(x-\overline s),\ \ \ s\not\in [x]
\end{split}
\]
and right poly slice monogenic Cauchy kernels of order $\ell$ is defined by
\[
\begin{split}
P_\ell S^{-1}_{R}(s,x):&=\frac{(Re(s-x))^\ell}{\ell!}S^{-1}_{R,\ell}(s,x)
\\
&
=-\frac{(Re(s-x))^\ell}{\ell!}(x-\bar s)(x^2-2\Re (s)x+|s|^2)^{-1},\ \ \ s\not\in [x].
\end{split}
\]
\end{definition}
\begin{remark}
Observe that, in both cases, we have used the slice monogenic Cauchy kernels written in form I so that the $PS$-functional calculus will work for paravector operators with noncommuting components. For the function theory it is also possible to use slice monogenic Cauchy kernels written in form II, but in this case to define the $PS$-functional calculus we are limited to commuting operators.
\end{remark}

\begin{lemma}
Let $\ell\in \mathbb{N}$ and let $x,s\in \mathbb{R}^{n+1}$ with $s\notin [x]$.
The kernel $P_\ell S_L^{-1}(s,x)$ is left poly slice monogenic in $x$ and right poly
slice monogenic in $s$ of order $\ell+1$.
The  kernel $P_\ell S_R^{-1}(s,x)$ is left poly  slice monogenic in $s$ and right poly
slice monogenic in $x$ of order $\ell+1$.
\end{lemma}
\begin{proof}
For $\ell\in \mathbb{N}$. For all $s,x\in\mathbb{R}^{n+1}$, we set
$$
P_\ell(s,x)=\frac{(Re(s-x))^\ell}{\ell!} .
$$
Then, $P_\ell(s,x)$ is a slice function, intrinsic and poly slice monogenic of order $\ell+1$ with respect to the variables $x$ and $s$.
 By a direct computation, for $s\notin [x]$, the product $P_\ell(s,x)S_L^{-1}(s,x)$
is left poly slice monogenic in $x$ and right poly slice monogenic in $s$ of order $\ell+1$.
 Using similar arguments  we treat  $P_\ell(s,x) S^{-1}_{R}(s,x)$.
\end{proof}

\begin{theorem}[The poly slice Cauchy formulas with kernels $P_\ell S_L^{-1}$ and $P_\ell S_R^{-1}$]
\label{Cauchygeneralep}
Let $U\subset\mathbb{R}^{n+1}$ be a slice Cauchy domain. For $\mathrm{j}\in\mathbb{S}$  set  $ds_\mathrm{j}=ds (-\mathrm{j})$  and  $\overline{\partial}_\mathrm{j}:=\frac{1}{2}(\partial_u+\mathrm{j}\partial_v)$.

(I) If $F$ is a left poly slice monogenic function on a set that contains $\overline{U}$ then
\begin{equation}\label{cauchynuovo}
 F(x)=\frac{1}{2 \pi}\int_{\partial (U\cap \mathbb{C}_\mathrm{j})}\sum_{\ell=0}^{M-1}(-2)^\ell P_\ell S^{-1}_{L}(s,x)\, ds_\mathrm{j}\, \overline{\partial_\mathrm{j}}^\ell F(s),\qquad\text{for any }\ \  x\in U.
\end{equation}

(II) If $F$ is a right poly slice monogenic function on a set that contains $\overline{U}$,
then
\begin{equation}\label{Cauchyright}
 F(x)=\frac{1}{2 \pi}\int_{\partial (U\cap \mathbb{C}_\mathrm{j})}\sum_{\ell=0}^{M-1}(-2)^\ell F(s)\overline{\partial_\mathrm{j}}^\ell\, ds_\mathrm{j}\,P_\ell S^{-1}_{R}(s,x),\qquad\text{for any }\ \  x\in U.
 \end{equation}
The integrals (\ref{cauchynuovo}) and (\ref{Cauchyright}) depend neither on $U$ nor on the imaginary unit $\mathrm{j}\in\mathbb{S}$.
\end{theorem}
\begin{proof} Consider (I).
It is a consequence of the Lemma \ref{polycaulem} and of Theorem \ref{formulpolyamon}.
We now use the poly splitting lemma (Ia).  Let $F\in \mathcal{PS}^M_L(U)$, and for every
$\mathrm{j}=\mathrm{j}_1\in\mathbb{S}$ let $\mathrm{j}_2,\ldots, \mathrm{j}_n$ be a completion to a
basis of  $\mathbb{R}_n$ satisfying the defining
relations $\mathrm{j}_r\mathrm{j}_s+\mathrm{j}_s\mathrm{j}_r=-2\delta_{rs}$.
Then there exist $2^{n-1}$
polyanalytic functions $\mathcal{F}_A\ : U\cap \mathbb{C}_\mathrm{j}\to \mathbb{C}_\mathrm{j}$ such that, for
every $z=u+\mathrm{j}v$, we have
$$
F_\mathrm{j}(z)=\sum_{|A|=0}^{n-1} \mathcal{F}_A(z) \mathrm{j}_A,\quad \mathrm{j}_A=\mathrm{j}_{i_1}\ldots \mathrm{j}_{i_s},
$$
where $A=i_1\ldots i_s$ is a subset of $\{2,\ldots ,n\}$, with ${i_1}<\ldots <{i_s}$,
or, when $|A|=0$, $\mathrm{j}_\emptyset=1$.
So we can write
$$
\mathcal{F}_A(z)=
\int_{\partial (U\cap \mathbb{C}_\mathrm{j})} \sum_{\ell=0}^{M-1}(-2)^\ell \frac{1}{2\pi \mathrm{j}}\frac{1}{(s-z)}\frac{({\rm Re}(s-z))^\ell}{\ell !}ds
\overline{\partial_\mathrm{j}}^\ell \mathcal{F}_A(s), \  \ \ \ z\in U\cap \mathbb{C}_\mathrm{j}
$$
and also
$$
F_\mathrm{j}(z)=\sum_{|A|=0}^{n-1}\mathcal{F}_A(z)\mathrm{j}_A=\frac{1}{2\pi}\int_{\partial (U\cap \mathbb{C}_\mathrm{j})} \sum_{\ell=0}^{M-1}(-2)^\ell
\frac{1}{   s-z}\frac{({\rm Re}(s-z))^\ell}{\ell !}\, ds_\mathrm{j} \,
\overline{\partial_\mathrm{j}}^\ell \, \sum_{|A|=0}^{n-1}\mathcal{F}_A(s)\mathrm{j}_A, \  \ \ \ z\in U\cap \mathbb{C}_\mathrm{j}.
$$
Thus we can write $F_\mathrm{j}(z)$ and $F_\mathrm{j}(\overline{z})$ as
$$
F_\mathrm{j}(z)=\frac{1}{2\pi}\int_{\partial (U\cap \mathbb{C}_\mathrm{j})} \sum_{\ell=0}^{M-1}(-2)^\ell
\frac{1}{   s-z}\frac{({\rm Re}(s-z))^\ell}{\ell !}\, ds_\mathrm{j} \,
\overline{\partial_\mathrm{j}}^\ell \,F(s), \  \ \ \ z\in U\cap \mathbb{C}_\mathrm{j},
$$
and
$$
F_\mathrm{j}(\overline{z})=\frac{1}{2\pi}\int_{\partial (U\cap \mathbb{C}_\mathrm{j})} \sum_{\ell=0}^{M-1}(-2)^\ell
\frac{1}{   s-\overline{z}}\frac{({\rm Re}(s-\overline{z}))^\ell}{\ell !}\, ds_\mathrm{j} \,
\overline{\partial_\mathrm{j}}^\ell \,F(s), \  \ \ \ \overline{z}\in U\cap \mathbb{C}_\mathrm{j},
$$
where we set  $x=u+\mathrm{j}_x v$,  $z=u+\mathrm{j}v$, $s=s_0+\mathrm{j}s_1$, and where $\mathrm{j}$ is the imaginary unit of the complex plane $\mathbb{C}_\mathrm{j}$ on which we integrate.
Now we use the poly monogenic representation formulas and the two Cauchy formulas
for  $F_\mathrm{j}(z)$ and $F_\mathrm{j}(\overline{z})$ on the complex plane $\mathbb{C}_\mathrm{j}$. Specifically,
we use
$$
F(x)=\frac{1}{2}\Big[1 -  \mathrm{j}_x \mathrm{j}\Big]F_\mathrm{j}(z)+\frac{1}{2}\Big[1 +  \mathrm{j}_x \mathrm{j}\Big]F_\mathrm{j}(\overline{z})
$$
 so we have
\[
\begin{split}
F(x)&=\frac{1}{2}\Big[1 -  \mathrm{j}_x \mathrm{j}\Big]\frac{1}{2\pi}\int_{\partial (U\cap \mathbb{C}_\mathrm{j})} \sum_{\ell=0}^{M-1}(-2)^\ell
\frac{1}{   s-z}\frac{({\rm Re}(s-z))^\ell}{\ell !}\, ds_\mathrm{j} \,
\overline{\partial_\mathrm{j}}^\ell \,F(s)
\\
&
+\frac{1}{2}\Big[1 +  \mathrm{j}_x \mathrm{j}\Big]\frac{1}{2\pi}\int_{\partial (U\cap \mathbb{C}_\mathrm{j})} \sum_{\ell=0}^{M-1}(-2)^\ell
\frac{1}{   s-\overline{z}}\frac{({\rm Re}(s-\overline{z}))^\ell}{\ell !}\, ds_\mathrm{j} \,
\overline{\partial_\mathrm{j}}^\ell \,F(s),
\end{split}
\]
and also
\[
\begin{split}
F(x)&=\frac{1}{2\pi}\int_{\partial (U\cap \mathbb{C}_\mathrm{j})} \sum_{\ell=0}^{M-1}(-2)^\ell
\frac{1}{2}\Big[1 -  \mathrm{j}_x \mathrm{j}\Big]\frac{1}{   s-z}\frac{({\rm Re}(s-z))^\ell}{\ell !}\, ds_\mathrm{j} \,
\overline{\partial_\mathrm{j}}^\ell \,F(s)
\\
&
+\frac{1}{2\pi}\int_{\partial (U\cap \mathbb{C}_\mathrm{j})} \sum_{\ell=0}^{M-1}(-2)^\ell
\frac{1}{2}\Big[1 +  \mathrm{j}_x \mathrm{j}\Big]\frac{1}{   s-\overline{z}}\frac{({\rm Re}(s-\overline{z}))^\ell}{\ell !}\, ds_\mathrm{j} \,
\overline{\partial_\mathrm{j}}^\ell \,F(s).
\end{split}
\]
Finally  we get
\[
\begin{split}
F(x)&=\frac{1}{2\pi}\int_{\partial (U\cap \mathbb{C}_\mathrm{j})} \sum_{\ell=0}^{M-1}(-2)^\ell
\frac{({\rm Re}(s-x))^\ell}{\ell !}S^{-1}_L(s,x)\, ds_\mathrm{j} \,
\overline{\partial_\mathrm{j}}^\ell \,F(s)
\end{split}
\]
where  we have replaced the poly slice monogenic Cauchy kernel
$$
\frac{({\rm Re}(s-x))^\ell}{\ell !}S^{-1}_L(s,x)=\frac{1}{2}\Big[1 -  \mathrm{j}_x \mathrm{j}\Big]\frac{({\rm Re}(s-z))^\ell}{\ell !}\frac{1}{s-z} +\frac{1}{2}\Big[1 +  \mathrm{j}_x \mathrm{j}\Big]\frac{({\rm Re}(s-\overline{z}))^\ell}{\ell !}\frac{1}{s-\overline{z}}
$$
 written via the poly slice monogenic representation formula.
\end{proof}

\begin{definition}
We say that a poly slice monogenic function is poly slice monogenic function at infinity if
its slice monogenic components are slice monogenic at infinity.
\end{definition}

\begin{theorem}[Poly slice monogenic Cauchy formulas on unbounded slice Cauchy domains with kernels $P_\ell S_L^{-1}$ and $P_\ell S_R^{-1}$]\label{CauchyOutsidepoly}
Let $U\subset \mathbb{R}^{n+1}$ be an unbounded slice Cauchy.
For $\mathrm{j}\in\mathbb{S}$  set  $ds_\mathrm{j}=ds (-\mathrm{j})$  and  $\overline{\partial}_\mathrm{j}:=\frac{1}{2}(\partial_u+\mathrm{j}\partial_v)$.

(I)Let $F\in\mathcal{PS}_{L}(\overline{U})$
 with poly slice decomposition
 $$
F(x)=\sum_{k=0}^{M-1}\overline{x}^kf_k(x)
$$
and the components $f_k\in\mathcal{SM}_{L}(\overline{U})$ are such that the limits
$
\lim_{|x|\to\infty}f_k(x)=f_k(\infty)
$
exist and are finite for all $k=0,...,M-1$.
Then, we have
\[
F(x) =\sum_{k=0}^{M-1} \overline{x}^k f_k(\infty) +\frac{1}{2 \pi}\int_{\partial (U\cap \mathbb{C}_\mathrm{j})}\sum_{\ell=0}^{M-1}(-2)^\ell P_\ell S^{-1}_{L}(s,x)\, ds_\mathrm{j}\, \overline{\partial_\mathrm{j}}^\ell F(s) \qquad\text{for any }\ \ x\in U.
\]

(II) Let $F\in\mathcal{PS}_{R}(\overline{U})$
 with poly slice decomposition
 $$
F(x)=\sum_{k=0}^{M-1}f_k(x)\overline{x}^k
$$
and the components $f_k\in\mathcal{SM}_{R}(\overline{U})$ are such that the limits
$
\lim_{|x|\to\infty}f_k(x)=f_k(\infty)
$
exist and are finite for all $k=0,...,M-1$.
Then we have
\[
F(x) =\sum_{k=0}^{M-1}  f_k(\infty)\overline{x}^k +\frac{1}{2 \pi}\int_{\partial (U\cap \mathbb{C}_\mathrm{j})}\sum_{\ell=0}^{M-1}(-2)^\ell  F(s)\overline{\partial_\mathrm{j}}^\ell\, ds_\mathrm{j}\, P_\ell S^{-1}_{R}(s,x)\qquad\text{for any }\ \ x\in U.
\]
 \end{theorem}
 \begin{proof}
  Let us consider (I).
  For sufficiently large $r>0$ , the set $U_r := U\cap B_{r}(0)$ is a bounded slice Cauchy domain with $x\in U_r$ and $\hh\setminus U_r \subset U$. By  Theorem \ref{Cauchygeneralep}
 \begin{align*}
 F(x) =& \frac{1}{2 \pi}\int_{\partial (U_r\cap \mathbb{C}_\mathrm{j})}\sum_{\ell=0}^{M-1}(-2)^\ell P_\ell S^{-1}_{L}(s,x)\, ds_\mathrm{j}\, \overline{\partial_\mathrm{j}}^\ell F(s)
 \\
=& \frac{1}{2\pi}\int_{\partial(U\cap\cc_{j})} \sum_{\ell=0}^{M-1}(-2)^\ell P_\ell S^{-1}_{L}(s,x)\, ds_\mathrm{j}\, \overline{\partial_\mathrm{j}}^\ell F(s)\\
 &+ \frac{1}{2\pi}\int_{\partial(B_r(0)\cap\cc_{j})}\sum_{\ell=0}^{M-1}(-2)^\ell P_\ell S^{-1}_{L}(s,x)\, ds_\mathrm{j}\, \overline{\partial_\mathrm{j}}^\ell F(s), \qquad\text{for any }\ \ x\in U.
 \end{align*}
 The Cauchy theorem for poly slice monogenic functions implies that we can vary $r$ without changing the value of the second integral. Letting $r$ tend to infinity, we find that the monogenic components converge to $f_k(\infty)$ and we obtain the statement, since
 $$
\lim_{r\to \infty}\frac{1}{2\pi}\int_{\partial(B_r(0)\cap\cc_{j})}\sum_{\ell=0}^{M-1}(-2)^\ell P_\ell S^{-1}_{L}(s,x)\, ds_\mathrm{j}\, \overline{\partial_\mathrm{j}}^\ell F(s)= \sum_{k=0}^{M-1} \overline{x}^k f_k(\infty).
 $$
 \end{proof}
\begin{remark}
There is a more direct proof of the Cauchy formulas with the Cauchy kernels $P_\ell S_L^{-1}$ and $P_\ell S_R^{-1}$.
In fact, the integrals
 $$
J_k(x):=  \frac{1}{2\pi}\int_{\partial(U\cap\cc_{j})}\sum_{\ell=0}^{k-1}(-2)^\ell P_\ell S^{-1}_{L}(s,x)\, ds_\mathrm{j}\, \overline{\partial_\mathrm{j}}^\ell F(s)
 $$
 can be computed directly. As an example consider
 $k=2$.
\[
\begin{split}
{J}_2(x)&=
\frac{1}{2 \pi}\int_{\partial(U\cap\cc_{j})} S^{-1}_{L}(s,x)\, ds_\mathrm{j}\, \Big(\overline{s}f_1(s)+f_0(s)\Big)
\\
&
+\frac{1}{2 \pi}\int_{\partial(U\cap\cc_{j})}(-2)
 (Re(s)-x_0) S_L^{-1}(s,x)\, ds_\mathrm{j}\,f_1(s)
\end{split}
\]
using the relation
$
S_L^{-1}(s,x)s - xS_L^{-1}(s,x) = 1
$
and the Cauchy theorem we get
 \[
\begin{split}
{J}_2(x)=
\frac{1}{2 \pi} \overline{x}\int_{\partial (U\cap\cc_{j})} S^{-1}_{L}(s,x)\, ds_\mathrm{j}\, f_1(s)
+\frac{1}{2 \pi}\int_{\partial(U\cap\cc_{j})}
 S_L^{-1}(s,x)\, ds_\mathrm{j}\,f_0(s)
\end{split}
\]
so  we obtain
$
J_2(x)=\overline{x}f_1(x)+f_0(x).
$

\end{remark}

\subsection{Cauchy formulas with kernels $\Pi_\ell S_L^{-1}$ and $\Pi_\ell S_R^{-1}$ }

As one may easily compute, the components of $P_\ell S_L^{-1}$ (and of $P_\ell S_R^{-1}$) in the poly slice monogenic decomposition do not have a finite limit for $|x|\to\infty$, for every fixed $s$ and so they are not poly slice monogenic at infinity. Thus we introduce other kernels whose components are poly slice monogenic at infinity.

\begin{definition}[The Cauchy kernels $\Pi_\ell S_L^{-1}$ and $\Pi_\ell S_R^{-1}$]
Let $\ell\in \mathbb{N}$. We define
the left poly slice  monogenic Cauchy kernels of order $\ell+1$ by
\[
\begin{split}
\Pi_\ell S^{-1}_{L}(s,x):&
=\frac{1}{\ell!}\sum_{k=0}^\ell \binom{\ell}{k}\overline{x}^kS_L^{-1}(s,x)(-\overline{s})^{\ell-k}
\\
&
=-\frac{1}{\ell!}\sum_{k=0}^\ell \binom{\ell}{k}\overline{x}^k(x^2 -2 \Re(s)x+|s|^2)^{-1}(x-\overline s) (-\overline{s})^{\ell-k}
,\ \ \ s\not\in [x]
\end{split}
\]
and right poly slice monogenic Cauchy kernels of order $\ell+1$ is defined by
\[
\begin{split}
\Pi_\ell S^{-1}_{R}(s,x):&=\frac{1}{\ell!}\sum_{k=0}^\ell \binom{\ell}{k}
(-\overline{s})^{\ell-k}S_R^{-1}(s,x)\overline{x}^k
\\
&
=-\frac{1}{\ell!}\sum_{k=0}^\ell \binom{\ell}{k}
(-\overline{s})^{\ell-k}(x-\bar s)(x^2-2\Re (s)x+|s|^2)^{-1} \overline{x}^k,\ \ \ s\not\in [x].
\end{split}
\]
\end{definition}

\begin{lemma}
Let $\ell\in \mathbb{N}$ and let $x,s\in \mathbb{R}^{n+1}$ with $s\notin [x]$.
The kernel $\Pi_\ell S_L^{-1}(s,x)$ is left poly slice monogenic in $x$ and right poly
slice monogenic in $s$ of order $\ell+1$.
The kernel $\Pi_\ell S_R^{-1}(s,x)$ is left poly  slice monogenic in $s$ and right poly
slice monogenic in $x$ of order $\ell+1$.
\end{lemma}
\begin{proof}
 Consider the kernel $\Pi_\ell S_L^{-1}(s,x)$.
The statement is a direct consequence of the definition of poly slice monogenicity
and of the poly decomposition of poly slice monogenic functions.
In fact, $\Pi_\ell S_L^{-1}(s,x)$  is of the form
$$
\Pi_\ell S^{-1}_{L}(s,x)
=\sum_{k=0}^\ell \psi_\ell(s,x)\overline{s}^{\ell-k}
$$
where the components
$$
\psi_\ell(s,x):=\frac{1}{\ell!}\binom{\ell}{k}\overline{x}^kS_L^{-1}(s,x)(-1)^{\ell-k}
$$
are right  poly slice monogenic function in  the variable $s$ for  $s\notin [x]$. The order $\ell+1$ is clear from the definition. To see that the kernel $\Pi_\ell S_L^{-1}(s,x)$ is left poly slice monogenic in $x$ it is more convenient to write $S_L^{-1}(s,x)$ in the form II and reasoning in a similar way.
The case of $\Pi_\ell S_R^{-1}(s,x)$  follows with similar considerations.
\end{proof}
\begin{lemma}\label{lemmauni}
For $s\notin [x]$ the kernel $\Pi_\ell S^{-1}_{L}(s,x)$ is the unique left (resp. right) poly slice monogenic extension in $x$ (resp. $s$) of the kernel $\pi_\ell(z,s)=\frac{1}{\ell!}(\bar z-\bar s)^\ell (s-z)^{-1}$, $z,s\in\mathbb{C}_{\mathrm{j}}$, for $z\not=s$. Analogously, for $s\notin [x]$, the kernel $\Pi_\ell S^{-1}_{R}(s,x)$ is the unique right (resp. left) poly slice monogenic extension in $x$ (resp. $s$) of the kernel $\pi_\ell(z,s)=\frac{1}{\ell!}(\bar z-\bar s)^\ell (s-z)^{-1}$, $z,s\in\mathbb{C}_{\mathrm{j}}$, $z\not=s$.
\end{lemma}
\begin{proof}
When $x=z\in\mathbb{C}_{\mathrm{j}}$ it is evident that the restriction of $\Pi_\ell S^{-1}_{L}(s,x)$ to $\mathbb{C}_{\mathrm{j}}$ is $\pi_\ell(z,s)$. Since
$\Pi_\ell S^{-1}_{L}(s,x)=\frac{1}{\ell!}\sum_{k=0}^\ell \binom{\ell}{k}\overline{x}^kS_L^{-1}(s,x)(-\overline{s})^{\ell-k}$ and its slice monogenic components are unique, the assertion follows.
\end{proof}
\begin{theorem}[Cauchy formulas with the kernels $\Pi_\ell S_L^{-1}$ and $\Pi_\ell S_R^{-1}$]
\label{CauchygeneralePI}
Let $U\subset\mathbb{R}^{n+1}$ be a slice Cauchy domain. For $\mathrm{j}\in\mathbb{S}$  set  $ds_\mathrm{j}=ds (-\mathrm{j})$  and  $\overline{\partial}_\mathrm{j}:=\frac{1}{2}(\partial_u+\mathrm{j}\partial_v)$.

(I) If $F$ is a left poly slice monogenic function on a set that contains $\overline{U}$ then
\begin{equation}\label{cauchynuovoPI}
 F(x)=\frac{1}{2 \pi}\int_{\partial (U\cap \mathbb{C}_\mathrm{j})}\sum_{\ell=0}^{M-1} \Pi_\ell S^{-1}_{L}(s,x)\, ds_\mathrm{j}\, \overline{\partial_\mathrm{j}}^\ell F(s),\qquad\text{for any }\ \  x\in U.
\end{equation}

(II) If $F$ is a right poly slice monogenic function on a set that contains $\overline{U}$,
then
\begin{equation}\label{CauchyrightPI}
 F(x)=\frac{1}{2 \pi}\int_{\partial (U\cap \mathbb{C}_\mathrm{j})}\sum_{\ell=0}^{M-1} F(s)\overline{\partial_\mathrm{j}}^\ell\, ds_\mathrm{j}\,\Pi_\ell S^{-1}_{R}(s,x),\qquad\text{for any }\ \  x\in U.
 \end{equation}
The integrals (\ref{cauchynuovo}) and (\ref{Cauchyright}) depend neither on $U$ nor on the imaginary unit $\mathrm{j}\in\mathbb{S}$.
\end{theorem}
\begin{proof} It follows in much the same way as the case of
the poly slice Cauchy formulas with kernels $P_\ell S_L^{-1}$ and $P_\ell S_R^{-1}$ (see Theorem \ref{Cauchygeneralep}) using the poly splitting lemmas by taking $x,s\in\mathbb{C}_{\mathrm{j}}$. Then the assertion follows by the poly representation formula, the polyanalytic Cauchy kernel $\pi_\ell$ used by Th\'eodoresco see (\ref{TEODOR}) and Lemma \ref{lemmauni}.

\end{proof}

\begin{theorem}[Poly slice monogenic Cauchy formulas on unbounded slice Cauchy domains]\label{CauchyOutsidepolyPIOUT}
Let $U\subset \mathbb{R}^{n+1}$ be an unbounded slice Cauchy.
For $\mathrm{j}\in\mathbb{S}$  set  $ds_\mathrm{j}=ds (-\mathrm{j})$  and  $\overline{\partial}_\mathrm{j}:=\frac{1}{2}(\partial_u+\mathrm{j}\partial_v)$.

(I) Let $F\in\mathcal{PS}_{L}(\overline{U})$
 with poly slice decomposition
 $$
F(x)=\sum_{k=0}^{M-1}\overline{x}^kf_k(x)
$$
and the components $f_k\in\mathcal{SM}_{L}(\overline{U})$ are such that the limits
$
\lim_{|x|\to\infty}f_k(x)=f_k(\infty)
$
exist and are finite for all $k=0,...,M-1$.
Then, we have
\[
F(x) =\sum_{k=0}^{M-1} \overline{x}^k f_k(\infty) +\frac{1}{2 \pi}\int_{\partial (U\cap \mathbb{C}_\mathrm{j})}\sum_{\ell=0}^{M-1}\Pi_\ell S^{-1}_{L}(s,x)\, ds_\mathrm{j}\, \overline{\partial_\mathrm{j}}^\ell F(s) \qquad\text{for any }\ \ x\in U.
\]

(II) Let $F\in\mathcal{PS}_{R}(\overline{U})$
 with poly slice decomposition
 $$
F(x)=\sum_{k=0}^{M-1}f_k(x)\overline{x}^k
$$
and the components $f_k\in\mathcal{SM}_{R}(\overline{U})$ are such that the limits
$
\lim_{|x|\to\infty}f_k(x)=f_k(\infty)
$
exist and are finite for all $k=0,...,M-1$.
Then we have
\[
F(x) =\sum_{k=0}^{M-1}  f_k(\infty)\overline{x}^k +\frac{1}{2 \pi}\int_{\partial (U\cap \mathbb{C}_\mathrm{j})}\sum_{\ell=0}^{M-1}  F(s)\overline{\partial_\mathrm{j}}^\ell\, ds_\mathrm{j}\, \Pi_\ell S^{-1}_{R}(s,x)\qquad\text{for any }\ \ x\in U.
\]
 \end{theorem}
 \begin{proof} Consider (I).
 Recall that we assume  the
 components $f_k\in\mathcal{SM}_{L}(\overline{U})$ of $F\in\mathcal{PS}_{R}(\overline{U})$ are such that the limits
$
\lim_{|x|\to\infty}f_k(x)=f_k(\infty)
$
exist and are finite for all $k=0,...,M-1$.
For sufficiently large $r>0$ , the set $U_r := U\cap B_{r}(0)$ is a bounded slice Cauchy domain with $x\in U_r$ and $\mathbb{R}^{n+1}\setminus U_r \subset U$. By \ref{Cauchygeneralep}
 \begin{align*}
 F(x) =& \frac{1}{2 \pi}\int_{\partial (U_r\cap \mathbb{C}_\mathrm{j})}\sum_{\ell=0}^{M-1} \Pi_\ell S^{-1}_{L}(s,x)\, ds_\mathrm{j}\, \overline{\partial_\mathrm{j}}^\ell F(s)
 \\
=& \frac{1}{2\pi}\int_{\partial(U\cap\cc_{j})} \sum_{\ell=0}^{M-1} \Pi_\ell S^{-1}_{L}(s,x)\, ds_\mathrm{j}\, \overline{\partial_\mathrm{j}}^\ell F(s)\\
 &+ \frac{1}{2\pi}\int_{\partial(U\cap\cc_{j})}\sum_{\ell=0}^{M-1} \Pi_\ell S^{-1}_{L}(s,x)\, ds_\mathrm{j}\, \overline{\partial_\mathrm{j}}^\ell F(s), \qquad\text{for any }\ \ x\in U.
 \end{align*}
 The Cauchy theorem for poly slice monogenic functions implies that we can vary $r$ without changing the value of the second integral. Letting $r$ tend to infinity, we find that it equals
  $
\sum_{k=0}^{M-1}\overline{x}^kf_k(\infty)
$
and we obtain the statement.
\end{proof}

\begin{remark}
Also for the
Cauchy formulas with the kernels $\Pi_\ell S_L^{-1}$ and $\Pi_\ell S_R^{-1}$ the integrals
 $$
I_k(r):=  \frac{1}{2\pi}\int_{\partial(U\cap\cc_{j})}\sum_{\ell=0}^{k-1} \Pi_\ell S^{-1}_{L}(s,x)\, ds_\mathrm{j}\, \overline{\partial_\mathrm{j}}^\ell F(s)
 $$
 can be computed directly. As an example consider the case $k=2$, so we have
\[
\begin{split}
{I}_2(x)&=
\frac{1}{2 \pi}\int_{\partial(U\cap\cc_{j})} S^{-1}_{L}(s,x)\, ds_\mathrm{j}\, \Big(\overline{s}f_1(s)+f_0(s)\Big)
\\
&
+\frac{1}{2 \pi}\int_{\partial(U\cap\cc_{j})}
 \Big(
S_L^{-1}(s,x)(-\overline{s})+\overline{x}S_L^{-1}(s,x)\Big)\, ds_\mathrm{j}\,f_1(s)
\end{split}
\]
and so we get
$
I_2(x)=\overline{x}f_1(x)+f_0(x)
$
using the Cauchy formula of slice monogenic functions.
\end{remark}

\subsection{Poly $\circledast_{L}$-product and $\circledast_{R}$-product }

We conclude this section with the product of poly slice monogenic functions.
We note that given two poly slice monogenic functions $F$ and $G$ of orders $N$ and $M$, we can use the poly decomposition formulas to define  a natural product. This product take out of the class but it will be useful for the product rule of the $PS$-functional calculus.

\begin{definition}[Poly $\circledast_{L}$-product and $\circledast_{R}$-product and pointwise product ]\label{plustratprd}
Let $U\subset \mathbb{R}^{n+1}$ be an axially symmetric  open set and $M,N\geq 1$.

(I)
Let $F\in \mathcal{PS}^N_L(U)$ and $G\in \mathcal{PS}^M_L(U)$ and let
 $$
F(x)=\displaystyle \sum_{k=0}^{N-1}\overline{x}^kf_k(x) \ \ \  \text{   and   } \ \ \
 G(x)=\displaystyle \sum_{k=0}^{M-1}\overline{x}^kg_k(x),  \ \ \text{ for all } x\in U,
 $$
be their poly decompositions  for $f_0,...,f_{N-1}$ and $g_0,...,g_{M-1}\in \mathcal{SM}_L(U)$.
We defined the poly $\circledast_{L}$-product of $F$ and $G$ by:
\begin{equation}\label{Gstar'}
(F\circledast_{L} G)(x):=\sum_{\ell=0}^{N+M-2}\overline{x}^\ell \left( \sum_{k+h=\ell} (f_k*_L g_h)(x)\right),
\end{equation}
where $*_L$ is the $*$-product of left slice monogenic functions.

(II)
Let $F\in \mathcal{PS}^N_R(U)$ and $G\in \mathcal{PS}^M_R(U)$ and let
 $$
F(x)=\displaystyle \sum_{k=0}^{N-1}f_k(x)\overline{x}^k \ \ \  \text{   and   } \ \ \
 G(x)=\displaystyle \sum_{k=0}^{M-1}g_k(x)\overline{x}^k,  \ \ \text{ for all } x\in U,
 $$
be their poly decompositions  for $f_0,...,f_{N-1}$ and $g_0,...,g_{M-1}\in \mathcal{SM}_R(U)$.
We defined the poly $\circledast_{R}$-product of $F$ and $G$ by:
\begin{equation}\label{Gstar'}
(F\circledast_{R} G)(x):=\sum_{\ell=0}^{N+M-2} \left( \sum_{k+h=\ell} (f_k*_R g_h)(x)\right) \overline{x}^\ell,
\end{equation}
where $*_R$ is the star-product of right slice monogenic functions.

(III)
Let $F\in \mathcal{PN}^N_L(U)$ and $G\in \mathcal{PS}^M_L(U)$ and let
 $$
F(x)=\displaystyle \sum_{k=0}^{N-1}\overline{x}^kf_k(x) \ \ \  \text{   and   } \ \ \
 G(x)=\displaystyle \sum_{k=0}^{M-1}\overline{x}^kg_k(x),  \ \ \text{ for all } x\in U,
 $$
be their poly decompositions  for $f_0,...,f_{N-1}  \in \mathcal{N}(U)$ and $g_0,...,g_{M-1}\in \mathcal{SM}_L(U)$ .
The  pointwise product of $F$ and $G$ is  defined by:
\begin{equation}\label{polyprod}
(F G)(x):=\sum_{\ell=0}^{N+M-2}\overline{x}^\ell \left( \sum_{k+h=\ell} (f_k g_h)(x)\right).
\end{equation}
Similarly we define the pointwise product for  $F\in \mathcal{PS}^N_R(U)$ and $G\in \mathcal{PN}^M_R(U)$.
\end{definition}

\begin{remark}
{\rm
It is clear that the pointwise product of two poly slice monogenic functions does not preserve the poly slice monogenicity, but, as we will see, it does when we consider intrinsic functions.
}
\end{remark}

\begin{theorem}
Let $U\in \mathbb{R}^{n+1}$ be an axially symmetric  open set and $M,N\geq 1$.

(I)
Let $F\in \mathcal{PS}^N_L(U)$ and $G\in \mathcal{PS}^M_L(U)$ and let

\begin{equation}\label{21}
(F\circledast_{L} G)(x):=\sum_{\ell=0}^{N+M-2}\overline{x}^\ell \left( \sum_{k+h=\ell} (f_k*_L g_h)(x)\right),
\end{equation}
 the poly $\circledast_{L}$-product of $F$ and $G$. Then $(F\circledast_{L} G)\in \mathcal{PS}^{M+N-1}_L(U)$.

(II)
Let $F\in \mathcal{PS}^N_R(U)$ and $G\in \mathcal{PS}^M_R(U)$ and let

\begin{equation}\label{Gstar'}
(F\circledast_{R} G)(x):=\sum_{\ell=0}^{N+M-2} \left( \sum_{k+h=\ell} (f_k*_R g_h)(x)\right) \overline{x}^\ell,
\end{equation}
be the poly $\circledast_{R}$-product of $F$ and $G$. Then $(F\circledast_{R} G)\in \mathcal{PS}^{M+N-1}_R(U)$.
\end{theorem}

\begin{proof} Let us prove (I). With similar computations we get (II).
Observe that, in (\ref{21}),  the function
$$
h_\ell(x):= \sum_{k+h=\ell} (f_k*_L g_h)(x)
$$
is slice monogenic in $U$ for all $\ell$, by definition, so the function
$$
(F\circledast_{L} G)(x):=\sum_{\ell=0}^{N+M-2}\overline{x}^\ell h_\ell(x),
$$
belongs to $ \mathcal{PS}^{M+N-1}_R(U)$ thanks to the poly decomposition theorem.
\end{proof}
As a consequence of the previous result we have the case of intrinsic functions which will be used for the product rule for the $PS$-functional calculus.

\begin{corollary}\label{poinprodpoly}
Let $U\subseteq \mathbb{R}^{n+1}$ be an axially symmetric open set and $M,N\geq 1$.

(I)
Let $F\in \mathcal{PN}^N_L(U)$ and $G\in \mathcal{PS}^M_L(U)$. Then we have

\begin{equation}\label{21dd}
(F\circledast_{L} G)(x)= F(x)G(x)=\sum_{\ell=0}^{N+M-2}\overline{x}^\ell \left( \sum_{k+h=\ell} f_k(x) g_h(x)\right),
\end{equation}
 and $FG\in \mathcal{PS}^{M+N-1}_L(U)$.

(II)
Let $F\in \mathcal{PS}^N_R(U)$ and $G\in \mathcal{PN}^M_R(U)$. Then we have

\begin{equation}\label{Gstar'ddd}
(F\circledast_{R} G)(x)=F(x)G(x)=\sum_{\ell=0}^{N+M-2} \left( \sum_{k+h=\ell} f_k(x)g_h(x)\right) \overline{x}^\ell,
\end{equation}
and  $F G\in \mathcal{PS}^{M+N-1}_R(U)$.
\end{corollary}
\begin{proof}
It is a direct consequence of the  product theorem for slice monogenic functions, i.e.,
it is the case when the slice monogenic $*$-product becomes the pointwise product.
\end{proof}

\section{Formulations of the $PS$-functional calculus via the $\Pi S$-resolvent operators}\label{PSSEC}
In the sequel, we will consider a Banach space $V$ over
$\mathbb{R}$
 with norm $\|\cdot \|$.
It is possible to endow $V$
with an operation of multiplication by elements of $\rr_n$ which gives
a two-sided module over $\rr_n$.
A two-sided module $V$ over $\rr_n$ is called a Banach module over $\rr_n$,
 if there exists a constant $C \geq 1$  such
that $\|va\|\leq C\| v\| |a|$ and $\|av\|\leq C |a|\| v\|$ for all
$v\in V$ and $a\in\rr_n$.
 By $V_n$ we denote $V\otimes \rr_n$ that turns out to be a
 two-sided Banach module over $\rr_n$.
 An element in $V_n$ is of the type $\sum_A v_A\otimes e_A$ (where
 $A=i_1\ldots i_r$, $i_\ell\in \{1,2,\ldots, n\}$, $i_1<\ldots <i_r$ is a multi-index).
The multiplications of an element $v\in V_n$ with a scalar
$a\in \rr_n$ are defined by $va=\sum_A v_A \otimes (e_A a)$ and $av=\sum_A v_A \otimes (ae_A )$.
For simplicity, we will write
$\sum_A v_A e_A$ instead of $\sum_A v_A \otimes e_A$. Finally, we define $\| v\|^2_{V_n}=
\sum_A\| v_A\|^2_V$.

We denote by
$\mathcal{B}(V)$  the space
of bounded $\mathbb{R}$-homomorphisms of the Banach space $V$ to itself
 endowed with the natural norm denoted by $\|\cdot\|_{\mathcal{B}(V)}$.
Given $T_A\in \mathcal{B}(V)$, we can introduce the Clifford operator $T=\sum_A e_A T_A$ and
its action on $v=\sum v_Be_B\in V_n$ as $T(v)=\sum_{A,B}
T_A(v_B)e_Ae_B$. The operator $T$ is a right-module homomorphism which is a bounded linear
map on $V_n$.
\\
In the sequel, we will consider an important subclass of Clifford operators, the ones of the form
$T=T_0+\sum_{j=1}^ne_jT_j$ are called paravector operators,
 where $T_j\in\mathcal{B}(V)$ for $j=0,1,\ldots ,n$.
The subset of paravector operators in ${\mathcal{B}(V_n)}$ will be denoted by $\mathcal{B}^{\small 0,1}(V_n)$.
For Clifford operators $T=\sum_A e_A T_A$ we define
$\|T\|_{\mathcal{B}(V_n)}=\sum_A \|T_A\|_{\mathcal{B}(V)}$
and in particular when $T$ is a paravector operator we have
$\|T\|_{\mathcal{B}^{\small 0,1}(V_n)}=\sum_j \|T_j\|_{\mathcal{B}(V)}$.
Note that, in the sequel, we will omit the subscripts $\mathcal{B}^{\small 0,1}(V_n)$ or $\mathcal{B}(V_n)$ in the norm of an operator. Note also that  $\|TS\|\leq \|T\| \|S\|$, finally we denote by $\mathcal{I}$ the identity operator.
We recall this crucial result which is the heart of the spectral theory on the $S$-spectrum because it shows the notion of $S$-spectrum and of $S$-resolvent operators.
The subset of $\mathcal{B}^{0,1}(V_n)$ that consists of those paravector operators
$T=T_0+\sum_{j=1}^ne_jT_j$ with commuting components $T_0,...,T_n$ will be denoted by $\mathcal{BC}^{0,1}(V_n)$. Finally we will use the notation
$\overline{T}=T_0-\sum_{j=1}^ne_jT_j$.

\medskip
We recall some facts on the $S$-functional calculus that we will use in the sequel.
\begin{theorem}\label{leftright}
Let $T\in\mathcal{B}^{\small 0,1}(V_n)$ and let $s\in \mathbb{R}^{n+1}$ with $\|T\|<|s|$.
\begin{enumerate}[(i)]
\item\label{CSSL}
The left $S$-resolvent series equals
\[
\sum_{m=0}^{+\infty}T^ms^{-m-1} = - (T^2 - 2\Re (s)T +|s|^2\id)^{-1}(T-\overline{s}\id).
\]
\item\label{CSSR}
The right $S$-resolvent series equals
\[
 \sum_{m=0}^{+\infty}s^{-m-1}T^m = - (T-\overline{s}\id)(T^2 - 2\Re(s)T +|s|^2\id)^{-1}.
 \]
\end{enumerate}
\end{theorem}

\begin{definition}\label{sspectrum}
Let $T\in\mathcal{B}^{\small 0,1}(V_n)$. For $s\in\mathbb{R}^{n+1}$, we set
\[
\Q_{s}(T):=T^2 - 2\Re(s)T + |s|^2\id.
\]
We define the $S$-resolvent set $\rho_S(T)$ of $T$ as
\[
\rho_S(T) := \{s\in \mathbb{R}^{n+1}: \Q_{s}(T)  \text{ is invertible in $\mathcal{B}(V_n)$}\}
\]
and we define the $S$-spectrum $\sigma_S(T)$ of $T$ as
\[
 \sigma_S(T) := \mathbb{R}^{n+1}\setminus\rho_S(T).
\]
For $s\in\rho_{S}(T)$, the operator $\Q_{s}(T)^{-1}$
is called the pseudo $S$-resolvent operator  of $T$ at $s$.
\end{definition}
\begin{definition}
Let $T\in\mathcal{B}^{\small 0,1}(V_n)$. For $s\in\rho_S(T)$, we define the {\em left $S$-resolvent operator} as
$$
 S_L^{-1}(s,T) = -\Q_{s}(T)^{-1}(T-\overline{s}\,\id),
 $$
and the {\em right $S$-resolvent operator} as
$$
S_R^{-1}(s,T) = -(T-\overline{s}\id)\Q_{s}(T)^{-1}.
$$
\end{definition}

\begin{lemma} Let $T\in\mathcal{B}^{\small 0,1}(V_n)$.
\label{ResolventRegular}

(I) The left $S$-resolvent operator $S_L^{-1}(s,T)$ is a $\mathcal{B}(V_n)$-valued right-slice mponogenic function of the variable $s$ on $\rho_S(T)$.

(II)  The right $S$-resolvent operator $S_R^{-1}(s,T)$ is a $\mathcal{B}(V_n)$-valued left-slice monogenic function of the variable~$s$ on $\rho_S(T)$.

\end{lemma}

\begin{definition}
For $T\in\mathcal{B}^{\small 0,1}(V_n)$,
we denote by $\mathcal{SM}_L(\sigma_{S}(T))$,
$\mathcal{SM}_R(\sigma_{S}(T))$
and $\mathcal{N}(\sigma_{S}(T))$, 
the set of all  left, right and intrinsic  slice monogenic functions
with $\sigma_{S}(T)\subset U$,  where $U$ is a slice Cauchy domain such that
$\overline{U}\subset\dom(f)$ and $\dom(f)$ is the domain of the function $f$.
\end{definition}

\begin{definition}[The $S$-functional calculus]\label{SCalc}
Let $T\in\mathcal{B}^{\small 0,1}(V_n)$  for
$\mathrm{j}\in\mathbb{S}$  set  $ds_\mathrm{j}=ds (-\mathrm{j})$.
Then we have the formulations of the $S$-functional calculus.
We define
\begin{equation}\label{SCalcL}
f(T) :=\frac{1}{2 \pi}\int_{\partial (U\cap \mathbb{C}_\mathrm{j})} S^{-1}_{L}(s,T)\, ds_\mathrm{j}\,  f(s),
 \ \ \ {\rm for\ all} \ \ f\in\mathcal{SM}_L(\sigma_S(T)),
\end{equation}
and
\begin{equation}\label{SRIGHT}
 f(T):=\frac{1}{2 \pi}\int_{\partial (U\cap \mathbb{C}_\mathrm{j})} F(s)\, ds_\mathrm{j}\,S^{-1}_{R}(s,T),
  \ \ \ {\rm for\ all} \ \ f\in\mathcal{SM}_R(\sigma_S(T)).
 \end{equation}
\end{definition}
\begin{remark}
The definition of the $S$-functional calculus is well posed because
the integrals (\ref{SCalcL}) and (\ref{SRIGHT})
  depend neither on $U$ nor on the imaginary unit $\mathrm{j}\in\mathbb{S}$.
  This is also independent of the fact that the components of the operator $T$ commute or not among themselves.
\end{remark}
The $S$-resolvent equation is useful to prove several properties of the S-functional calculus. So it is natural to ask if it is possible to obtain an analog of the classical resolvent equation
\begin{equation}\label{resolveqclassic}
(\lambda I-A)^{-1}(\mu I-A)^{-1}=\big( (\lambda I-A)^{-1}-(\mu I-A)^{-1}\big)\big( \mu-\lambda\big)^{-1},\ \ \lambda, \mu\in \mathbb{C}\setminus \sigma(A),
\end{equation}
where $A$ is a complex operator on a Banach space.
The generalization  to this non commutative setting,
involves both the left and the right S-resolvent operators and the analogue of the term
$\big( (\lambda I-A)^{-1}-(\mu I-A)^{-1}\big)\big( \mu-\lambda\big)^{-1}$,
 which is the difference of the resolvent operators
$(\lambda I-A)^{-1}-(\mu I-A)^{-1}$ multiplied by the Cauchy kernel $(\mu-\lambda)^{-1}$
for the $S$-functional calculus  becomes the difference of the $S$-resolvent operators
$S_{R}^{-1}(s,T)-S_{L}^{-1}(p,T)$ entangled in a suitable way with the slice monogenic Cauchy kernels. In fact we have:

\begin{theorem}[The $S$-resolvent equation, see \cite{acgs}]\label{SREQa}
Let $T\in\boundOP^{0,1}(V_n)$ and let $s,q\in \rho_S(T)$ with $q\notin[s]$.
Set $
\Q_{s}(q)^{-1}:=(q^2-2\Re(s)q+|s|^2)^{-1}.$
 Then the equation
\begin{equation}\label{SREQ1}
S_R^{-1}(s,T)S_L^{-1}(q,T)=\left[\left(S_R^{-1}(s,T)-S_L^{-1}(q,T)\right)q
-\overline{s}\left(S_R^{-1}(s,T)-S_L^{-1}(q,T)\right)\right]\Q_{s}(q)^{-1}
\end{equation}
holds true. Equivalently, it can also be written as
\begin{equation}\label{SREQ2}
S_R^{-1}(s,T)S_L^{-1}(q,T)=\Q_{q}(s)^{-1}
\left[\left(S_L^{-1}(q,T) - S_R^{-1}(s,T)\right)\overline{q}-s\left(S_L^{-1}(q,T) - S_R^{-1}(s,T
)\right)
 \right].
\end{equation}

\end{theorem}
The $S$-resolvent equation is a consequence of the
left and the right $S$-resolvent equations:
\begin{theorem}
Let $T\in\mathcal{B}^{\small 0,1}(V_n)$ and let $s\in\rho_S(T)$. The left $S$-resolvent operator satisfies the {\em left $S$-resolvent equation}
\begin{equation}\label{LeftSREQ}
S_L^{-1}(s,T)s - TS_L^{-1}(s,T) = \id
\end{equation}
and the right $S$-resolvent operator satisfies the {\em right $S$-resolvent equation}
\begin{equation}\label{RightSREQ}
sS_R^{-1}(s,T) - S_R^{-1}(s,T)T = \id.
\end{equation}
\end{theorem}
We point out that
the equations (\ref{LeftSREQ}) and (\ref{RightSREQ})
 cannot be considered the generalizations of the classical resolvent equation.
 Only the equations in Theorem \ref{SREQa} have the properties of the classical resolvent equation
 (\ref{resolveqclassic}).
 The product rule is a consequence of the $S$-resolvent equation.
 \begin{theorem}[Product rule]\label{ProdRule}
Let $T\in\boundOP^{0,1}(V_n)$ and let   $f\in\mathcal{N}(\sigma_S(T))$ and $g\in  \mathcal{SM}_{L}(\sigma_S(T))$ or let $f\in\mathcal{SM}_{R}(\sigma_S(T))$ and $g\in\mathcal{N}(\sigma_S(T))$.  Then
$$
(f g)(T)=f(T)g(T).
$$
\end{theorem}

In this section we give  the formulations of the poly slice monogenic version of the $S$-functional  based on the poly slice Cauchy formulas. This calculus that will be indicated by
$PS$-functional calculus

\begin{definition}[The $\Pi S$-resolvent operators $\Pi_\ell S_L^{-1}(s,T)$ and $\Pi_\ell S_R^{-1}(s,T)$]
Let
$T\in\mathcal{B}^{\small 0,1}(V_n)$
 and $ s \in \rho_S(T)$.
We define the left poly $S$-resolvent operator (for short $\Pi S$-resolvent operator) of order $\ell+1$, for $\ell\in \mathbb{N}$ as
\[
\begin{split}
\Pi_\ell S^{-1}_{L}(s,T):&
=\frac{1}{\ell!}\sum_{k=0}^\ell \binom{\ell}{k}\overline{T}^kS_L^{-1}(s,T)(-\overline{s})^{\ell-k}
\\
&
=-\frac{1}{\ell!}\sum_{k=0}^\ell \binom{\ell}{k}\overline{T}^k\Q_{s}(T)^{-1}(T-\overline{s}\mathcal{I}) (-\overline{s})^{\ell-k}
,
\end{split}
\]
and the right poly $S$-resolvent operator (for short right $\Pi S$-resolvent operator)
\[
\begin{split}
\Pi_\ell S^{-1}_{R}(s,T):&=\frac{1}{\ell!}\sum_{k=0}^\ell \binom{\ell}{k}
(-\overline{s})^{\ell-k}S_R^{-1}(s,T)\overline{T}^k
\\
&
=-\frac{1}{\ell!}\sum_{k=0}^\ell \binom{\ell}{k}
(-\overline{s})^{\ell-k}(T-\overline{s}\id) \Q_{s}(T)^{-1} \overline{T}^k.
\end{split}
\]
where
$\Q_{s}(T):=T^2 - 2\Re(s)T + |s|^2\id.$
\end{definition}

\begin{lemma}\label{POLYResolventRegular} Let $T\in\mathcal{B}^{\small 0,1}(V_n)$ and let $\ell\in \mathbb{N}$. Then we have:

(I) The  poly left $S$-resolvent $\Pi_\ell S^{-1}_{L}(s,T)$ is a $\mathcal{B}(V_n)$-valued right  poly slice monogenic function of the variable $s$ on $\rho_S(T)$ of order $\ell+1$.

(II) \label{RResHolpoly}  The  poly right $S$-resolvent  $\Pi_\ell S^{-1}_{R}(s,T)$ is a $\mathcal{B}(V_n)$-valued left poly slice  monogenic function of the variable $s$ on $\rho_S(T)$
 of order $\ell+1$.
\end{lemma}
\begin{proof} Consider case (I).
It is a direct consequence of the definition because it is of the form
$$
\Pi_\ell S^{-1}_{L}(s,T)
=\sum_{k=0}^\ell \psi_\ell(s,T)\overline{s}^{\ell-k}
$$
where the components
$$
\psi_\ell(s,T):=\frac{1}{\ell!}\binom{\ell}{k}\overline{T}^kS_L^{-1}(s,T)(-1)^{\ell-k}
$$
are right $\mathcal{B}(V_n)$-valued  slice monogenic function in  the variable $s$ on $\rho_S(T)$, recalling that  $V_n=V\otimes \rr_n$. Apply Theorem \ref{carac3} which still holds for function
$\mathcal{B}(V_n)$-valued we get the statement.
 The order $\ell+1$ is clear from the definition.
Case (II) follows with the same considerations.
\end{proof}

We can define the left and the right $PS$-resolvent equations
  for $T\in\mathcal{BC}^{0,1}(V_n)$, that is when $T$ has commuting components.
 This will have consequences on the product rules.

\begin{definition}
Let $M\in \mathbb{N}$ and let $T\in\mathcal{B}^{\small 0,1}(V_n)$.
We denote by $\mathcal{PS}^M_L(\sigma_{S}(T))$,
$\mathcal{PS}^M_R(\sigma_{S}(T))$
and $\mathcal{PN}^M(\sigma_{S}(T))$
the set of all  left, right and intrinsic poly  slice monogenic functions $F$ or order $M$, respectively,
with $\sigma_{S}(T)\subset U$,  where $U$ is a slice Cauchy domain such that
$\overline{U}\subset\dom(F)$ and $\dom(F)$ is the domain of the function $F$.
\end{definition}
Using the Cauchy formula of poly slice monogenic functions we give the definition of the $PS$-functional calculus.
\begin{definition}[The $PS$-functional calculus (I)]\label{POLYSCalc}
Let $T\in\mathcal{B}^{\small 0,1}(V_n)$, $M\in \mathbb{N}$, for
$\mathrm{j}\in\mathbb{S}$,  set  $ds_\mathrm{j}=ds (-\mathrm{j})$, $\overline{\partial}_\mathrm{j}:=\frac{1}{2}(\partial_u+\mathrm{j}\partial_v)$
let $\Pi_\ell S_L^{-1}(s,T)$ and $\Pi_\ell S_R^{-1}(s,T)$
 be the $\Pi S$-resolvent operators,
 for $\ell=0,...,M-1$.
We define
\begin{equation}\label{POLYSLEFT}
F(T) :=\frac{1}{2 \pi}\int_{\partial (U\cap \mathbb{C}_\mathrm{j})}\sum_{\ell=0}^{M-1} \Pi_\ell S^{-1}_{L}(s,T)\, ds_\mathrm{j}\, \overline{\partial_\mathrm{j}}^\ell F(s),
 \ \ \ {\rm for\ all} \ \ F\in\mathcal{PS}^M_L(\sigma_S(T)),
\end{equation}
and
\begin{equation}\label{POLYSRIGHT}
 F(T):=\frac{1}{2 \pi}\int_{\partial (U\cap \mathbb{C}_\mathrm{j})}\sum_{\ell=0}^{M-1} F(s)\overline{\partial_\mathrm{j}}^\ell\, ds_\mathrm{j}\,\Pi_\ell S^{-1}_{R}(s,T),
  \ \ \ {\rm for\ all} \ \ F\in\mathcal{PS}^M_R(\sigma_S(T)).
 \end{equation}
\end{definition}
The following theorem shows that the definitions of the $PS$-functional calculus are well posed.
\begin{theorem}\label{INDEP}
Let $T\in\mathcal{B}^{\small 0,1}(V_n)$, $M\in \mathbb{N}$, for
$\mathrm{j}\in\mathbb{S}$,  set  $ds_\mathrm{j}=ds (-\mathrm{j})$,
$\overline{\partial}_\mathrm{j}:=\frac{1}{2}(\partial_u+\mathrm{j}\partial_v)$ and
 let $\Pi_\ell S_L^{-1}(s,T)$ and $\Pi_\ell S_R^{-1}(s,T)$
 be the $\Pi S$-resolvent operators,
 for $\ell=0,...,M-1$.
 Then the integrals (\ref{POLYSLEFT}) and (\ref{POLYSRIGHT}) depend
neither on $U$ nor on the imaginary unit $\mathrm{j}\in\mathbb{S}$.
\end{theorem}
\begin{proof}
Recall that we work under Assumption \ref{assumption}.
The independence of the integrals (\ref{POLYSLEFT}) and (\ref{POLYSRIGHT}) from the open set $U$ is standard.  We treat the  case
of $F\in\mathcal{PS}_L(\sigma_{S}(T))$, for functions in  $F\in\mathcal{PS}_R(\sigma_{S}(T))$ the proof is similar with obvious changes.
If $\overline{U'}\not\subset U$, then $O:= U \cap U'$ is a slice Cauchy domain that contains $\sigma_{S}(T)$.
We can hence find a third  slice Cauchy domain $U''$ with $\sigma_{S}(T)\subset U''$ and $\overline{U''} \subset O = U \cap U'$. The above arguments show that the integrals over the boundaries of all three sets agree.

To show the independence of $\mathrm{j}\in \mathbb{S}$ we choose
 two units $\mathrm{i},\mathrm{j}\in\mathbb{S}$ and two slice Cauchy domains $U_q,U_s\subset \dom(F)$ with $\sigma_{S}(T)\subset U_q$ and $\overline{U_q}\subset U_s$. (The subscripts $q$ and $s$ are chosen  to indicate the respective variable of integration in the following computation).
We start from the definition of the $PS$-functional calculus integrating on $\partial (U_s\cap \mathbb{C}_\mathrm{j})$:
$$
F(T) :=\frac{1}{2 \pi}\int_{\partial (U_s\cap \mathbb{C}_\mathrm{j})}\sum_{\ell=0}^{M-1} \Pi_\ell S^{-1}_{L}(s,T)\, ds_\mathrm{j}\, \overline{\partial_\mathrm{j}}^\ell F(s),
$$
where
$
F(s)=\sum_{k=0}^{M-1}\overline{s}^kf_k(s).
$
We recall that
$$
\overline{\partial_\mathrm{j}}^\ell \big(\overline{s}^k f_k(s)\big)=
\frac{k!}{(k-\ell)!}\overline{s}^{k-\ell} f_k(s),  \ \ \  {\rm for} \ \ \  k\geq \ell
$$
and $\overline{\partial_\mathrm{j}}^\ell \big(\overline{s}^k f_k(s)\big)=0$, for $k< \ell$.
So we have that
$$
\overline{\partial_\mathrm{j}}^\ell \Big(\sum_{k=0}^{M-1}\overline{s}^k f_k(s)\Big)
= \sum_{k=0}^{M-1}\overline{\partial_\mathrm{j}}^\ell\Big(\overline{s}^k f_k(s)\Big)
= \sum_{k=\ell}^{M-1}\frac{k!}{(k-\ell)!}\overline{s}^{k-\ell} f_k(s),
$$
since for $k< \ell$ the terms $\overline{\partial_\mathrm{j}}^\ell \Big(\sum_{k=0}^{M-1}\overline{s}^k f_k(s)\Big)$ are zero.
Now consider $F(T)$ written as
$$
F(T) =\frac{1}{2 \pi}\int_{\partial (U_s\cap \mathbb{C}_\mathrm{j})}\sum_{\ell=0}^{M-1} \Pi_\ell S^{-1}_{L}(s,T)\, ds_\mathrm{j}\, \sum_{k=\ell}^{M-1}\frac{k!}{(k-\ell)!}\overline{s}^{k-\ell} f_k(s).
$$
We can write it more explicitly (putting a label $(\ell=1,2,....,M-1)$ in front of the integrals to identify them in the sequel) as
\begin{equation}
\begin{split}
F(T) &=(\ell=0)\frac{1}{2 \pi}\int_{\partial (U_s\cap \mathbb{C}_\mathrm{j})}\Pi_0 S^{-1}_{L}(s,T)\, ds_\mathrm{j}\, \sum_{k=0}^{M-1}\overline{s}^{k} f_k(s),
\\
&
+(\ell=1)\frac{1}{2 \pi}\int_{\partial (U_s\cap \mathbb{C}_\mathrm{j})}\Pi_1 S^{-1}_{L}(s,T)\, ds_\mathrm{j}\, \sum_{k=1}^{M-1}\frac{k!}{(k-1)!}\overline{s}^{k-1} f_k(s),
\\
&
+(\ell=2)\frac{1}{2 \pi}\int_{\partial (U_s\cap \mathbb{C}_\mathrm{j})}
\Pi_2 S^{-1}_{L}(s,T)\, ds_\mathrm{j}\, \sum_{k=2}^{M-1}\frac{k!}{(k-2)!}\overline{s}^{k-2} f_k(s),
\\
& \ldots
\\
&
+(\ell=M-1)\frac{1}{2 \pi}\int_{\partial (U_s\cap \mathbb{C}_\mathrm{j})}
\Pi_{M-1} S^{-1}_{L}(s,T)\, ds_\mathrm{j}\,(M-1)! f_{M-1}(s).
\end{split}
\end{equation}

Now we replace the explicit expressions of the $\Pi_\ell S$-resolvent operators
\begin{equation}\label{EXPRES}
\begin{split}
&
\Pi_0 S^{-1}_{L}(s,T)=S_L^{-1}(s,T),
 \\
 &
\Pi_1 S^{-1}_{L}(s,T) =S^{-1}_{L}(s,T)(-\overline{s}) +\overline{T}S_L^{-1}(s,T),
\\
&
\Pi_2 S^{-1}_{L}(s,T)=\frac{1}{2!}\Big(S_L^{-1}(s,T)(-\overline{s})^{2}
+ 2 \overline{T}^1S_L^{-1}(s,T)(-\overline{s})
+\overline{T}^2S_L^{-1}(s,T)\Big),
\\
&
\ldots
\\
&
\Pi_{M-1} S^{-1}_{L}(s,T)
=\frac{1}{(M-1)!}\sum_{k=0}^{M-1} \binom{M-1}{k}\overline{T}^kS_L^{-1}(s,T)(-\overline{s})^{M-1-k}.
\end{split}
\end{equation}
and we get
\[
\begin{split}
F(T) &=(\ell=0)\frac{1}{2 \pi}\int_{\partial (U_s\cap \mathbb{C}_\mathrm{j})} S^{-1}_{L}(s,T)\, ds_\mathrm{j}\, \sum_{k=0}^{M-1}\overline{s}^{k} f_k(s),
\\
&
+(\ell=1)\frac{1}{2 \pi}\int_{\partial (U_s\cap \mathbb{C}_\mathrm{j})}
\Big(S^{-1}_{L}(s,T)(-\overline{s}) +\overline{T}S_L^{-1}(s,T)\Big)\, ds_\mathrm{j}\, \sum_{k=1}^{M-1}\frac{k!}{(k-1)!}\overline{s}^{k-1} f_k(s),
\\
&
+(\ell=2)\frac{1}{2 \pi}\int_{\partial (U_s\cap \mathbb{C}_\mathrm{j})}
\frac{1}{2!}\Big(S_L^{-1}(s,T)(-\overline{s})^{2}
+ 2 \overline{T}^1S_L^{-1}(s,T)(-\overline{s})
+\overline{T}^2S_L^{-1}(s,T)\Big)\, ds_\mathrm{j}\, \\
&
 \ \ \ \ \times\sum_{k=2}^{M-1}\frac{k!}{(k-2)!}\overline{s}^{k-2} f_k(s)
\\
& \ldots
\\
&
+(\ell=M-1)\frac{1}{2 \pi}\int_{\partial (U_s\cap \mathbb{C}_\mathrm{j})}
\Big( \frac{1}{(M-1)!}\sum_{k=0}^{M-1} \binom{M-1}{k}\overline{T}^kS_L^{-1}(s,T)(-\overline{s})^{M-1-k}\Big)\, ds_\mathrm{j}\,
\\
&
\ \ \ \ \times (M-1)!  f_{M-1}(s).
\end{split}
\]
 Now in the integral for $(\ell=0)$
 $$
 J_{(\ell=0)}:=\frac{1}{2 \pi}\int_{\partial (U_s\cap \mathbb{C}_\mathrm{j})} S^{-1}_{L}(s,T)\, ds_\mathrm{j}\, \sum_{k=0}^{M-1}\overline{s}^{k} f_k(s)
$$
we separate
the term with $f_0(s)$ and we write it as
$$
 J_{(\ell=0)}:=R_0-\frac{1}{2 \pi}\int_{\partial (U_s\cap \mathbb{C}_\mathrm{j})} S^{-1}_{L}(s,T)\, ds_\mathrm{j}\, f_0(s)
$$
where $R_0$ contains all the other terms. The terms for $(\ell=1)$
are separated as
$$
 J_{(\ell=1)}:= R_1+ \frac{1}{2 \pi}\int_{\partial (U_s\cap \mathbb{C}_\mathrm{j})}
\overline{T}S_L^{-1}(s,T)\, ds_\mathrm{j}\,  f_1(s)
$$
where $R_1$ contains all the other terms. We proceed in the same manner also for $(\ell=2)$
and for the rest of the terms to get
$$
 J_{(\ell=2)}:= R_2+ \frac{1}{2 \pi}\int_{\partial (U_s\cap \mathbb{C}_\mathrm{j})}
\overline{T}^2S_L^{-1}(s,T)\, ds_\mathrm{j}\,  f_2(s)
$$
$$
\ldots
$$
$$
 J_{(\ell=M-1)}:= R_{M-1}+ \frac{1}{2 \pi}\int_{\partial (U_s\cap \mathbb{C}_\mathrm{j})}
\overline{T}^{M-1}S_L^{-1}(s,T)\, ds_\mathrm{j}\,  f_{M-1}(s).
$$
Finally consider the sum of the  terms in $\sum_{j=0}^{M-1}R_j$.
It turns out to be zero and this can be seen by gathering in the
sum $\sum_{j=0}^{M-1}R_j$ the terms that contain $f_0$, $f_1$ .... $f_{N-1}$.
In each of these sums have the addends that cancel.
This shows that we are left with
$$
F(T)=\sum_{k=0}^{M-1}\overline{T}^k
\frac{1}{2 \pi}\int_{\partial (U_s\cap \mathbb{C}_\mathrm{j})} S^{-1}_{L}(s,T)\, ds_\mathrm{j}\, f_k(s).
$$
Repeating the above computation with the imaginary  unit $\mathrm{i}\in \mathbb{S}$ on
$\partial (U_q\cap \mathbb{C}_\mathrm{i})$ we obtain
$$
F(T)=\sum_{k=0}^{M-1}\overline{T}^k
\frac{1}{2 \pi}\int_{\partial (U_q\cap \mathbb{C}_\mathrm{i})} S^{-1}_{L}(q,T)\, dq_\mathrm{i}\, f_k(q)
$$
but since the $S$-functional calculus is independent from the imaginary units $\mathrm{i}$ and $\mathrm{j}$ in $\mathbb{S}$ we get the statement.

\end{proof}

\begin{remark}[The case $M=3$]
For the reader's convenience, in order
to understand in which way the terms in the proof of the previous result combine,
we write explicitly the case $M=3$.
We recall that the explicit expressions of the $\Pi_\ell S$-resolvent operators are given in (\ref{EXPRES}) and the function $F(s)$ is of the form
$$
F(s)=\overline{s}^2f_2(s)+\overline{s}f_1(s)+f_0(s).
$$
So we have
$$
F(T) :=\frac{1}{2 \pi}\int_{\partial (U_s\cap \mathbb{C}_\mathrm{j})}
\sum_{\ell=0}^{2} \Pi_\ell S^{-1}_{L}(s,T)\, ds_\mathrm{j}\, \overline{\partial_\mathrm{j}}^\ell
\Big(\overline{s}^2f_2(s)+\overline{s}f_1(s)+f_0(s)\Big)
$$
and we split in 3 terms
\[
\begin{split}
F(T) &= (\ell=0)
\frac{1}{2 \pi}\int_{\partial (U_s\cap \mathbb{C}_\mathrm{j})}
 \Pi_0 S^{-1}_{L}(s,T)\, ds_\mathrm{j}\,
\Big(\overline{s}^2f_2(s)+\overline{s}f_1(s)+f_0(s)\Big)
\\
&
+(\ell=1)
\frac{1}{2 \pi}\int_{\partial (U_s\cap \mathbb{C}_\mathrm{j})}
 \Pi_1 S^{-1}_{L}(s,T)\, ds_\mathrm{j}\,
\Big(2\overline{s}f_2(s)+f_1(s)\Big)
\\
&+(\ell=2)
\frac{1}{2 \pi}\int_{\partial (U_s\cap \mathbb{C}_\mathrm{j})}
 \Pi_2 S^{-1}_{L}(s,T)\, ds_\mathrm{j}\,
\Big(2f_2(s)\Big)
\end{split}
\]
and, replacing the $\Pi_\ell S$-resolvent operators, the addends become
\[
\begin{split}
J_{(\ell=0)}:&=
\frac{1}{2 \pi}\int_{\partial (U_s\cap \mathbb{C}_\mathrm{j})}
 S^{-1}_{L}(s,T)\, ds_\mathrm{j}\,
\Big(\overline{s}^2f_2(s)+\overline{s}f_1(s)\Big)
\\
&
+\frac{1}{2 \pi}\int_{\partial (U_s\cap \mathbb{C}_\mathrm{j})}
 S^{-1}_{L}(s,T)\, ds_\mathrm{j}\,
\Big(f_0(s)\Big)
\end{split}
\]
where the last integral gives $f_0(T)$, thanks to the $S$-functional calculus. Now consider
\[
\begin{split}
J_{(\ell=1)}:&=
\frac{1}{2 \pi}\int_{\partial (U_s\cap \mathbb{C}_\mathrm{j})}
 \Big(S^{-1}_{L}(s,T)(-\overline{s}) \Big)\, ds_\mathrm{j}\,
\Big(2\overline{s}f_2(s)\Big)
\\
&
+
\frac{1}{2 \pi}\int_{\partial (U_s\cap \mathbb{C}_\mathrm{j})}
 \Big(S^{-1}_{L}(s,T)(-\overline{s}) \Big)\, ds_\mathrm{j}\,
\Big(f_1(s)\Big)
\\
&
+
\frac{1}{2 \pi}\int_{\partial (U_s\cap \mathbb{C}_\mathrm{j})}
 \Big(\overline{T}S_L^{-1}(s,T)\Big)\, ds_\mathrm{j}\,
\Big(2\overline{s}f_2(s)\Big)
\\
&
+
\frac{1}{2 \pi}\int_{\partial (U_s\cap \mathbb{C}_\mathrm{j})}
 \Big(\overline{T}S_L^{-1}(s,T)\Big)\, ds_\mathrm{j}\,
  \Big(f_1(s)\Big),
\end{split}
\]
where the last integral gives $\overline{T}f_1(T)$. For $J_{(\ell=2)}$ we have

\[
\begin{split}
&
J_{(\ell=2)}:=
\frac{1}{2 \pi}\int_{\partial (U_s\cap \mathbb{C}_\mathrm{j})}
 \Big(\frac{1}{2!}S_L^{-1}(s,T)(-\overline{s})^{2}
\Big)\, ds_\mathrm{j}\,
\Big(2 f_2(s)\Big)
\\
&
+
\frac{1}{2 \pi}\int_{\partial (U_s\cap \mathbb{C}_\mathrm{j})}
 \Big(\frac{1}{2!} \cdot 2 \overline{T}S_L^{-1}(s,T)(-\overline{s})
\Big)\, ds_\mathrm{j}\,
\Big(2 f_2(s)\Big)
\\
&
+
\frac{1}{2 \pi}\int_{\partial (U_s\cap \mathbb{C}_\mathrm{j})}
 \Big(\frac{1}{2!}\overline{T}^2S_L^{-1}(s,T)\Big)\, ds_\mathrm{j}\,
\Big(2f_2(s)\Big),
\end{split}
\]
where the last integral give the term $\overline{T}^2f_2(T)$.
Finally consider all the terms that remain.
Observe that all the terms in the sum with $f_1$ cancel and
the terms that contain $f_2$ cancel out as well.
\end{remark}

\begin{remark}
To prove that the integrals (\ref{POLYSLEFT}) and (\ref{POLYSRIGHT}) do not depend
 on the imaginary unit $\mathrm{j}\in\mathbb{S}$
we can also use the poly slice monogenic Cauchy formula.
This strategy to show that the $PS$-functional calculus is well posed
 is more similar to the one used for the $S$-functional calculus.
 We just give the hints because the computations are longer with respect to the proof that we have given above.
 With the same notations on the domains as in the proof above we observe the following: the set $U_q^c := \mathbb{R}^{n+1} \setminus U_q$  is an unbounded axially symmetric domain with $\overline{U_q^c}\subset \rho_{S}(T)$. The left $PS$-resolvent  operator is right poly slice monogenic on $\rho_{S}(T)$ and  at infinity, since
\[
\lim_{s\to\infty}  S_L^{-1}(s,T) = \lim_{s\to\infty}
\sum_{n=0}^{+\infty}T^n s^{-n-1} = 0
\]
and the left $S$-resolvent operator
 is slice monogenic.
Now observe that the $\Pi S$-resolvent operator
$$
\Pi_\ell S^{-1}_{L}(s,T)
=\frac{1}{\ell!}\sum_{k=0}^\ell \binom{\ell}{k}\overline{T}^kS_L^{-1}(s,T)(-\overline{s})^{\ell-k}
$$
has operator valued slice monogenic components, and they go to zero, i.e.,
$$
\frac{(-1)^{\ell-k}}{\ell!}\binom{\ell}{k}\overline{T}^kS_L^{-1}(s,T)\to 0, \ \ {\rm for} \ \ s\to \infty.
$$
So we can represent it with the Cauchy formula for unbounded slice Cauchy domains
and it is right poly slice monogenic on the $S$-resolvent set $\rho_S(T)$. We have
\[
\Pi_\ell S^{-1}_{L}(s,T) =\frac{1}{2 \pi}\int_{\partial (U_q^c\cap \mathbb{C}_\mathrm{i})}\sum_{m=0}^{M-1}
 \Pi_\ell S^{-1}_{L}(q,T)\overline{\partial_\mathrm{i}}^m\, dq_\mathrm{i}\, \Pi_m S^{-1}_{R}(q,s)
\]
because $\Pi_\ell S^{-1}_{L}(s,T)$ is poly slice monogenic up to order $M\geq\ell$,
for any $s\in U_q^c$. From $S_R^{-1}(q,s) = - S_L^{-1}(s,q)$ we find that the Cauchy kernel
$$
\Pi_m S^{-1}_{R}(s,q)=\frac{1}{m!}\sum_{k=0}^m \binom{m}{k}
(-\overline{s})^{m-k}S_R^{-1}(s,q)\overline{q}^k
$$
can be written as
$$
\Pi_m S^{-1}_{R}(q,s)=\frac{1}{m!}\sum_{k=0}^m \binom{m}{k}
(-\overline{q})^{m-k}S_R^{-1}(q,s)\overline{s}^k
$$
$$
=-\frac{1}{m!}\sum_{k=0}^\ell \binom{m}{k}
(-\overline{q})^{m-k}S_L^{-1}(s,q)\overline{s}^k.
$$
So keeping in mind that $\partial(U_q^c\cap\cc_\mathrm{i}) = - \partial (U_q\cap\cc_\mathrm{i})$ we have
the chain of equalities
\[
\begin{split}
&F(T) :=\frac{1}{2 \pi}\int_{\partial (U_s\cap \mathbb{C}_\mathrm{j})}\sum_{\ell=0}^{M-1} \Pi_\ell S^{-1}_{L}(s,T)\, ds_\mathrm{j}\, \overline{\partial_\mathrm{j}}^\ell F(s)
\\
&
=\frac{1}{2 \pi}\int_{\partial (U_s\cap \mathbb{C}_\mathrm{j})}\sum_{\ell=0}^{M-1}
\Big( \frac{1}{2 \pi}\int_{\partial (U_q^c\cap \mathbb{C}_\mathrm{i})}
\sum_{m=0}^{M-1} \Pi_\ell S^{-1}_{L}(q,T)\overline{\partial_\mathrm{i}}^m \, dq_\mathrm{i}\,
\Pi_m S^{-1}_{R}(q,s) \Big)\, ds_\mathrm{j}\, \overline{\partial_\mathrm{j}}^\ell F(s)
\\
&
=
\frac{1}{2 \pi}\int_{\partial (U_s\cap \mathbb{C}_\mathrm{j})}\sum_{\ell=0}^{M-1}
\Big( \frac{1}{2 \pi}\int_{\partial (U_q\cap \mathbb{C}_\mathrm{i})}
\sum_{m=0}^{M-1}  \Pi_\ell S^{-1}_{L}(q,T)\,\overline{\partial_\mathrm{i}}^m dq_\mathrm{i}\, \Pi_m S^{-1}_{L}(s,q) \Big)\, ds_\mathrm{j}\, \overline{\partial_\mathrm{j}}^\ell F(s)
\\
&
=
\frac{1}{2 \pi}\int_{\partial (U_q\cap \mathbb{C}_\mathrm{i})}\sum_{m=0}^{M-1}
  \Pi_\ell S^{-1}_{L}(q,T)\,\overline{\partial_\mathrm{i}}^m  dq_\mathrm{i}\,
  \Big(
\frac{1}{2 \pi}\int_{\partial (U_s\cap \mathbb{C}_\mathrm{j})}
\sum_{\ell=0}^{M-1}\Pi_m S^{-1}_{L}(s,q) \, ds_\mathrm{j}\, \overline{\partial_\mathrm{j}}^\ell F(s)\Big)
\end{split}
\]
where we have used Fubini's theorem  and with some very long computations we obtain
$$
F(T) =\frac{1}{2 \pi}\int_{\partial (U_q\cap \mathbb{C}_\mathrm{i})}\sum_{m=0}^{M-1}
  \Pi_m S^{-1}_{L}(q,T)\,\overline{\partial_\mathrm{i}}^m dq_\mathrm{i}\,
   F(q)
  =\frac{1}{2 \pi}\int_{\partial (U_q\cap \mathbb{C}_\mathrm{i})}\sum_{m=0}^{M-1}
  \Pi_m S^{-1}_{L}(q,T)\, dq_\mathrm{i}\,
  \overline{\partial_\mathrm{i}}^m F(q).
$$
So we conclude that the poly slice monogenic Cauchy formula gives the independence on the imaginary
units in the sphere $\mathbb{S}$
because we chose $\overline{U_q}\subset U_s$.
The computations above to show all the cancellations
 are slightly more complicated with respect to the one in the prove that we have given above.

\end{remark}

We start by proving some basic results related to the $PS$-functional calculus introduced in Definition \ref{POLYSCalc}
\begin{proposition}[Linearity of the $PS$-functional calculus]
Let $T\in\mathcal{B}^{\small 0,1}(V_n)$ and $M\geq 1$.  Then, we have

\begin{enumerate}
\item If  $F,G\in\mathcal{PS}^M_L(\sigma_S(T))$, then we have:

$$(F+G)(T)=F(T)+G(T), \quad \textbf{    }  (F \lambda)(T)=F(T) \lambda, \ \  \text{   for all  } \lambda\in\mathbb{R}_n.$$

\item If  $F,G\in\mathcal{PS}^M_R(\sigma_S(T))$, then we have:

$$(F+G)(T)=F(T)+G(T), \quad \textbf{    }  (\lambda F)(T)=\lambda F(T), \ \ \ \text{   for all  } \lambda\in\mathbb{R}_n.$$
\end{enumerate}

\end{proposition}
\begin{proof}
This follows directly by construction of the PS-functional calculus in Definition \ref{POLYSCalc}.
\end{proof}

\begin{remark}\label{RMSER}
{\rm
So far, we always made use of the poly decomposition (\ref{strtpolyleft}) of a function
  $F\in \mathcal{PS}^M_L(U)$ in terms of $f_0,...,f_{M-1}\in \mathcal{SM}_L(U)$.
When the components are polynomials or power series
$f_\ell(x)$ for $\ell=0,...,M-1$, then
$$
f_\ell(x)=\sum_{u=0}^\infty x^u A_{\ell,u}, \ \ \ {\rm with}\ \ \ \  \  A_{\ell,u}\in \mathbb{R}_n,\ \ \ {\rm for} \  \ u\in  \mathbb{N}_0, \ \ \ \  \ell=0,...,M-1.
$$
So, for $f_{M-1}\not=0$, we have the decomposition
\begin{equation}\label{poweserpolyleft}
F(x)=\sum_{\ell=0}^{M-1} \overline{x}^\ell \sum_{u=0}^\infty x^u A_{\ell,u}, \ \ \ \ \forall x\in U,
\end{equation}
where $U$ is the set of convergence of all the series $f_\ell(x)$, for  $\ell=0,...,M-1$.
Analogous considerations hold for the right case (\ref{strtpolyright}), and we have
\begin{equation}\label{poweserpolyright}
F(x)=\sum_{\ell=0}^{M-1}  \Big(\sum_{u=0}^\infty A_{\ell,u} x^u \Big) \overline{x}^\ell, \ \ \ \ \forall x\in U.
\end{equation}
}
\end{remark}

\medskip
The following theorem shows the compatibility of the $PS$-functional calculus with the
poly slice monogenic polynomials and series with respect to the slice monogenic components. That is when we consider, for example, the expansion
(\ref{poweserpolyleft}),
the $PS$-functional calculus gives
$$
F(T)=\sum_{\ell=0}^{M-1} \overline{T}^\ell \sum_{u=0}^\infty T^u A_{\ell,u},
$$
when $F$ is defined on the $S$-spectrum of $T$.
\begin{theorem}\label{ywdejtybx}
Let $T\in\mathcal{B}^{\small 0,1}(V_n)$ where we set
 $\overline{T}=T_0-\sum_{j=1}^ne_jT_j$.
 Let $M\in \mathbb{N}$, for
$\mathrm{j}\in\mathbb{S}$,  set  $ds_\mathrm{j}=ds (-\mathrm{j})$  and  $\overline{\partial}_\mathrm{j}:=\frac{1}{2}(\partial_u+\mathrm{j}\partial_v)$.
Let $\Pi_\ell S_L^{-1}(s,T)$ and $\Pi_\ell S_R^{-1}(s,T)$
 be the $\Pi S$-resolvent operators,
 for $\ell=0,...,M-1$,
and suppose  $\sigma_S(T) \subset U$, where $U$ is a slice Cauchy domain.

(I) Then when $F$ is the series expansion in (\ref{poweserpolyleft})
 we have
$$
\frac{1}{2 \pi}\int_{\partial (U\cap \mathbb{C}_\mathrm{j})}\sum_{\ell=0}^{M-1} \Pi_\ell S^{-1}_{L}(s,T)\, ds_\mathrm{j}\, \overline{\partial_\mathrm{j}}^\ell F(s)
=\sum_{\ell=0}^{M-1} \overline{T}^\ell \sum_{u=0}^\infty T^u A_{k,u}.
$$

(II)
 Then when $F$ is the series expansion in (\ref{poweserpolyright})
 we have
$$
\frac{1}{2 \pi}\int_{\partial (U\cap \mathbb{C}_\mathrm{j})}\sum_{\ell=0}^{M-1} F(s)\overline{\partial_\mathrm{j}}^\ell\, ds_\mathrm{j}\,\Pi_\ell S^{-1}_{R}(s,T)=
\sum_{\ell=0}^{M-1}  \Big(\sum_{u=0}^\infty A_{k,u} T^u \Big) \overline{T}^\ell.
$$
\end{theorem}
\begin{proof}
It follows from the definition with some computations.
\end{proof}
\begin{remark}
From Theorem \ref{ywdejtybx} we see a first difference with respect to the $S$-functional calculus for intrinsic functions. In fact, for the $PS$-functional calculus it is not enough that
$F$ is an intrinsic function to have that
$$
\sum_{\ell=0}^{M-1} \overline{T}^\ell \sum_{u=0}^\infty T^u A_{k,u}
=\sum_{\ell=0}^{M-1}  \Big(\sum_{u=0}^\infty A_{k,u} T^u \Big) \overline{T}^\ell
$$
but we have to require that $T$ has commuting components, i.e. $T\in\mathcal{BC}^{\small 0,1}(V_n)$.
\end{remark}
The following theorem is important to prove the product rule.
\begin{theorem}[The intrinsic $PS$-functional calculus (I)]\label{PSfuncintrin}
Let $T\in\mathcal{BC}^{\small 0,1}(V_n)$ and set  $ds_\mathrm{j}=ds (-\mathrm{j})$
and  $\overline{\partial}_\mathrm{j}:=\frac{1}{2}(\partial_u+\mathrm{j}\partial_v)$. Let   $F\in\mathcal{PN}^M(\sigma_{S}(T))$.
Let $\Pi_\ell S_L^{-1}(s,T)$ and $\Pi_\ell S_R^{-1}(s,T)$
 be the $\Pi S$-resolvent operators,
 for $\ell=0,...,M-1$,
and suppose  $\sigma_S(T) \subset U$, where $U$ is a slice Cauchy domain.
Then the left and the right formulations of the $PS$-functional calculus define the same operators; i.e., we have
\begin{equation}\label{UGUALIPOLYSCalcL}
\begin{split}
F(T)&=\frac{1}{2 \pi}\int_{\partial (U\cap \mathbb{C}_\mathrm{j})}\sum_{\ell=0}^{M-1} \Pi_\ell S^{-1}_{L}(s,T)\, ds_\mathrm{j}\, \overline{\partial_\mathrm{j}}^\ell F(s)
\\
&
=\frac{1}{2 \pi}\int_{\partial (U\cap \mathbb{C}_\mathrm{j})}\sum_{\ell=0}^{M-1} F(s)\overline{\partial_\mathrm{j}}^\ell\, ds_\mathrm{j}\,\Pi_\ell S^{-1}_{R}(s,T).
\end{split}
\end{equation}
\end{theorem}
\begin{proof}
It is a consequence of the Theorem \ref{smpRunge}, when we approximate the intrinsic function $F$ by a sequence of rational functions $F_n$
we have
\[
\begin{split}
 \frac{1}{2 \pi}\int_{\partial (U\cap \mathbb{C}_\mathrm{j})}\sum_{\ell=0}^{M-1} \Pi_\ell S^{-1}_{L}(s,T)\, ds_\mathrm{j}\, \overline{\partial_\mathrm{j}}^\ell F_n(s)
=\frac{1}{2 \pi}\int_{\partial (U\cap \mathbb{C}_\mathrm{j})}\sum_{\ell=0}^{M-1} F_n(s)\overline{\partial_\mathrm{j}}^\ell\, ds_\mathrm{j}\,\Pi_\ell S^{-1}_{R}(s,T)
\end{split}
\]
for $T\in\mathcal{BC}^{\small 0,1}(V_n)$  and from the continuity of the $PS$-functional calculus we get the statement.
\end{proof}

\section{Formulations of the $PS$-functional calculus via the modified $S$-resolvent operators}\label{modofresr}
There is an alternative way to define the $PS$-functional calculus:
instead of using the Cauchy formula for poly slice monogenic functions, we can use
the integral representation
of slice monogenic function.
In order to apply this strategy we need to define the modified $S$-resolvent operators and their related modified $S$-resolvent equation.
As we have seen the $PS$-resolvent operators are poly slice monogenic, while the
modified $S$-resolvent operators are slice monogenic.

\medskip
We define the modified $S$-resolvent operator series.
\begin{definition}
Let $T\in\mathcal{B}^{\small 0,1}(V_n)$ and let $s\in \mathbb{R}^{n+1}$ with $\|T\|<|s|$ and
$B \in\boundOP(V_n)$. We call
$$
\sum_{m=0}^{+\infty}T^m Bs^{-m-1}, \ \ \ {\rm and} \ \ \ \sum_{m=0}^{+\infty}s^{-m-1}B T^m
$$
the modified left (resp. right) $S$-resolvent operator series expansions.
\end{definition}
\begin{theorem}\label{leftrightmod}
Let $T\in\mathcal{B}^{\small 0,1}(V_n)$ and let $s\in \mathbb{R}^{n+1}$ with $\|T\|<|s|$ and
let $B \in\boundOP(V_n)$.
\begin{enumerate}[(i)]
\item\label{CSSLmod}
The modified left $S$-resolvent operator series equals
\[
\sum_{m=0}^{+\infty}T^m B s^{-m-1} = - (T^2 - 2\Re (s)T +|s|^2\id)^{-1}(TB-B\overline{s}).
\]
\item\label{CSSRmod}
The modified right $S$-resolvent operator  series equals
\[
 \sum_{m=0}^{+\infty}s^{-m-1}BT^m = - (BT-\overline{s}B)(T^2 - 2\Re(s)T +|s|^2\id)^{-1}.
 \]
\end{enumerate}
\end{theorem}
\begin{proof}
It follows directly from the relations
\begin{equation}\label{SERLEFT}
\sum_{m = 0}^{+\infty} q^m B s^{-1-m}=-(q^2-2\Re(s) q+|s|^2)^{-1}(qB-B\overline{s}),
\end{equation}
and
\begin{equation}\label{SERRIGHT}
\sum_{m = 0}^{+\infty} s^{-1-m} Bq^m =-(Bq-\overline{s}B)(q^2-2\Re(s) q+|s|^2)^{-1},
\end{equation}
which hold
for $s,q \in \mathbb{R}^{n+1}$ with $|q|< |s|$, replacing $q$ by $T$ and $\|T\|< |s|$.
See \cite{acgs} for more details.
\end{proof}

\begin{definition}[The modified $S$-resolvent operators]
Let $T\in\mathcal{B}^{\small 0,1}(V_n)$ and
let $B \in\boundOP(V_n)$. For $s\in\rho_S(T)$, we define the {\em modified  left $S$-resolvent operator} as
$$
 S_L^{-1}(s,T;B):= -\Q_{s}(T)^{-1}(TB-B\overline{s}),
 $$
and the {\em modified  right $S$-resolvent operator} as
$$
S_R^{-1}(s,T;B) := -(BT-\overline{s}B)\Q_{s}(T)^{-1},
$$
where
$
\Q_{s}(T):=T^2 - 2\Re(s)T + |s|^2\id.
$
\end{definition}
\begin{lemma}\label{holmodres} Let $T\in\mathcal{B}^{\small 0,1}(V_n)$ and
let $B \in\boundOP(V_n)$.

(I) The modified left $S$-resolvent operator $S_L^{-1}(s,T;B)$ is a $\mathcal{B}(V_n)$-valued right-slice monogenic function of the variable $s$ on $\rho_S(T)$.

(II)  The modified right $S$-resolvent operator $S_R^{-1}(s,T;B)$ is a $\mathcal{B}(V_n)$-valued left-slice monogenic function of the variable~$s$ on $\rho_S(T)$.

\end{lemma}
\begin{proof}
It is can be proved by a direct computation as in the case of the $S$-resolvent operators.
\end{proof}

Some more interesting relation can be obtained when  $B \in\boundOP(V_n)$ commutes with $T$ and this fact will be crucial in the sequel.
\begin{lemma}\label{Scoomution}
Let $T\in\mathcal{B}^{\small 0,1}(V_n)$ and suppose that $B \in\boundOP(V_n)$ commutes with $T$.
Let $s\in\rho_S(T)$.
Then we have
$$
 S_L^{-1}(s,T;B)= BS_L^{-1}(s,T),
 $$
and
 $$
  S_R^{-1}(s,T;B) = S_R^{-1}(s,T)B.
$$
\end{lemma}
\begin{proof}
Since $TB=BT$ the statement follows from the fact that $(TB-B\overline{s})=B(T-\overline{s}\mathcal{I})$ and $\Q_{s}(T)^{-1}B=B\Q_{s}(T)^{-1}$.
\end{proof}
\begin{theorem}
Let $T\in\mathcal{B}^{\small 0,1}(V_n)$ and suppose that $B \in\boundOP(V_n)$ commutes with $T$.
Let $s\in\rho_S(T)$. The modified left $S$-resolvent operator satisfies the {\em modified left $S$-resolvent equation}
\begin{equation}\label{LeftSREQmod}
S_L^{-1}(s,T;B)s - TS_L^{-1}(s,T;B)=BS_L^{-1}(s,T)s - TBS_L^{-1}(s,T) = B
\end{equation}
and the modified right $S$-resolvent operator satisfies the {\em modified right $S$-resolvent equation}
\begin{equation}\label{RightSREQmod}
sS_R^{-1}(s,T;B) -S_R^{-1}(s,T;B)T =sS_R^{-1}(s,T)B - S_R^{-1}(s,T)BT = B.
\end{equation}
\end{theorem}
\begin{proof}
Since $2\Re(s)$ and $|s|^2$ are real, they commute with the operator $T$. Therefore
\[
 T\Q_{s}(T) = \Q_{s}(T)T
\]
and in turn
\[
 \Q_{s}(T)^{-1}T = T\Q_{s}(T)^{-1}.
 \]
Thus we have
\[
\begin{split}
S_L^{-1}(s,T;B)s - TS_L^{-1}(s,T;B) =
& - \Q_{s}(T)^{-1}(TB-B \overline{s})s +  T\Q_{s}(T)^{-1}(TB-B\overline{s})\\
=&\Q_{s}(T)^{-1}\left( - (TB-B\overline{s})s +  T(TB-B\overline{s}) \right)\\
=& \Q_{s}(T)^{-1}B\Q_{s}(T) = B,
\end{split}
\]
where we have used the fact that $B$ and $T$ commute.
The modified right $S$-resolvent equation follows by similar computations.
\end{proof}

\begin{theorem}[The modified $S$-resolvent equation]\label{SREQ}\index{$S$-resolvent equation}
Let $T\in\mathcal{B}^{\small 0,1}(V_n)$ and suppose that $B \in\boundOP(V_n)$ commutes with $T$ and let $s,q\in \rho_S(T)$ with $q\notin[s]$. Then the equation
\begin{multline}\label{SREQ1mod}
S_R^{-1}(s,T)B S_L^{-1}(q,T)=\left[\left(S_R^{-1}(s,T)B-B S_L^{-1}(q,T)\right)q-\overline{s}\left(S_R^{-1}(s,T)B-BS_L^{-1}(q,T)\right)\right]\Q_{s}(q)^{-1}
\end{multline}
holds true. Equivalently, it can also be written as
\begin{multline}\label{SREQ2}
S_R^{-1}(s,T)BS_L^{-1}(q,T)=\Q_{q}(s)^{-1}
\left[\left(BS_L^{-1}(q,T) - S_R^{-1}(s,T)B\right)\overline{q}-s\left(BS_L^{-1}(q,T) - S_R^{-1}(s,T)B\right)
 \right],
\end{multline}
where
$
\Q_{s}(q):=q^2 - 2\Re(s)q + |s|^2.
$
\end{theorem}
\begin{proof}
We show that
\begin{multline}\label{SREQ1Var}
S_R^{-1}(s,T)BS_L^{-1}(q,T)\Q_{s}(q)
=\left[\left(S_R^{-1}(s,T)B-B S_L^{-1}(q,T)\right)q-\overline{s}\left(S_R^{-1}(s,T)B-BS_L^{-1}(q,T)\right)\right],
\end{multline}
which is equivalent to \eqref{SREQ1mod}. The modified left $S$-resolvent equation
\eqref{LeftSREQmod}
 implies
\[
BS_L^{-1}(q,T)q=BTS_L^{-1}(q,T)+B.
\]
Applying this identity in the third and fifth equality and using the fact the $TB=BT$:
\[
\begin{split}
  S_R^{-1}&(s,T)BS_L^{-1}(q,T)\Q_{s}(q)
 \\
 &
= S_R^{-1}(s,T)BS_L^{-1}(q,T)(q^2 - 2s_0q + |s|^2)
\\
&
=S_R^{-1}(s,T)\big[BS_L^{-1}(q,T)q\big]q - 2s_0S_R^{-1}(s,T)\big[BS_L^{-1}(q,T)q\big] + |s|^2S_R^{-1}(s,T)BS_L^{-1}(q,T)
 \\
&
=S_R^{-1}(s,T) \big[BTS_L^{-1}(q,T)+B\big]q -2s_0S_R^{-1}(s,T)\big[BTS_L^{-1}(q,T)+B\big]
+|s|^2S_R^{-1}(s,T)BS_L^{-1}(q,T)
\\
&
=S_R^{-1}(s,T)T \big[BS_L^{-1}(q,T)q\big]+S_R^{-1}(s,T)Bq
-2s_0S_R^{-1}(s,T)BTS_L^{-1}(q,T)-2s_0S_R^{-1}(s,T)B
\\
&
+|s|^2S_R^{-1}(s,T)BS_L^{-1}(q,T)
\\
&
=S_R^{-1}(s,T)T \big[BTS_L^{-1}(q,T)+B\big]+S_R^{-1}(s,T)Bq
-2s_0S_R^{-1}(s,T)BTS_L^{-1}(q,T)-2s_0S_R^{-1}(s,T)B
\\
&
+|s|^2S_R^{-1}(s,T)BS_L^{-1}(q,T).
\end{split}
\]
So, the above relation becomes
\[
\begin{split}
S_R^{-1}(s,T)BS_L^{-1}(q,T)\Q_{s}(q)&=
S_R^{-1}(s,T)T BTS_L^{-1}(q,T)+S_R^{-1}(s,T)T B+S_R^{-1}(s,T)Bq
\\
&
-2s_0S_R^{-1}(s,T)BTS_L^{-1}(q,T)-2s_0S_R^{-1}(s,T)B
+|s|^2S_R^{-1}(s,T)BS_L^{-1}(q,T).
\end{split}
\]
Now we replace the term $S_R^{-1}(s,T)BT$ by $sS_R^{-1}(s,T)B-B$, two times,
 in the above formula
using the right $S$-resolvent equation \eqref{RightSREQmod} written as
$$
S_R^{-1}(s,T)BT=sS_R^{-1}(s,T)B-B.
$$
We get the equality
\[
\begin{split}
S_R^{-1}(s,T)B&S_L^{-1}(q,T)\Q_{s}(q)=
\big[S_R^{-1}(s,T)BT\big] TS_L^{-1}(q,T)+\big[S_R^{-1}(s,T)T B\big]+S_R^{-1}(s,T)Bq
\\
&
-2s_0 \big[S_R^{-1}(s,T)BT\big]S_L^{-1}(q,T)-2s_0S_R^{-1}(s,T)B
+|s|^2S_R^{-1}(s,T)BS_L^{-1}(q,T)
\end{split}
\]
and replacing \eqref{RightSREQmod} in the above equation
 one more time in two places indicated in parenthesis
 we get
\[
\begin{split}
S_R^{-1}(s,T)&BS_L^{-1}(q,T)\Q_{s}(q)=
\big[sS_R^{-1}(s,T)B-B\big] TS_L^{-1}(q,T)+\big[sS_R^{-1}(s,T)B-B\big]+S_R^{-1}(s,T)Bq
\\
&
-2s_0 \big[sS_R^{-1}(s,T)B-B\big]S_L^{-1}(q,T)-2s_0S_R^{-1}(s,T)B
+|s|^2S_R^{-1}(s,T)BS_L^{-1}(q,T)
\end{split}
\]
so we get
\[
\begin{split}
S_R^{-1}(s,T)&BS_L^{-1}(q,T)\Q_{s}(q)=
s\big[S_R^{-1}(s,T)BT\big] S_L^{-1}(q,T)-BTS_L^{-1}(q,T)
\\
&
+\big[sS_R^{-1}(s,T)B-B\big]+S_R^{-1}(s,T)Bq
\\
&
-2s_0 \big[sS_R^{-1}(s,T)B-B\big]S_L^{-1}(q,T)-2s_0S_R^{-1}(s,T)B
+|s|^2S_R^{-1}(s,T)BS_L^{-1}(q,T).
\end{split}
\]
Replacing one more time we have
\[
\begin{split}
S_R^{-1}(s,T)&BS_L^{-1}(q,T)\Q_{s}(q)=
s\big[sS_R^{-1}(s,T)B-B\big] S_L^{-1}(q,T)
\\
&
-BTS_L^{-1}(q,T) +\big[sS_R^{-1}(s,T)B-B\big]+S_R^{-1}(s,T)Bq
\\
&
-2s_0 \big[sS_R^{-1}(s,T)B-B\big]S_L^{-1}(q,T)-2s_0S_R^{-1}(s,T)B
\\
&
+|s|^2S_R^{-1}(s,T)BS_L^{-1}(q,T)
\end{split}
\]
so we obtain
\[
\begin{split}
S_R^{-1}(s,T)BS_L^{-1}(q,T)\Q_{s}(q)&=
s^2S_R^{-1}(s,T)BS_L^{-1}(q,T)-s BS_L^{-1}(q,T)
\\
&
-BTS_L^{-1}(q,T) + sS_R^{-1}(s,T)B-B+S_R^{-1}(s,T)Bq
\\
&
-2s_0 sS_R^{-1}(s,T)BS_L^{-1}(q,T)+2s_0 BS_L^{-1}(q,T)-2s_0S_R^{-1}(s,T)B
\\
&
+|s|^2S_R^{-1}(s,T)BS_L^{-1}(q,T).
\end{split}
\]
Now we gather some terms
in order to get
\[
\begin{split}
S_R^{-1}(s,T)BS_L^{-1}(q,T)\Q_{s}(q)&=
(s^2-2s_0s+|s|^2)S_R^{-1}(s,T)BS_L^{-1}(q,T)
\\
&
-s BS_L^{-1}(q,T)
-BTS_L^{-1}(q,T)+ sS_R^{-1}(s,T)B-B+S_R^{-1}(s,T)Bq
\\
&+2s_0 BS_L^{-1}(q,T)-2s_0S_R^{-1}(s,T)B
\end{split}
\]
and
\[
\begin{split}
S_R^{-1}(s,T)BS_L^{-1}(q,T)\Q_{s}(q)&=
(s^2-2s_0s+|s|^2)S_R^{-1}(s,T)BS_L^{-1}(q,T)
\\
&
-s BS_L^{-1}(q,T)
-BTS_L^{-1}(q,T) +sS_R^{-1}(s,T)B-B+S_R^{-1}(s,T)Bq
\\
&+2s_0 BS_L^{-1}(q,T)-2s_0S_R^{-1}(s,T)B.
\end{split}
\]
Since $-s BS_L^{-1}(q,T)+2s_0 BS_L^{-1}(q,T)=\overline{s} BS_L^{-1}(q,T)$
we have
\[
\begin{split}
S_R^{-1}(s,T)BS_L^{-1}(q,T)\Q_{s}(q)&=
(s^2-2s_0s+|s|^2)S_R^{-1}(s,T)BS_L^{-1}(q,T)
\\
&
+\overline{s} BS_L^{-1}(q,T)
-BTS_L^{-1}(q,T) +sS_R^{-1}(s,T)B-B+S_R^{-1}(s,T)Bq
\\
&-2s_0S_R^{-1}(s,T)B.
\end{split}
\]
Recalling the $S$-resolvent equation
$$
BS_L^{-1}(q,T)q=BTS_L^{-1}(q,T)+B
$$
we obtain
\[
\begin{split}
S_R^{-1}(s,T)BS_L^{-1}(q,T)\Q_{s}(q)&=
(s^2-2s_0s+|s|^2)S_R^{-1}(s,T)BS_L^{-1}(q,T)
\\
&
+\overline{s} BS_L^{-1}(q,T)
-BS_L^{-1}(q,T)q +sS_R^{-1}(s,T)B+S_R^{-1}(s,T)Bq
\\
&-2s_0S_R^{-1}(s,T)B
\end{split}
\]
and so
\[
\begin{split}
S_R^{-1}(s,T)BS_L^{-1}(q,T)\Q_{s}(q)&=
(s^2-2s_0s+|s|^2)S_R^{-1}(s,T)BS_L^{-1}(q,T)
\\
&
+\overline{s} BS_L^{-1}(q,T)
-BS_L^{-1}(q,T)q -\overline{s}S_R^{-1}(s,T)B+S_R^{-1}(s,T)Bq.
\end{split}
\]
Finally we obtain
\[
\begin{split}
S_R^{-1}(s,T)BS_L^{-1}(q,T)\Q_{s}(q)&=
(s^2-2s_0s+|s|^2)S_R^{-1}(s,T)BS_L^{-1}(q,T)
\\
&
+\overline{s} BS_L^{-1}(q,T)
-BS_L^{-1}(q,T)q -\overline{s}S_R^{-1}(s,T)B+S_R^{-1}(s,T)Bq
\\
&=
(s^2-2s_0s+|s|^2)S_R^{-1}(s,T)BS_L^{-1}(q,T)
\\
&
+
(S_R^{-1}(s,T)B-BS_L^{-1}(q,T))q   - \overline{s} ( S_R^{-1}(s,T)B- BS_L^{-1}(q,T)
\end{split}
\]
and since $s^2-2s_0s+|s|^2=0$ we get \eqref{SREQ1Var}.
With similar computations we can show that also \eqref{SREQ2} holds.
\end{proof}

\begin{definition}(The modified $S$-resolvent operators)
Let $T\in\mathcal{B}^{\small 0,1}(V_n)$  and $ s \in \rho_S(T)$.
We define the  modified left $S$-resolvent operator of order $\ell\in \mathbb{N}$ as
\begin{equation}\label{MODYPLYSRESOLLFT}
\begin{split}
\overline{T}^\ell S^{-1}_{L}(s,T):=-\overline{T}^\ell\Q_{s}(T)^{-1}(T-\overline{s}\,\id),
\end{split}
\end{equation}
and the  modified right $S$-resolvent operator of order $\ell\in \mathbb{N}$ as
\begin{equation}\label{MODYPLYSRESORIGHT}
\begin{split}
 S^{-1}_{R}(s,T)\overline{T}^\ell:
=-(T-\overline{s}\id) \Q_{s}(T)^{-1}\overline{T}^\ell ,
\end{split}
\end{equation}
where
$\Q_{s}(T):=T^2 - 2\Re(s)T + |s|^2\id.$
\end{definition}

\begin{remark}
{\rm
Observe that when $T$ has commuting components the operator
$B:=\overline{T}^\ell$ commutes with $T$, so in this case, by
Lemma \ref{Scoomution}, we have
$$
 S_L^{-1}(s,T;\overline{T}^\ell)= \overline{T}^\ell S_L^{-1}(s,T),
 $$
and
 $$
  S_R^{-1}(s,T;\overline{T}^\ell) = S_R^{-1}(s,T)\overline{T}^\ell.
$$
}
\end{remark}
We are ready for an alternative, but equivalent, definition of the $PS$-functional calculus.
\begin{definition}[The $PS$-functional calculus (II)]\label{modySCalc}
Let $T\in\mathcal{B}^{\small 0,1}(V_n)$  for
$\mathrm{j}\in\mathbb{S}$  set  $ds_\mathrm{j}=ds (-\mathrm{j})$ and let $M\in \mathbb{N}$.
Then we have the formulations (II) of the $PS$-functional calculus.
We define
\begin{equation}\label{modySCalcL}
F(T) :=\frac{1}{2 \pi}\int_{\partial (U\cap \mathbb{C}_\mathrm{j})}
\sum_{\ell=0}^{M-1}\overline{T}^\ell S^{-1}_{L}(s,T)\, ds_\mathrm{j}\,  f_\ell(s),
 \ \ \ {\rm for\ all} \ \ F(s)=\sum_{\ell=0}^{M-1} \bar s^\ell f_\ell(s)\in\mathcal{PS}^M_L(\sigma_S(T)),
\end{equation}
and
\begin{equation}\label{modySRIGHT}
 F(T):=\frac{1}{2 \pi}\int_{\partial (U\cap \mathbb{C}_\mathrm{j})}
 \sum_{\ell =0}^{M-1} f_\ell(s)\, ds_\mathrm{j}\,S^{-1}_{R}(s,T)\overline{T}^\ell,
  \ \ \ {\rm for\ all} \ \ F(s)=\sum_{\ell=0}^{M-1} f_\ell(s) \bar s^\ell \in\mathcal{PS}^M_R(\sigma_S(T)).
 \end{equation}
\end{definition}

\begin{theorem}
The definition of the $PS$-functional calculus  (II) in (\ref{modySCalcL}) and (\ref{modySRIGHT})
 is well posed because
the integrals
  depend neither on $U$ nor on the imaginary unit $\mathrm{j}\in\mathbb{S}$.
  \end{theorem}
  \begin{proof}
The Definitions (\ref{modySCalcL}) and (\ref{modySRIGHT})  are well posed since
they are based on the $S$-functional calculus and
Lemma \ref{holmodres} that assures that the modified
$S$-resolvent operators preserve the slice monogenicity.
\end{proof}

The linearity of the $PS$-functional calculus (II) is a direct consequence of the definition.
\begin{proposition}
Let $T\in\mathcal{B}^{\small 0,1}(V_n)$ and $M\in \mathbb{N}$.

\begin{enumerate}
\item If  $F,G\in\mathcal{PS}^M_L(\sigma_S(T))$, then we have:

$$(F+G)(T)=F(T)+G(T), \quad \textbf{    }  (F \lambda)(T)=F(T) \lambda, \ \  \text{   for all  } \lambda\in\mathbb{R}_n.$$

\item If  $F,G\in\mathcal{PS}^M_R(\sigma_S(T))$, then we have:

$$(F+G)(T)=F(T)+G(T), \quad \textbf{    }  (\lambda F)(T)=\lambda F(T), \ \ \ \text{   for all  } \lambda\in\mathbb{R}_n.$$
\end{enumerate}

\end{proposition}

\begin{theorem}
Let $T\in\mathcal{B}^{\small 0,1}(V_n)$ and let $M\in \mathbb{N}$. Then:

(I) When $F$ is the series expansion in (\ref{poweserpolyleft}),
 we have
$$
\frac{1}{2 \pi}\int_{\partial (U\cap \mathbb{C}_\mathrm{j})}
\sum_{\ell=0}^{M-1}\overline{T}^\ell S^{-1}_{L}(s,T)\, ds_\mathrm{j}\,  f_\ell(s)
=\sum_{\ell=0}^{M-1} \overline{T}^\ell \sum_{u=0}^\infty T^u A_{k,u}.
$$

(II)
 When $F$ is the series expansion in (\ref{poweserpolyright}),
 we have
$$
\frac{1}{2 \pi}\int_{\partial (U\cap \mathbb{C}_\mathrm{j})}
 \sum_{\ell =0}^{M-1} f_\ell(s)\, ds_\mathrm{j}\,S^{-1}_{R}(s,T)\overline{T}^\ell=
\sum_{\ell=0}^{M-1}  \Big(\sum_{u=0}^\infty A_{k,u} T^u \Big) \overline{T}^\ell.
$$
We recall that $\overline{T}=T_0-\sum_{j=1}^ne_jT_j$.
\end{theorem}
\begin{proof}
It is a consequence of the $S$-functional calculus.
\end{proof}

\begin{theorem}[The intrinsic $PS$-functional calculus (II)]\label{PSfuncintrin}
Let $T\in\mathcal{BC}^{\small 0,1}(V_n)$ and set  $ds_\mathrm{j}=ds (-\mathrm{j})$ and let $F$ be an intrinsic  poly monogenic function  $F\in\mathcal{PN}^M(\sigma_{S}(T))$, for $M\in \mathbb{N}$ with components $f_\ell$. Then the left and the right formulations of the $PS$-functional calculus (II) define the same operator, i.e., we have
\begin{equation}\label{UGUALIPOLYSCalcL}
\begin{split}F(T) =\frac{1}{2 \pi}\int_{\partial (U\cap \mathbb{C}_\mathrm{j})}
\sum_{\ell=0}^{M-1}\overline{T}^\ell S^{-1}_{L}(s,T)\, ds_\mathrm{j}\,  f_\ell(s)
=\frac{1}{2 \pi}\int_{\partial (U\cap \mathbb{C}_\mathrm{j})}
 \sum_{\ell =0}^{M-1} f_\ell(s)\, ds_\mathrm{j}\,S^{-1}_{R}(s,T)\overline{T}^\ell.
\end{split}
\end{equation}
\end{theorem}
\begin{proof}
It is a consequence of the fact that $T\in\mathcal{BC}^{\small 0,1}(V_n)$ and of the $S$-functional calculus for intrinsic functions, i.e.,
\[
\begin{split}
 f_\ell(T)={{1}\over{2\pi }} \int_{\partial (U\cap \mathbb{C}_\mathrm{j})} S_L^{-1} (s,T)\  ds_\mathrm{j} \ f_\ell(s)
={{1}\over{2\pi }} \int_{\partial (U\cap \mathbb{C}_\mathrm{j})} \  f_\ell(s)\ ds_\mathrm{j} \ S_R^{-1} (s,T),
\end{split}
\]
holds.
\end{proof}

\section{Equivalence of the two definitions of the $PS$-functional calculus  and the product rules}\label{equivprodrule}

According to Proposition \ref{Fgprop} we have some possibilities to
 obtain the pointwise product of a poly slice monogenic function and a slice monogenic function
 that preserves the poly slice monogenicity.
According to such proposition we obtain the related product rules.
Using a product that takes out of the class of poly slice monogenic function of a given order, by
Corollary \ref{poinprodpoly}, we obtain a more general product rule.

First we establish the equivalence of the definitions of the $PS$-functional calculus (I) and (II).

\begin{theorem}[Equivalence of the definitions of the $PS$-functional calculus]
Let $T\in\mathcal{B}^{\small 0,1}(V_n)$, $M\in \mathbb{N}$, for
$\mathrm{j}\in\mathbb{S}$,  set  $ds_\mathrm{j}=ds (-\mathrm{j})$, $\overline{\partial}_\mathrm{j}:=\frac{1}{2}(\partial_u+\mathrm{j}\partial_v)$.
Let $\Pi_\ell S_L^{-1}(s,T)$ and $\Pi_\ell S_R^{-1}(s,T)$
 be the $\Pi S$-resolvent operators.
Then we have
\begin{equation}\label{POLYSLEFTdd}
\begin{split}
F(T) &=\frac{1}{2 \pi}\int_{\partial (U\cap \mathbb{C}_\mathrm{j})}\sum_{\ell=0}^{M-1} \Pi_\ell S^{-1}_{L}(s,T)\, ds_\mathrm{j}\, \overline{\partial_\mathrm{j}}^\ell F(s)
\\
&
=
\frac{1}{2 \pi}\int_{\partial (U\cap \mathbb{C}_\mathrm{j})}
\sum_{\ell=0}^{M-1}\overline{T}^\ell S^{-1}_{L}(s,T)\, ds_\mathrm{j}\,  f_\ell(s),
 \ \ \ {\rm for\ all} \ \ F(s)=\sum_{\ell=0}^{M-1}\bar s^\ell f_\ell(s)\in\mathcal{PS}^M_L(\sigma_S(T)),
 \end{split}
\end{equation}
and
\begin{equation}\label{POLYSRIGHTdd}
\begin{split}
 F(T)&=\frac{1}{2 \pi}\int_{\partial (U\cap \mathbb{C}_\mathrm{j})}\sum_{\ell=0}^{M-1} F(s) \overline{\partial_\mathrm{j}}^\ell\, ds_\mathrm{j}\,\Pi_\ell S^{-1}_{R}(s,T)
  \\
  &=
  \frac{1}{2 \pi}\int_{\partial (U\cap \mathbb{C}_\mathrm{j})}
 \sum_{\ell =0}^{M-1} f_\ell(s)\, ds_\mathrm{j}\,S^{-1}_{R}(s,T)\overline{T}^\ell,
  \ \ \ {\rm for\ all} \ \ F(s)=\sum_{\ell=0}^{M-1}f_\ell(s)\bar s^\ell \in\mathcal{PS}^M_R(\sigma_S(T)).
\end{split}
 \end{equation}
\end{theorem}
\begin{proof}
It is a consequence of Runge's theorems for slice monogenic and poly slice monogenic functions
and the continuity properties of the two formulations of the $PS$-functional calculus.
\end{proof}

As a corollary of the above theorem we have the four possible representations of the operator $f(T)$
when $F$ is an intrinsic poly slice monogenic function, i.e., when $F\in\mathcal{PN}^M(\sigma_{S}(T))$.

\begin{theorem}[The formulations of $PS$-functional calculus for intrinsic functions]
Let $T\in\mathcal{BC}^{\small 0,1}(V_n)$, $M\in \mathbb{N}$, for
$\mathrm{j}\in\mathbb{S}$,  set  $ds_\mathrm{j}=ds (-\mathrm{j})$, $\overline{\partial}_\mathrm{j}:=\frac{1}{2}(\partial_u+\mathrm{j}\partial_v)$ and
Let $M\in \mathbb{N}$, for
$\mathrm{j}\in\mathbb{S}$,  set  $ds_\mathrm{j}=ds (-\mathrm{j})$  and  $\overline{\partial}_\mathrm{j}:=\frac{1}{2}(\partial_u+\mathrm{j}\partial_v)$.
Let $\Pi_\ell S_L^{-1}(s,T)$ and $\Pi_\ell S_R^{-1}(s,T)$
 be the $\Pi S$-resolvent operators.
Then, for every $F\in\mathcal{PN}^M(\sigma_{S}(T))$, we have
\begin{equation}\label{POLYSLEFTGreat}
\begin{split}
F(T) &=\frac{1}{2 \pi}\int_{\partial (U\cap \mathbb{C}_\mathrm{j})}\sum_{\ell=0}^{M-1} \Pi_\ell S^{-1}_{L}(s,T)\, ds_\mathrm{j}\, \overline{\partial_\mathrm{j}}^\ell F(s)
\\
&
=
\frac{1}{2 \pi}\int_{\partial (U\cap \mathbb{C}_\mathrm{j})}
\sum_{\ell=0}^{M-1}\overline{T}^\ell S^{-1}_{L}(s,T)\, ds_\mathrm{j}\,  f_\ell(s)
 \\
 &=\frac{1}{2 \pi}\int_{\partial (U\cap \mathbb{C}_\mathrm{j})}\sum_{\ell=0}^{M-1} F(s)\overline{\partial_\mathrm{j}}^\ell\, ds_\mathrm{j}\,\Pi_\ell S^{-1}_{R}(s,T)
  \\
  &=
  \frac{1}{2 \pi}\int_{\partial (U\cap \mathbb{C}_\mathrm{j})}
 \sum_{\ell =0}^{M-1} f_\ell(s)\, ds_\mathrm{j}\,S^{-1}_{R}(s,T)\overline{T}^\ell.
\end{split}
 \end{equation}
\end{theorem}
\begin{proof}
It is a consequence of Runge's theorems for slice monogenic and poly slice monogenic functions
and the continuity properties of the two formulations of the $PS$-functional calculus.
\end{proof}

The following lemma will be crucial in the proof of the product rules.
\begin{lemma}[See \cite{acgs}]\label{HelpProdRule}
Let $B\in \boundOP(V_n)$. For any $q,s\in \mathbb{R}^{n+1}$ with $q\notin[s]$ and let  $f$ be an intrinsic slice monogenic function and let $U$ be a slice Cauchy domain with $\overline{U}\subset \dom(f)$. By setting
$$
\Q_{s}(q)^{-1}:=(q^2-2\Re(s)q+|s|^2)^{-1}
$$
then we have
\[
\frac{1}{2\pi}\int_{\partial(U\cap\cc_\mathrm{j})}f(s)\, ds_\mathrm{j}\, (\overline{s}B-Bq)\Q_{s}(q)^{-1}=Bf(q)
\]
for any $q$ that belongs to $U$ and any $\mathrm{j}\in\mathbb{S}$.

\end{lemma}

\medskip
Using Proposition \ref{Fgprop} we obtain the first set of product rules for the $PS$-functional calculus.
\begin{theorem}[Product rule (first case)]\label{PRDRULEFgprop}
 Let $T\in\mathcal{BC}^{\small 0,1}(V_n)$ and $M\in \mathbb{N}$.

(Ia)
Let $F\in \mathcal{PN}^M(\sigma_S(T))$ and $g\in \mathcal{SM}_L(\sigma_S(T))$.
Then $(Fg)(T)=F(T)g(T)$.

(Ib)
Let $F\in \mathcal{PM}^M_L(\sigma_S(T))$ and $g\in \mathcal{N}(\sigma_S(T))$.
Then  $(gF)(T)=g(T)F(T)$.

(IIa)
Let $F\in \mathcal{PN}^M(\sigma_S(T))$ and $g\in \mathcal{SM}_R(\sigma_S(T))$.
Then  $(gF)(T)=g(T)F(T)$.

(IIb)
 Let
$F\in \mathcal{PM}^M_R(\sigma_S(T))$ and $g\in \mathcal{N}(\sigma_S(T))$.
Then  $(Fg)(T)=F(T)g(T)$.

(III)
Let $F\in \mathcal{PN}^M(\sigma_S(T))$ and $g\in \mathcal{N}(\sigma_S(T))$.
Then we have
$$
(gF)(T)=(Fg)(T)=F(T)g(T)=g(T)F(T).
$$

\end{theorem}
\begin{proof}
We will show point (Ia). The other claims follow in much the same way.
So for $F\in \mathcal{PN}^M(\sigma_S(T))$ and $g\in \mathcal{SM}_L(\sigma_S(T))$
using the $PS$-functional calculus (II) for $F$ and the $S$-functional calculus for
$g$
we have
$$
g(T) :=\frac{1}{2 \pi}\int_{\partial (U\cap \mathbb{C}_\mathrm{j})} S^{-1}_{L}(s,T)\, ds_\mathrm{j}\,  g(s),
 \ \ \ {\rm for\ all} \ \ f\in\mathcal{SM}_L(\sigma_S(T)).
 $$
Let $U_q$ and $U_s$ be bounded slice Cauchy domains that contain $\sigma_{S}(T)$ such that
$\overline{U_q}\subset U_s $ and  $\overline{U_s}\subset \dom(F)\cap\dom(g)$. The subscripts $q$ and $s$ refer to  the respective variable of integration in the following computation. We choose $\mathrm{j}\in\mathbb{S}$ and we set $\Gamma_s := \partial(U_s \cap\cc_\mathrm{j})$ and $\Gamma_q := \partial(U_q\cap\cc_\mathrm{j})$ for neatness.
So for $F\in \mathcal{PN}^M(\sigma_S(T))$ by Theorem \ref{PSfuncintrin} we can represent
$F$ as
$$
F(T)=\frac{1}{2 \pi}\int_{\partial (U\cap \mathbb{C}_\mathrm{j})} \sum_{k=0}^{M-1} f_k(s)\, ds_\mathrm{j}\, S_R^{-1}(s,T) \overline{T}^k.
$$
When we consider the product
$$
F(T)g(T)=\frac{1}{2 \pi}\int_{\Gamma_s} \sum_{k=0}^{M-1} f_k(s)\, ds_\mathrm{j}\, S_R^{-1}(s,T) \overline{T}^k
\frac{1}{2 \pi}\int_{\Gamma_q}S^{-1}_{L}(q,T)\, dq_\mathrm{j}\,  g(q)
$$
so we write
\begin{equation}
\label{sdyktfbx}
F(T)g(T)=\frac{1}{2 \pi}\int_{\Gamma_s} \sum_{k=0}^{M-1} f_k(s)\, ds_\mathrm{j}\,
\frac{1}{2 \pi}\int_{\Gamma_q} S_R^{-1}(s,T) \overline{T}^kS^{-1}_{L}(q,T)\, dq_\mathrm{j}\,  g(q).
\end{equation}
Now we use the modified $S$-resolvent equation with $B=\overline{T}^k$ that commute with $T$, since  the paravector operator $T$ has commuting components we have the modified $S$-resolvent equation
\begin{equation}
\label{diwajk}
S_R^{-1}(s,T)\overline{T}^k S_L^{-1}(q,T)=\left[\left(S_R^{-1}(s,T)\overline{T}^k-\overline{T}^k S_L^{-1}(q,T)\right)q-\overline{s}\left(S_R^{-1}(s,T)\overline{T}^k-\overline{T}^kS_L^{-1}(q,T)\right)
\right]\Q_{s}(q)^{-1}
\end{equation}
we replace it (\ref{diwajk}) in (\ref{sdyktfbx}) to get
$$
F(T)g(T)=\frac{1}{2 \pi}\int_{\Gamma_s} \sum_{k=0}^{M-1} f_k(s)\, ds_\mathrm{j}\,
\frac{1}{2 \pi}\int_{\Gamma_q} \Big[\Big(S_R^{-1}(s,T)\overline{T}^k-\overline{T}^k S_L^{-1}(q,T)\Big)q
$$
$$
-\overline{s}\Big(S_R^{-1}(s,T)\overline{T}^k-\overline{T}^kS_L^{-1}(q,T)\Big)
\Big]\Q_{s}(q)^{-1}\, dq_\mathrm{j}\,  g(q)
$$
but also
\[
\begin{split}
F(T)g(T)&=\frac{1}{2 \pi}\int_{\Gamma_s} \sum_{k=0}^{M-1} f_k(s)\, ds_\mathrm{j}\,
\frac{1}{2 \pi}\int_{\Gamma_q}
S_R^{-1}(s,T)\overline{T}^k q\Q_{s}(q)^{-1}\, dq_\mathrm{j}\,  g(q)
\\
&-\frac{1}{2 \pi}\int_{\Gamma_s} \sum_{k=0}^{M-1} f_k(s)\, ds_\mathrm{j}\,
\frac{1}{2 \pi}\int_{\Gamma_q}
\overline{T}^k S_L^{-1}(q,T)q
\Q_{s}(q)^{-1}\, dq_\mathrm{j}\,  g(q)
\\
&-\frac{1}{2 \pi}\int_{\Gamma_s} \sum_{k=0}^{M-1} f_k(s)\, ds_\mathrm{j}\,
\frac{1}{2 \pi}\int_{\Gamma_q}
\overline{s}S_R^{-1}(s,T)\overline{T}^k
\Q_{s}(q)^{-1}\, dq_\mathrm{j}\,  g(q)
\\
&+\frac{1}{2 \pi}\int_{\Gamma_s} \sum_{k=0}^{M-1} f_k(s)\, ds_\mathrm{j}\,
\frac{1}{2 \pi}\int_{\Gamma_q}
\overline{s}\overline{T}^kS_L^{-1}(q,T)
\Q_{s}(q)^{-1}\, dq_\mathrm{j}\,  g(q).
\end{split}
\]
Now the integrals with the right $S$-resolvent operator are zero,
in fact we can write them as
\[
\begin{split}
&\frac{1}{2 \pi}\int_{\Gamma_s} \sum_{k=0}^{M-1} f_k(s)\, ds_\mathrm{j}\,
\frac{1}{2 \pi}\int_{\Gamma_q}
S_R^{-1}(s,T)\overline{T}^k q\Q_{s}(q)^{-1}\, dq_\mathrm{j}\,  g(q)
 \\
 =
 &
 \frac{1}{2 \pi}\int_{\Gamma_s} \sum_{k=0}^{M-1} f_k(s)\, ds_\mathrm{j}\,
S_R^{-1}(s,T)\overline{T}^k \frac{1}{2 \pi}\int_{\Gamma_q}\Big[q\Q_{s}(q)^{-1}\, dq_\mathrm{j}\,  g(q)\Big]=0
\end{split}
\]
by the Cauchy theorem because $q\Q_{s}(q)^{-1}$ and $g(q)$ are slice monogenic. Similarly it is
$$
-\frac{1}{2 \pi}\int_{\Gamma_s} \sum_{k=0}^{M-1} f_k(s)\, ds_\mathrm{j}\,
\frac{1}{2 \pi}\int_{\Gamma_q}
\overline{s}S_R^{-1}(s,T)\overline{T}^k
\Q_{s}(q)^{-1}\, dq_\mathrm{j}\,  g(q)=0
$$
so
we remain with
\[
\begin{split}
F(T)g(T)&=-\frac{1}{2 \pi}\int_{\Gamma_s} \sum_{k=0}^{M-1} f_k(s)\, ds_\mathrm{j}\,
\frac{1}{2 \pi}\int_{\Gamma_q}
\overline{T}^k S_L^{-1}(q,T)q
\Q_{s}(q)^{-1}\, dq_\mathrm{j}\,  g(q)
\\
&+\frac{1}{2 \pi}\int_{\Gamma_s} \sum_{k=0}^{M-1} f_k(s)\, ds_\mathrm{j}\,
\frac{1}{2 \pi}\int_{\Gamma_q}
\overline{s}\overline{T}^kS_L^{-1}(q,T)
\Q_{s}(q)^{-1}\, dq_\mathrm{j}\,  g(q)
\end{split}
\]
and also we obtain
\[
\begin{split}
F(T)g(T)&=\frac{1}{2 \pi}\int_{\Gamma_s} \sum_{k=0}^{M-1} f_k(s)\, ds_\mathrm{j}\,
\frac{1}{2 \pi}\int_{\Gamma_q}
\Big(\overline{s}\overline{T}^kS_L^{-1}(q,T)
-\overline{T}^k S_L^{-1}(q,T)q\Big)
\Q_{s}(q)^{-1}\, dq_\mathrm{j}\,  g(q) .
\end{split}
\]
Using Fubini's theorem we get
\[
\begin{split}
F(T)g(T)&=\frac{1}{2 \pi}\int_{\Gamma_q}
\Big[\frac{1}{2 \pi}\int_{\Gamma_s} \sum_{k=0}^{M-1} f_k(s)\, ds_\mathrm{j}\,
\Big(\overline{s}\overline{T}^kS_L^{-1}(q,T)
-\overline{T}^k S_L^{-1}(q,T)q\Big)
\Q_{s}(q)^{-1}\Big]\, dq_\mathrm{j}\,  g(q)
\end{split}
\]
and by Lemma \ref{HelpProdRule}, setting $B:=\overline{T}^k S_L^{-1}(q,T)$, we obtain
\[
\begin{split}
F(T)g(T)&=\frac{1}{2 \pi}\int_{\Gamma_q}
 \sum_{k=0}^{M-1} \overline{T}^kS_L^{-1}(q,T) f_k(q)
\, dq_\mathrm{j}\,  g(q)
\end{split}
\]
and finally
\[
\begin{split}
F(T)g(T)&=\frac{1}{2 \pi}\int_{\Gamma_q}
 \sum_{k=0}^{M-1} \overline{T}^kS_L^{-1}(q,T)
\, dq_\mathrm{j}\,  f_k(q)g(q)
\\
&
=\sum_{k=0}^{M-1} \overline{T}^k (f_kg)(T)
=(Fg)(T),
\end{split}
\]
and this concludes the proof.
\end{proof}

Based on the product of poly slice monogenic functions in  Corollary \ref{poinprodpoly} we  prove the following product rule. This product rule requires that the paravector operator $T$ has commuting components.

\begin{theorem}[Product rule (second case)]\label{PRULEpoly}
Let $T\in\mathcal{BC}^{\small 0,1}(V_n)$ and $M, N\in \mathbb{N}$.

(I)
Let $F\in \mathcal{N}^N_L(\sigma_S(T))$ and $G\in \mathcal{PS}^M_L(\sigma_S(T))$, then we have

\begin{equation}\label{21RULE}
(FG)(T)=F(T)G(T).
\end{equation}

(II)
Let $F\in \mathcal{PS}^N_R(\sigma_S(T))$ and $G\in \mathcal{N}^M_R(\sigma_S(T))$, then we have

\begin{equation}\label{RULE55}
(FG)(T)=F(T)G(T).
\end{equation}
\end{theorem}
\begin{proof}
 We reason as in Theorem \ref{PRDRULEFgprop} taking the same contours of integration and pointing out just the differences.
We show point (I). The other point (II) follows in much the same way.
Let $U_q$ and $U_s$ be bounded slice Cauchy domains that contain $\sigma_{S}(T)$ such that
$\overline{U_q}\subset U_s $ and  $\overline{U_s}\subset \dom(F)\cap\dom(G)$. The subscripts $q$ and $s$ refer to  the respective variable of integration in the following computation. We choose $\mathrm{j}\in\mathbb{S}$ and we set $\Gamma_s := \partial(U_s \cap\cc_\mathrm{j})$ and $\Gamma_q := \partial(U_q\cap\cc_\mathrm{j})$ for neatness.
So by Theorem \ref{PSfuncintrin} we can represent
$F$ as
$$
F(T)=\frac{1}{2 \pi}\int_{\Gamma_s} \sum_{k=0}^{N-1} f_k(s)\, ds_\mathrm{j}\, S_R^{-1}(s,T) \overline{T}^k, \ \ \ {\rm for\ all} \ \ F\in \mathcal{PN}^N(\sigma_S(T))
$$
and
$$
G(T)=
\frac{1}{2 \pi}\int_{\Gamma_q}
\sum_{\ell=0}^{M-1}\overline{T}^\ell S^{-1}_{L}(q,T)\, dq_\mathrm{j}\,  g_\ell(q),
 \ \ \ {\rm for\ all} \ \ F\in\mathcal{PS}^M_L(\sigma_S(T))
 $$
so when we consider the product
$$
F(T)G(T)=\frac{1}{2 \pi}\int_{\Gamma_s} \sum_{k=0}^{N-1} f_k(s)\, ds_\mathrm{j}\, S_R^{-1}(s,T) \overline{T}^k
\frac{1}{2 \pi}\int_{\Gamma_q}
\sum_{\ell=0}^{M-1}\overline{T}^\ell S^{-1}_{L}(q,T)\, dq_\mathrm{j}\,  g_\ell(q)
 $$
 and also.
 $$
F(T)G(T)=\frac{1}{2 \pi}\int_{\Gamma_s} \sum_{k=0}^{N-1} f_k(s)\, ds_\mathrm{j}\,
\frac{1}{2 \pi}\int_{\Gamma_q}
\sum_{\ell=0}^{M-1}S_R^{-1}(s,T) \overline{T}^{k+\ell} S^{-1}_{L}(q,T)\, dq_\mathrm{j}\,  g_\ell(q).
 $$
 We use the modified $S$-resolvent equation to get
 \[
\begin{split}
F(T)G(T)&=\frac{1}{2 \pi}\int_{\Gamma_s} \sum_{k=0}^{N-1} f_k(s)\, ds_\mathrm{j}\,
\frac{1}{2 \pi}\int_{\Gamma_q}
\sum_{\ell=0}^{M-1}S_R^{-1}(s,T)\overline{T}^{k+\ell} q\Q_{s}(q)^{-1}\, dq_\mathrm{j}\,  g_\ell(q)
\\
&
-\frac{1}{2 \pi}\int_{\Gamma_s} \sum_{k=0}^{N-1} f_k(s)\, ds_\mathrm{j}\,
\frac{1}{2 \pi}\int_{\Gamma_q}
\sum_{\ell=0}^{M-1}
\overline{T}^{k+\ell} S_L^{-1}(q,T)q
\Q_{s}(q)^{-1}\, dq_\mathrm{j}\,  g_\ell(q)
\\
&-\frac{1}{2 \pi}\int_{\Gamma_s} \sum_{k=0}^{N-1} f_k(s)\, ds_\mathrm{j}\,
\frac{1}{2 \pi}\int_{\Gamma_q}
\sum_{\ell=0}^{M-1}
\overline{s}S_R^{-1}(s,T)\overline{T}^{k+\ell}
\Q_{s}(q)^{-1}\, dq_\mathrm{j}\,  g_\ell(q)
\\
&+\frac{1}{2 \pi}\int_{\Gamma_s} \sum_{k=0}^{N-1} f_k(s)\, ds_\mathrm{j}\,
\frac{1}{2 \pi}\int_{\Gamma_q}
\sum_{\ell=0}^{M-1}
\overline{s}\overline{T}^{k+\ell}S_L^{-1}(q,T)
\Q_{s}(q)^{-1}\, dq_\mathrm{j}\,  g_\ell(q).
\end{split}
\]
Also here the two integrals that contain the right $S$-resolvent operators are zero
so we obtain
 \[
\begin{split}
F(T)G(T)=
&
-\frac{1}{2 \pi}\int_{\Gamma_s} \sum_{k=0}^{N-1} f_k(s)\, ds_\mathrm{j}\,
\frac{1}{2 \pi}\int_{\Gamma_q}
\sum_{\ell=0}^{M-1}
\overline{T}^{k+\ell} S_L^{-1}(q,T)q
\Q_{s}(q)^{-1}\, dq_\mathrm{j}\,  g_\ell(q)
\\
&+\frac{1}{2 \pi}\int_{\Gamma_s} \sum_{k=0}^{N-1} f_k(s)\, ds_\mathrm{j}\,
\frac{1}{2 \pi}\int_{\Gamma_q}
\sum_{\ell=0}^{M-1}
\overline{s}\overline{T}^{k+\ell}S_L^{-1}(q,T)
\Q_{s}(q)^{-1}\, dq_\mathrm{j}\,  g_\ell(q).
\end{split}
\]
Using Fubini's theorem we finally have
\[
\begin{split}
&F(T)G(T)=
\\
&=\frac{1}{2 \pi}\int_{\Gamma_q}
\Big[\frac{1}{2 \pi}\int_{\Gamma_s} \sum_{k=0}^{N-1} \sum_{\ell=0}^{M-1}f_k(s)\, ds_\mathrm{j}\,
\Big(\overline{s}\overline{T}^{k+\ell}S_L^{-1}(q,T)
-\overline{T}^{k+\ell} S_L^{-1}(q,T)q\Big)
\Q_{s}(q)^{-1}\Big]\, dq_\mathrm{j}\,  g_\ell(q)
\end{split}
\]
setting $B:=\overline{T}^{k+\ell} S_L^{-1}(q,T)$ in  Lemma \ref{HelpProdRule} we get
\[
\begin{split}
F(T)G(T)&=\frac{1}{2 \pi}\int_{\Gamma_q}
 \sum_{k=0}^{N-1} \sum_{\ell=0}^{M-1}\overline{T}^{k+\ell}S_L^{-1}(q,T) f_k(q)
\, dq_\mathrm{j}\,  g_\ell(q)
\\
&
=\frac{1}{2 \pi}\int_{\Gamma_q}
 \sum_{k=0}^{N-1} \sum_{\ell=0}^{M-1}\overline{T}^{k+\ell}S_L^{-1}(q,T)
\, dq_\mathrm{j}\,  (f_kg_\ell)(q)
\\
&
=
 \sum_{k=0}^{N-1} \sum_{\ell=0}^{M-1}\overline{T}^{k+\ell}  (f_kg_\ell)(T)
=(FG)(T),
\end{split}
\]
and this concludes the proof.
\end{proof}


\end{document}